\documentclass[11pt]{amsart}

\newif\iffinalrun
\finalruntrue

\newif\ifkuvio

\iffinalrun
  \kuviotrue
\fi

\usepackage{amssymb,eucal}
\iffinalrun
\else
  \usepackage[notref,notcite]{showkeys}
\fi
\usepackage{epsfig}

\newcommand{\color}[2][1]{}

\ifkuvio
  \usepackage[forcekdg,arrsy]{kuvio}
\fi

\usepackage{hyperref}
\DeclareMathAlphabet{\mathpzc}{OT1}{pzc}{m}{it}

\newcommand{\noi}{\noindent}

\newcommand{\bs}{\backslash}
\def\varddots{\mathinner{\raise7pt\vbox{\kern3pt\hbox{.}}\mkern1mu\smash{\raise4pt\hbox{.}}
\mkern1mu\smash{\raise1pt\hbox{.}}}}

\newcommand{\p}{\mathpzc p}

\newcommand{\B}{\mathcal B}
\newcommand{\oB}{{\overline\B}}
\newcommand{\C}{\mathbb C}

\newcommand{\F}{\mathbb F}

\newcommand{\fpb}{\overline \F_p}

\newcommand{\G}{\mathbb G}
\newcommand{\fG}{\mathfrak G}

\newcommand{\HH}{\mathcal H}

\newcommand{\hg}{\HH_G(V)}
\newcommand{\hgd}{\HH_G(V^*)}

\newcommand{\htd}{\HH_T((V^*)_{\o U(k)})}

\newcommand{\hm}{\HH_M(\vonk)}

\newcommand{\I}{\mathcal I}
\renewcommand{\k}{{\bar k}}

\newcommand{\mm}{\mathfrak m}

\renewcommand{\O}{{\mathcal O}}

\renewcommand{\P}{\mathcal P}
\newcommand{\oP}{{\overline\P}}
\newcommand{\Q}{\mathcal Q}
\newcommand{\oQ}{{\overline\Q}}
\newcommand{\QQ}{\mathbb Q}

\newcommand{\qp}{\QQ_p}

\newcommand{\R}{\mathbb R}
\renewcommand{\SS}{\mathcal S}

\newcommand{\sg}{\SS_G}
\newcommand{\sm}{\SS_M}
\newcommand{\smg}{\SS^M_G}
\newcommand{\smgop}{\SS^{M,\mathrm{op}}_G}
\newcommand{\psg}{{}'\SS_G}
\newcommand{\psm}{{}'\SS_M}
\newcommand{\psmg}{{}'\SS^M_G}

\newcommand{\fT}{\mathfrak T}

\newcommand{\fU}{\mathfrak U}

\newcommand{\vp}{\varphi}

\newcommand{\vd}{V^*}
\newcommand{\vnk}{V^{N(k)}}
\newcommand{\vonk}{V_{\o N(k)}}
\newcommand{\vuk}{V^{U(k)}}
\newcommand{\vouk}{V_{\o U(k)}}

\newcommand{\Z}{\mathbb Z}

\renewcommand{\o}[1]{\overline{#1}}
\newcommand{\wt}[1]{\widetilde{#1}}

\newcommand{\INTO}{\hookrightarrow}
\newcommand{\onto}{\twoheadrightarrow}
\newcommand{\congto}{\xrightarrow{\,\sim\,}}
\newcommand{\tocong}{\xleftarrow{\,\sim\,}}

\newcommand{\s}{^\times}

\newcommand{\dual}{^\vee}

\newcommand{\kalg}{_{\text{$\k$-alg}}}

\DeclareMathOperator{\End}{End}
\DeclareMathOperator{\Hom}{Hom}

\DeclareMathOperator{\GL}{GL}

\DeclareMathOperator{\PGL}{PGL}
\DeclareMathOperator{\GSp}{GSp}
\DeclareMathOperator{\Special}{Sp}
\newcommand{\Sp}[1]{\Special_{#1}}
\newcommand{\oSp}[1]{\overline{\Special}_{#1}}

\DeclareMathOperator{\diag}{diag}

\DeclareMathOperator{\Gal}{Gal}
\DeclareMathOperator{\Stab}{Stab}
\DeclareMathOperator{\tr}{tr}

\DeclareMathOperator{\ord}{ord}

\DeclareMathOperator{\Ind}{Ind}
\DeclareMathOperator{\ind}{c-Ind}

\DeclareMathOperator{\soc}{soc}
\DeclareMathOperator{\red}{{red}}
\DeclareMathOperator{\ch}{ch}
\DeclareMathOperator{\im}{im}

\DeclareMathOperator{\Sym}{Sym}
\DeclareMathOperator{\supp}{supp}

\DeclareMathOperator{\id}{id}
\DeclareMathOperator{\Ord}{Ord}
\DeclareMathOperator{\JH}{JH}

\newcommand{\opInd}{\Ind_{\o P}^G}
\newcommand{\oqInd}{\Ind_{\o Q}^G}
\newcommand{\obInd}{\Ind_{\o B}^G}
\newcommand{\kind}{\ind_K^G}

\newcommand{\lind}{\ind_{L(\O)}^L}
\newcommand{\mind}{\ind_{M(\O)}^M}
\newcommand{\mzind}{\ind_{M(\O)Z_M}^M}
\newcommand{\pind}{\ind_\P^G}
\newcommand{\opind}{\ind_\oP^G}

\iffinalrun
  \newcommand{\need}[1]{}
  \newcommand{\mar}[1]{}
\else
  \newcommand{\need}[1]{{\tiny *** #1}}
  \newcommand{\mar}[1]{\marginpar{\tiny #1}}
\fi

\newcommand{\rr}{\rrbracket}

\renewcommand{\)}{\textup{)}}

\theoremstyle{plain} 
\newtheorem{lm}[equation]{Lemma}
\newtheorem{sublm}[equation]{Sublemma}
\newtheorem{prop}[equation]{Proposition}
\newtheorem{thm}[equation]{Theorem}
\newtheorem{coroll}[equation]{Corollary}

\theoremstyle{definition}
\newtheorem{df}[equation]{Definition}
\newtheorem{rk}[equation]{Remark}
\newtheorem{qu}[equation]{Question}
\newtheorem{ex}[equation]{Example}

\numberwithin{equation}{section}
\numberwithin{figure}{section}

\def\RCS$#1: #2 ${\expandafter\def\csname RCS#1\endcsname{#2}}
\RCS$Revision: 28 $
\RCS$Date: 2011-02-19 15:56:57 -0500 (Sat, 19 Feb 2011) $
\newcommand{\versioninfo}{Version \RCSRevision; Last commit \RCSDate}

\begin{document}

\title[Irreducible mod $p$ representations of a $p$-adic $\GL_n$]
{The classification of irreducible admissible mod $p$ representations of a $p$-adic $\GL_n$}
\author{Florian Herzig}
\address{Department of Mathematics, 2033 Sheridan Road, Evanston IL 60208-2730, USA}
\email{herzig@math.northwestern.edu}
\thanks{Partially supported by NSF grant DMS-0902044}
\date{\today}
\maketitle

\iffinalrun
\else
\versioninfo
\fi

\begin{abstract}
Let $F$ be a finite extension of~$\qp$. Using the mod~$p$ Satake transform, we define what it means for an
irreducible admissible smooth representation of an $F$-split $p$-adic reductive group over~$\fpb$ to be supersingular.
We then give the classification of irreducible admissible smooth $\GL_n(F)$-representations over~$\fpb$ in 
terms of supersingular representations. As a consequence we deduce that \emph{supersingular} is the same
as \emph{supercuspidal}. These results generalise the work of Barthel--Livn\'e for $n = 2$. For general
split reductive groups we obtain similar results under stronger hypotheses.
\end{abstract}

\section{Introduction}
\label{sec:introduction}

Let~$F$ be a finite extension of~$\qp$ with ring of integers~$\O$ and residue field~$k$.
The hypothetical mod $p$ Langlands correspondence is expected to associate to an $n$-dimensional mod $p$ Galois
representation $\rho : \Gal(\o F/F) \to \GL_n(\fpb)$ an admissible smooth representation $\pi(\rho)$ of $\GL_n(F)$ over
$\fpb$ (or maybe an equivalence class of such representations). So far this is understood in the case $n = 2$ and $F = \qp$
(see \cite{bib:Breuil2}, \cite{bib:Colmez}), with some recent progress on the case $n = 2$ and $F/\qp$ unramified (see \cite{bib:BP}). While mod
$p$ Galois representations are relatively easy to understand, the theory of mod $p$ smooth representations of $p$-adic
reductive groups is only at its beginnings. It takes its origin with the fundamental work of Barthel--Livn\'e
\cite{bib:BL-unram}, \cite{bib:BL-general} that classifies irreducible representations of $\GL_2(F)$ over $\fpb$ that have a
central character into four classes: (i) irreducible principal series, (ii) one-dimensional representations, (iii) twists of the Steinberg representation,
and (iv) ``supersingular'' representations. While the first three classes are explicit, the
classification of supersingular representations has only been completed so far for $F = \qp$ \cite{bib:Breuil}.  When $F \ne
\qp$ it turned out to be much more complicated \cite{bib:BP}, \cite{bib:Hu}.

In this paper we generalise the work of Barthel--Livn\'e to $\GL_n$, giving the classification of irreducible admissible representations
of $\GL_n(F)$ over $\fpb$ in terms of supersingular representations. (From now on ``admissible'' is short for ``admissible smooth''.)
We first define what it means for an irreducible
admissible representation to be supersingular (generalising the definition for $\GL_2$) by means of the mod $p$ Satake isomorphism
\cite{bib:satake}. To state our main theorem we need to recall the definition of the generalised Steinberg representations.
Let $P$ denote any standard parabolic subgroup, that is a parabolic subgroup containing the Borel subgroup $B$ of
upper triangular matrices. Let $\o P$ denote the opposite parabolic and let
\begin{equation*}
  \Sp P = \frac{\Ind_{\o P(F)}^{\GL_n(F)} 1}{\sum_{Q \supsetneq P}\Ind_{\o Q(F)}^{\GL_n(F)} 1}.
\end{equation*}
It is known to be irreducible and admissible \cite{bib:grosse-kloenne}.

\begin{thm}\label{thm:class-intro}
  The irreducible admissible $\GL_n(F)$-representation are given by
  $\Ind_{\o P(F)}^{\GL_n(F)} (\sigma_1 \otimes \cdots \otimes \sigma_r)$, where
  \begin{enumerate}
  \item $P$ is a standard parabolic with Levi $\prod_{i=1}^r \GL_{n_i}$;
  \item $\sigma_i$ is an irreducible admissible $\GL_{n_i}(F)$-representation such that either
    \begin{itemize}
    \item $\sigma_i$ is supersingular and $n_i > 1$, or
    \item $\sigma_i \cong \Sp{Q_i} \otimes (\eta_i \circ \det)$ for some smooth character $\eta_i : F\s \to \fpb\s$ and some
      standard parabolic $Q_i \subset \GL_{n_i}$;
    \end{itemize}
  \item $\eta_i \ne \eta_{i+1}$ whenever both $\sigma_i$ and $\sigma_{i+1}$ fall into the second case.
  \end{enumerate}
  Moreover, $P$ is uniquely determined and each $\sigma_i$ is unique up to isomorphism.
\end{thm}

Note that Ollivier \cite{bib:ollivier} proved the irreducibility of this representation when it is a principal series, i.e.,
when $n_1 = \cdots = n_r = 1$.  In this case $\sigma_i = \eta_i$ for all $i$ and consecutive characters are
distinct. When $n = 2$ the result recovers the classification of Barthel--Livn\'e: when $P = B$ one
obtains the irreducible principal series; when $P = G$ one obtains twists of generalised
Steinberg representations (either trivial or Steinberg in this case) and supersingular representations.

\begin{coroll}\label{cor:class-intro}
  Suppose that $\pi$ is an irreducible admissible $\GL_n(F)$-representation.
  Suppose that $Q = LN'$ is a standard parabolic and that $\tau$ is an irreducible admissible $L(F)$-representation.
  \begin{enumerate}
  \item $\Ind_{\o Q(F)}^{\GL_n(F)} \tau$ is of finite length, and all constituents occur with multiplicity one.
  \item $\pi$ is supersingular if and only if $\pi$ is supercuspidal.
  \end{enumerate}
\end{coroll}

Note that the constituents in part (i) can be described explicitly.  By a \emph{supercuspidal} representation in part (ii) we
mean an irreducible admissible representation that does not occur among the constituents of any parabolic induction in part
(i) of the corollary with $Q \ne \GL_n$.

These results follow from the work of Barthel--Livn\'e when $n = 2$. Part (i) was also known in the case of the
trivial principal series $\Ind_{\o B(F)}^{\GL_n(F)} 1$ by \cite{bib:grosse-kloenne} (the constituents are the $\Sp P$) and for all principal
series when $n = 3$ \cite{bib:vigneras}.

It is not hard to determine the submodule structure of the representations in part (i) using the methods
of this paper. See Section~\ref{sec:submodule-structure}.

\subsection{Comparison with classical results}
\label{sec:classical-results}

It is interesting to compare the results with the results of Bernstein and Zelevinsky over the complex numbers \cite{bib:BZ},
\cite{bib:BZ2}.  The most striking difference is the complete lack of intertwining operators in our context, by the uniqueness part of
Theorem~\ref{thm:class-intro}.  Thus for example the principal series $\obInd(\chi_1 \otimes \cdots \otimes \chi_n)$ and
$\obInd(\chi_1' \otimes \cdots \otimes \chi_n')$ do not share any common constituents unless $\chi_i = \chi_i'$ for all $i$.
Note that there is no obvious way to produce intertwining operators over $\fpb$ due to the lack of an $\fpb$-valued Haar
measure.  One would however expect non-trivial extensions between the irreducible representations and the knowledge of these
will play some role in the understanding of a mod~$p$ Langlands correspondence. (For $\GL_2(\qp)$ extensions between
irreducibles were computed in \cite{bib:Colmez}, \cite{bib:emerton-ordinary2}, \cite{bib:BP}, \cite{bib:pask}.)

Another difference, maybe related to the lack of intertwining operators, is that parabolic inductions tend to be irreducible.
In some sense the only reducibilities arise for ``obvious'' reasons, namely when $\eta_i = \eta_{i+1}$ for some $i$ in
Theorem~\ref{thm:class-intro}.  And in those cases constituents occur with multiplicity one, unlike over the complex numbers
(in fact Zelevinsky formulated an analogue of the Kazhdan--Lusztig conjecture to predict the multiplicities \cite{bib:zele-kl}).

\subsection{Methods used}
\label{sec:methods}

We now want to explain in more detail what goes into the proofs. Most of our methods work for more general groups, so let $G$
be a split connected reductive group over $F$. It is known that $G$ extends to a connected reductive group over the ring of
integers $\O$. Fix such an integral model $G_{/\O}$, and fix a maximal split torus $T_{/\O}$ and a Borel subgroup $B_{/\O}$
containing $T$.  Then $K := G(\O)$ is a hyperspecial maximal compact subgroup of $G(F)$ \cite[3.8.1]{bib:Tits}. We usually
write $G$ for $G(F)$, etc. This should not lead to any confusion.

\subsubsection{Weights, Hecke eigenvalues, and supersingular representations}
\label{sec:weights-and-hecke}

Let $\pi$ be any admissible smooth $G$-representation over $\fpb$. It is easy to see that $\pi|_K$ contains an irreducible
$K$-subrepresentation $V$.  Since the kernel of $K = G(\O) \to G(k)$ is a pro-$p$ group and the coefficient field has
characteristic~$p$, it follows that $V$ factors through a representation of the finite group $G(k)$ over $\fpb$.  Such a
representation is called a \emph{$K$-weight}, or simply \emph{weight}. The multiplicity space $\Hom_K(V,\pi)$ is
finite-dimensional (as $\pi$ is admissible). Let $\kind V$ denote the representation compactly induced from $V$. By Frobenius
reciprocity we have $\Hom_K(V,\pi) \cong \Hom_G(\kind V, \pi)$, so that the \emph{Hecke algebra} $\hg := \End_G(\kind V)$
naturally acts on it.  Let $X_*(T)_-$ denote the subset of antidominant coweights in $\Hom(\G_m,T)$. In~\cite{bib:satake} we
established an analogue of the Satake isomorphism, showing that $\hg \cong \fpb[X_*(T)_-]$.  (The map depends on the choice
of positive roots and on the choice of a uniformiser.) In particular, $\hg$ is commutative. Thus the $\hg$-module
$\Hom_K(V,\pi)$ is a direct sum of generalised eigenspaces, and there exists an eigenvector in each of them.  Note that an
algebra homomorphism $\chi : \hg \to \fpb$ occurs as set of Hecke eigenvalues of $\hg$ on $\Hom_K(V,\pi)$ if and only if
there is a non-zero $G$-linear map
\begin{equation}\label{eq:11}
  \kind V \otimes_{\hg,\chi} \fpb \to \pi.
\end{equation}

Suppose that $\pi$ is moreover irreducible. For each set of Hecke eigenvalues $\chi$ occurring on $\Hom_K(V,\pi)$ we get a
monoid homomorphism $\chi' : X_*(T)_- \to \fpb$ via the Satake isomorphism. Note that in our setup above we fixed a reductive
integral structure $G_{/\O}$, giving rise to $K$. We say that $\pi$ is \emph{supersingular} if, for any reductive
integral structure $G_{/\O}$, we have that $\chi'$ is zero on all non-invertible elements of $X_*(T)_-$ for all $\chi$ and
all $V$. (When $G = \GL_n$, all hyperspecial maximal subgroups are conjugate, so that we can work with a fixed $G_{/\O}$
and $K$. Moreover, as a corollary to Theorem~\ref{thm:class-intro} we will see that $\chi'$ is independent of $\chi$ and $V$
above, so that the condition only needs to be checked once.)

We now discuss two of the key tools that go into the proofs of the main results.

\subsubsection{Comparison of compact and parabolic inductions}
\label{sec:comp-inductions}

Suppose we are given a weight $V$ and an algebra homomorphism $\chi : \hg \to \fpb$.
the support of $\chi' : X_*(T)_- \to \fpb$ determines a facet of the antidominant Weyl chamber which in turn
determines a standard parabolic $P_\chi$. Let $P = MN$ be a standard parabolic containing $P_\chi$.
If $V$ is sufficiently ``regular'' (depending on $P$), a condition that is satisfied by most weights, we show that
there exists a natural isomorphism
\begin{equation}\label{eq:24}
  \kind V\otimes_{\hg,\chi} \fpb \congto \opInd (\pi_M),
\end{equation}
for a certain smooth $M$-representation $\pi_M$ defined in terms of $V$ and $\chi$. Note that $\pi_M$ is not admissible in
general. It is a representation of the same form as the left-hand side. (See Theorem~\ref{thm:cptind-parabind} for the
precise statement.) When $P = G$, the isomorphism~(\ref{eq:24}) is the identity map.

When $G = \GL_2$ and $P = B$ the isomorphism~(\ref{eq:24}) was established by Barthel--Livn\'e using a calculation with the
Bruhat--Tits tree. In that case $\pi_M$ is a character. In our proof we compare compact and parabolic inductions via a
parahoric induction, building on the ideas of Schneider--Stuhler \cite{bib:schneider-stuhler} and Vign\'eras
\cite{bib:vigneras2}.

\subsubsection{Changing the weight}
\label{sec:changing-weight}

In certain situations we construct isomorphisms of the form
\begin{equation*}
  \kind V_1\otimes_{\HH_G(V_1),\chi_1} \fpb \congto \kind V_2\otimes_{\HH_G(V_2),\chi_2} \fpb,
\end{equation*}
where $V_i$ are weights and $\chi_i : \HH_G(V_i) \to \fpb$ are algebra homomorphisms ($i = 1$, $2$).
Suppose we have such an isomorphism. By~(\ref{eq:11}), if $V_1$ occurs in $\pi$ with eigenvalues $\chi_1$,
then $V_2$ occurs in $\pi$ with eigenvalues $\chi_2$.

While the method is general, for the moment we can obtain an explicit criterion mainly when $G = \GL_n$.
The idea is to study the Hecke bimodule $\Hom_G(\kind V_1,\kind V_2)$. Our proof requires an
explicit version of the Satake transform, which we deduce from the corresponding result over the complex
numbers (the formula of Lusztig and Kato).

\subsubsection{Irreducibility proof}
\label{sec:irreduc}

Let $G = \GL_n$.  Suppose $\pi := \Ind_{\o P}^{G} (\sigma_1 \otimes \cdots \otimes \sigma_r)$ as in
Theorem~\ref{thm:class-intro}.  To prove irreducibility it is enough to show that any weight $V \subset \pi$ generates $\pi$
as $G$-representation.  Without loss of generality $V \INTO \pi$ is a Hecke eigenvector. If $V$ is sufficiently regular we
find a map $\pi_M \onto \sigma_1 \otimes \cdots \otimes \sigma_r$ and use the isomorphism~(\ref{eq:24}) to obtain a
surjection $\kind V \onto \pi$. This shows that $V$ generates $\pi$. Otherwise we use the criterion
of~\S\ref{sec:changing-weight} to change to a sufficiently regular weight and we are done as before. (To verify the criterion
we use conditions (ii) and (iii) in Theorem~\ref{thm:class-intro}.)

\subsubsection{Classification proof}
\label{sec:class}

Let $G = \GL_n$.  Suppose $\pi$ is an irreducible admissible $G$-representation.  Let $V$ be a weight of $\pi$ and let $\chi$
be a set of Hecke eigenvalues on $\Hom_K(V,\pi)$.  Then $\pi$ is a quotient of $\kind V \otimes_{\hg,\chi} \fpb$
by~(\ref{eq:11}). Suppose we can satisfy the conditions in \S\ref{sec:comp-inductions} for some proper parabolic $P \ne G$.
Then $\opInd (\pi_M) \onto \pi$ by \eqref{eq:24}. Using Emerton's theory of ordinary parts \cite{bib:emerton-ordinary1} one
gets a map $\opInd \sigma \onto \pi$ for some irreducible admissible $M$-representation $\sigma$ and one can proceed by
induction. Otherwise one tries to change the weight using the criterion of \S\ref{sec:changing-weight}. If both tactics fail,
it turns out that $\pi$ is either supersingular or $\pi$ looks like a one-dimensional representation in the sense that it
contains the same weight with the same Hecke eigenvalues as such a representation. A calculation with the Iwahori Hecke
algebra allows us to conclude in the latter case.

Both the irreducibility and the classification proofs go through for general~$G$ provided one puts conditions on the weights of the representation
so that all weights that occur in the arguments are sufficiently regular. See Theorems~\ref{thm:irred-parab-ind} and \ref{thm:irr-rep-as-ind-of-ss} for
precise statements.

\subsection{Other results}
\label{sec:other}

We complete the irreducibility proof of the generalised Steinberg representations $\Sp P$ of split reductive groups. This was
proved by Gro\ss e-Kl\"onne in case the root system is of type A, B, C, or D. In general we use his work to show that
$\Sp P$ contains a unique weight (with multiplicity one) and we determine the Hecke eigenvalues in that weight.
We then use \S\ref{sec:comp-inductions} to show that the unique weight generates the representation.

\subsection{Arrangement of the paper}
\label{sec:arrangement-paper}

In \S\ref{sec:hecke-actions-satake}--\ref{sec:parab-ind-comp-ind} we discuss Hecke actions, variants of the mod~$p$ Satake
transform, and the comparison result between compact and parabolic inductions.  In \S\ref{sec:hecke-eigenv} we give a
parameterisation of Hecke eigenvalues which is analogous to the classical parameterisation by unramified characters of the
torus. We also define what it means for an irreducible admissible representation to be supersingular. In
\S\ref{sec:compute-satake} we adapt the result of Lusztig and Kato to give an explicit version of the mod~$p$ Satake
isomorphism. In \S\ref{sec:maps-betw-cpt-ind} we study maps between compact inductions and deduce our ``change of weight''
criterion. In \S\ref{sec:gen-steinberg} we discuss generalised Steinberg representations. In
\S\ref{sec:irreducibility}--\ref{sec:classification} we deduce the main results. Finally in \S\ref{sec:submodule-structure}
we prove results on the submodule structure of parabolically induced representations.

\subsection{Notation}
\label{sec:notation}

Let~$F$ be a finite extension of~$\qp$ with ring of integers~$\O$, uniformiser~$\varpi$, and residue field~$k$ of order~$q$.
Let $G$ be a split connected reductive group over $\O$ and fix a maximal split torus $T$. Let $\Phi \subset X^*(T)$ denote
the set of roots. For each root $\alpha \in \Phi$ denote by $\alpha\dual$ the associated coroot.
Choose a system of positive roots~$\Phi^+$ and let $B = TU$ denote the associated Borel subgroup.
Let~$W$ be the Weyl group and let $K = G(\O)$, a hyperspecial maximal compact subgroup of $G(F)$.  For example, we could take
$G = \GL_n$ with the diagonal maximal torus $T$ and the Borel subgroup $B$ consisting of upper triangular matrices.

The letters $P$, $Q$ usually denote standard parabolic subgroups, i.e., parabolic subgroups containing the Borel $B$.
A Levi decomposition $P = M N$ of a standard parabolic is implicitly assumed to be standard, i.e., $M$ is 
the unique Levi subgroup of $P$ that contains $T$. If $P = MN$ is a parabolic we denote by $\o P = M \o N$ the opposite
parabolic (determined by our choice of Levi $M$). 
We will sometimes use that the multiplication map $\o N \times M \times N \to G$ is injective on $F$-points
(in fact it is even an open immersion on the level of $\O$-schemes \cite[\S II.1.11]{bib:Jan-reps}).

We denote by $Z$ the \emph{connected} centre of $G$. Similarly $Z_M$ denotes the connected centre of a Levi subgroup $M$.
We denote by $X_*(T)_-$ the set of antidominant coweights of $T$.
The set of simple roots in $\Phi^+$ is denoted by $\Delta$. If $P = MN$ is a standard parabolic, we similarly
define $\Delta_M$ (with respect to $\Phi^+_M = \Phi_M \cap \Phi^+$).

We denote by $\red : K = G(\O) \to G(k)$ the reduction map and by $K(1)$ its kernel, which is a pro-$p$ group.
For a standard parabolic subgroup $P$, we let $\P := \red^{-1}(P(k))$ be the corresponding parahoric subgroup.
We will also write $I$ for the standard Iwahori subgroup $\B$
and define $I(1) := \red^{-1}(U(k))$ (the pro-$p$ Sylow subgroup of $I$).

If $H$ is a group and $\sigma$ an $H$-representation, we denote by $\soc_H \sigma$ the $H$-socle, i.e., the largest
semisimple subrepresentation.

We usually write $G$, $P$, $M$, \ldots\ when we really mean $G(F)$, $P(F)$, $M(F)$, \ldots\ This should cause no
confusion.

We use the following abbreviated notations for the conjugation action: ${}^t K = t K t^{-1}$ and $K^t = t^{-1} K t$.

\emph{All representations in this paper, unless otherwise stated, live on $\k$-vector spaces.}

\subsection{Acknowledgements}
\label{sec:acknowledgements}

I am grateful to Matthew Emerton for his comments and particularly for making a very helpful observation
related to Proposition~\ref{prop:triv-wt-hecke-evals}.
I thank Christophe Breuil, Guy Henniart, Vytautas Pa\v sk\= unas, Peter Schneider, Marie-France Vign\'eras, and the
referee for useful comments.
I also thank the mathematics department at UCLA for the excellent working conditions during my stay
in the spring of 2009.

\tableofcontents

\section{Hecke actions and the Satake transform}
\label{sec:hecke-actions-satake}

\subsection{Background}
\label{sec:background}

\subsubsection{Weights}
\label{sec:wts}

\begin{df}\label{df:wts}
  A \emph{weight} is an irreducible representation~$V$ of the finite group~$G(k)$ over~$\k$.
\end{df}

The set of $q$-restricted weights is defined to be:
\begin{equation*}
  X_q(T) = \{\lambda \in X^*(T): 0 \le \langle \lambda, \alpha^\vee\rangle < q \quad \forall \alpha \in \Delta \}.
\end{equation*}
We also define:
\begin{equation*}
  X^0(T) = \{\lambda \in X^*(T): \langle \lambda, \alpha^\vee\rangle = 0 \quad \forall\alpha\in \Phi\}.
\end{equation*}
For $\nu \in X^*(T)$ dominant, let $F(\nu)$ denote the irreducible $G_{/\k}$-module of highest weight $\nu$.
Via the inclusion $G(k) \to G(\k)$, we can consider $F(\nu)$ as $G(k)$-representation.

\begin{prop}\label{prop:wts}
  Suppose that the derived subgroup of $G$ is simply connected. 
  Then all weights are of the form $F(\nu)$, where $\nu \in X_q(T)$. Moreover for $\nu$, $\nu' \in X_q(T)$,
  we have $F(\nu) \cong F(\nu')$ as $G(k)$-representations if and only if $\nu-\nu' \in (q-1)X^0(T)$.
\end{prop}

This goes back to Steinberg if $G$ is semisimple; in general, see Prop.~1.3 in the appendix to \cite{bib:thesis}.  If the
derived subgroup of $G$ fails to be simply connected, weights can be described using a $z$-extension of $G_{/k}$ (as in the
proof of \cite[Lemma~2.5]{bib:satake}).

We will sometimes denote by $1$ the trivial weight. If $P = MN$ is a standard parabolic, we denote
by $X^0_M(T)$ and $F^M(\nu)$ the analogues of $X^0(T)$ and $F(\nu)$ for the Levi $M$.

For the following, very useful Lemma see~\cite[Lemma 2.5]{bib:satake}. It was first proved by Smith and Cabanes.

\begin{lm}\label{lm:wts-invts}
  Suppose that $V$ is a weight and that $P = MN$ is a standard parabolic. Then $\vnk$ and $\vonk$ are weights for
  $M$ and the natural, $M(k)$-linear map $\vnk \to \vonk$ is an isomorphism. In particular, $\vuk \cong \vouk$ is
  one-dimensional.

  Suppose the derived subgroup of $G$ is simply connected. If $V \cong F(\nu)$ for some $\nu \in X_q(T)$, then $\vnk \cong
  F^M(\nu)$. Moreover, as a subspace of $F(\nu)$, $\vnk$ is the sum of all weight spaces $F(\nu)_{\nu'}$ with
  $\nu-\nu' \in \Z\Phi_M$.
\end{lm}

\begin{df}
  A weight $V$ is said to be \emph{$M$-regular} if $\Stab_W(V^{U(k)}) \subset W_M$.
\end{df}

Here $\Stab_W(V^{U(k)})$ denotes the set of $w \in W$ that preserve the one-dimensional, $T(k)$-stable subspace $V^{U(k)}
\subset V$.  Let us make this condition more explicit in case $V = F(\nu)$ for some $\nu \in X_q(T)$.  Note that
$\Stab_W(V^{U(k)}) = \Stab_W(\nu)$, which is generated by the simple reflections $s_\alpha$ for $\alpha \in \Delta$ such
that $\langle \nu,\alpha\dual\rangle = 0$. Thus $V = F(\nu)$ is $M$-regular if and only if $0 <
\langle{\nu,\alpha\dual}\rangle \le q-1$ for all $\alpha \in \Delta-\Delta_M$.

\begin{lm}\label{lm:lift-to-m-reg-wt}
  The map $V \mapsto \vnk$ from $M$-regular weights for $G$ to weights for $M$ is a bijection.
\end{lm}
  
\begin{proof}
  Suppose first that the derived subgroup of $G$ is simply connected. Then the same is true for $M$. Let $\o V$ be a
  weight for $M$.  Then $\o V \cong F^M(\nu)$ for some $\nu \in X^*(T)$ such that $0 \le \langle \nu,\alpha\dual\rangle \le
  q-1$ for all $\alpha \in \Delta_M$. We need to find $\nu' \in X^*(T)$ such that (a) $0 \le \langle \nu',\alpha\dual\rangle
  \le q-1$ for all $\alpha \in \Delta_M$, (b) $0 < \langle \nu',\alpha\dual\rangle \le q-1$ for all $\alpha \in \Delta -
  \Delta_M$, and (c) $\nu -\nu' \in (q-1)X^0_M(T)$. The first two conditions express that $\nu'$ is $q$-restricted and that
  the weight $V \cong F(\nu')$ is $M$-regular. Condition (c), which in fact implies (a), expresses that $\vnk \cong \o V$.  Clearly there is such a
  $\nu'$ and it is uniquely determined up to $(q-1)X^0(T)$. This completes the proof.
  
  In the general case, pick a $z$-extension $1 \to R \to \wt G\to G \to 1$ in the special fibre, just as in
  \cite[Lemma~2.5]{bib:satake}. We know there is a unique weight $V$ for $\wt G$ such that $V^{\wt N(k)} \cong \o V$ as $\wt
  M(k)$-representation. Since $R(k)$ acts trivially on $\o V$ and since it acts on $V$ via the central character, we see that $V$ descends to
  a $G(k)$-representation. The uniqueness of $V$ is even easier.
\end{proof}

\subsubsection{Smooth representations}
\label{sec:smooth}

We recall that a smooth $G$-representation $\pi$ is said to be \emph{admissible} if $\pi^H$ is finite-dimensional for all
open subgroups $H$ of $G$. It is sufficient to verify this condition for one open pro-$p$ subgroup $H$. An irreducible admissible $G$-representation $\pi$
has a central character, which we denote by $\omega_\pi : Z \to \k\s$.

Any irreducible smooth $K$-representation factors through $G(k)$ (as $K(1)$ is pro-$p$); it is thus a weight. If $\pi$ is
admissible, $\soc_K \pi \subset \pi^{K(1)}$ is finite-dimensional and non-zero, so $\pi$ contains a weight $V$.
(It is clear that even any smooth $G$-representation contains a weight.) We will
also say that $V$ is a weight of $\pi$.

For a closed subgroup $H \subset G$ and a smooth $H$-representation $\sigma$ we denote by $\Ind_H^G \sigma$ (resp., $\ind_H^G
\sigma$) the representation that is induced (resp., compactly induced) from $\sigma$.  If $H$ is open, then $\ind_H^G \sigma
\cong \k[G] \otimes_{\k[H]} \sigma$, so $\ind_H^G$ is a left adjoint to the forgetful functor. In this case, we will
denote by $[g,x] \in \ind_H^G \sigma$ the element that is supported on $Hg^{-1}$ and sends $g^{-1}$ to $x \in \sigma$. If
$P$ is a parabolic subgroup, then $\Ind_P^G$ is \emph{exact}. (This is because the map $G \to G/P$ has continuous sections;
see \cite[Prop.~4.1.5]{bib:emerton-ordinary1}.)

Suppose that $H$ is a compact open subgroup and that $V$ is a finite-dimensional smooth $H$-representation. We define the
\emph{Hecke algebra} of $V$ to be $\HH_H(V) := \End_G(\ind_H^G V)$.  Using the above adjunction, we can and usually will
think of it as the $\k$-algebra of compactly supported functions $\vp : G\to \End_{\k} V$ satisfying $\vp(h_1gh_2) = h_1
\circ \vp(g)\circ h_2$ for all $h_1$, $h_2 \in H$, $g \in G$, where the multiplication is given by convolution.
Explicitly, for $\vp \in \HH_H(V)$ and $f \in \ind_H^G V$, we have
\begin{equation*}
  \vp(f)(g) = \sum_{\gamma \in G/H} \vp(g\gamma)f(\gamma^{-1}).
\end{equation*}
Note that if $\pi$ is a smooth $G$-representation, $\HH_H(V)$ naturally acts on the multiplicity space $\Hom_H(V,\pi) \cong
\Hom_G(\ind_H^G V,\pi)$. (This is a right action.) Explicitly, if $\vp \in \HH_H(V)$ and $f : V \to \pi$ is $H$-linear, then
\begin{equation*}
(f \ast \vp)(v) = \sum_{\gamma \in H\bs G} \gamma^{-1} f(\vp(\gamma)v).
\end{equation*}

Suppose now that $V$ is a weight. We will usually write $\HH_G(V)$ instead of $\HH_K(V)$. We will see just below that
$\hg$ is commutative. If $V = 1$ then $\HH_G(V)$ is the usual unramified Hecke algebra $\k[K\bs G/K]$.

\subsection{The mod $p$ Satake transform}
\label{sec:satake-mod-p}

We begin by recalling some results of~\cite{bib:satake}.
Let~$T^-$ denote the submonoid of~$T$,
\begin{equation*}
  T^- = \{t \in T : \ord_F(\alpha(t)) \le 0 \quad \forall \alpha \in \Delta\},
\end{equation*}
and let $\HH^-_T(V^{U(k)})$ denote the subalgebra of~$\HH_T(V^{U(k)})$ consisting of those functions $\varphi : T \to \k$ that are
supported on~$T^-$.

\begin{thm}\label{thm:satake} Suppose that $V$ is a weight.
  Then
  \begin{align*}
    \SS_G : \HH_G(V) &\to \HH_T(V^{U(k)})\\ 
    \vp &\mapsto \left(t \mapsto \sum_{u \in U/U(\O)} \vp(tu)\Big|_{V^{U(k)}}\right)
  \end{align*}
  is an injective $\k$-algebra homomorphism with image $\HH^-_T(V^{U(k)})$.
\end{thm}

In particular, $\hg \cong \k[X_*(T)_-]$ is commutative and noetherian (Gordan's lemma shows that $X_*(T)_-$ is finitely generated). We recall that $G = \coprod
K\lambda(\varpi)K$, where $\lambda$ ranges over $X_*(T)_-$ (refined Cartan decomposition). Moreover, $\hg$ has a basis given
by $T_\lambda$ for $\lambda \in X_*(T)_-$, where $T_\lambda$ has support $K\lambda(\varpi)K$ and sends $\lambda(\varpi)$ to the
endomorphism $V \onto V_{N_\lambda(k)} \tocong V^{N_{-\lambda}(k)} \INTO V$ (see~\S\ref{sec:various-lemmas} for the
definition of $P_\lambda = M_\lambda N_\lambda$). We also denote by $\tau_\lambda \in \HH_T(V^{U(k)})$ the element supported on
$\lambda(\varpi)T(\O)$ that sends $\lambda(\varpi)$ to 1. We have:
\begin{equation}\label{eq:10}
  \SS_G(T_\lambda) = \sum_{\substack{\mu \in X_*(T)_- \\ \mu \ge_{\R} \lambda}} a_\lambda(\mu) \tau_\mu, \quad
  \text{with $a_\lambda(\mu) \in \k$ and $a_\lambda(\lambda) = 1$}.
\end{equation}
Here, $\mu \ge_\R \lambda$ means that $\mu-\lambda$ is a non-negative real linear combination of the simple coroots.

For $z \in Z$ we also define $T_z \in \hg$ such that $\supp(T_z) = Kz$ and $T_z(z) = \id_V$.
Thus $T_{\lambda(\varpi)} = T_\lambda$ whenever $\lambda(\varpi) \in Z$.
The following formulae for $z$, $z_1$, $z_2 \in Z$, $z_0 \in Z(\O)$ will be useful later.
\begin{gather}
  T_{z_1z_2} = T_{z_1} T_{z_2}, \ T_1 = 1, \label{eq:2} \\
  T_{z_0 z} = \omega_V(z_0)^{-1} T_z, \label{eq:5} \\
  T_z = z^{-1} \quad \text{on $\kind V$}, \label{eq:3}
\end{gather}
where $\omega_V : Z(k) \to \k\s$ denotes the central character of $V$.

\subsection{Variants of the Satake transform}
\label{sec:satake-variants}

Let $P = MN$ be a standard parabolic subgroup.  There is a ``partial'' Satake homomorphism $\psmg : \hg \to \HH_{M} (V_{\o
  N(k)})$. It is defined by
\begin{align*}
  \psmg : \hg &\to \HH_{M} (V_{\o N(k)}) \\
  \vp &\mapsto \left (m \mapsto p_{\o N} \sum_{\o N(\O)\bs \o N} \vp(\o n m) \right),
\end{align*}
where $p_{\o N}$ denotes the projection $V \onto \vonk$. It is easy to check that this map is well defined (just as for $\sg$
\cite{bib:satake}).  We also have $\smg : \hg \to \HH_{M} (\vnk)$, defined in the same way as $\SS_G$.

We have a simple compatibility between $\psmg$ and $\smg$. There is an algebra isomorphism $\hg \congto \hgd$ which sends
$\vp$ to $\vp'$ with $\vp'(g) = \vp(g^{-1})^*$. (Recall that $\hg$ is commutative.)  Similarly we have an isomorphism
$\HH_M(V^{\o N(k)}) \congto \HH_M((V^*)_{\o N(k)})$. In the following lemma we denote by $\smgop$ the analogue of $\smg$ that is
computed with respect to $\o P = M \o N$ (a standard parabolic for $\Phi^-$).

\begin{lm}\label{lm:satake-duality} We have the following commutative diagram.
  \ifkuvio
  \[ \cellwidth=45pt \Diag
  \hg & \rTo^{\smgop} & \HH_M(V^{\o N(k)}) \\
  \dTo_\cong  && \dTo^\cong  \dy{-3mm} \\
  \hgd & \rTo^{\psmg} & \HH_M((V^*)_{\o N(k)}) \\ 
  \endDiag \]
  \fi
  In particular $\psg := {}'\SS^T_G$ is injective. Its image consists of those elements of $\htd$ that are supported on $T^-$
  and is isomorphic to $\k[X_*(T)_-]$.
\end{lm}

\begin{proof}
  This is straightforward, noting that $V^{\o N(k)} \to V$ is dual to $V^* \to (V^*)_{\o N(k)}$.
\end{proof}

We will see in Cor.~\ref{cor:smg-equals-psmg} that $\smg$ and $\psmg$ are actually the same, under the natural identification
of $\vnk$ and $\vonk$ (at least when the derived subgroup of $G$ is simply connected). This does not seem to be obvious from
the definition. \mar{What about using compat. with $p$-adic Satake?}

\begin{prop}\label{prop:partial-satake}
  The map $\psmg$ is an injective algebra homomorphism. It is a localisation
  map of integral domains. Also $\psg = \psm \circ \psmg$.
\end{prop}

The same of course also holds for $\smg$.

\begin{proof}
  It is easy to see that $\psmg$ is an algebra homomorphism (see also \cite{bib:satake}). A direct calculation shows that
  $\psg = \psm \circ \psmg$. 
  Since $\psg$ is injective, so is $\psmg$. As we have fixed a uniformiser $\varpi$, we can identify $\HH_T(V_{\o U(k)})$
  with $\k[X_*(T)]$. By Lemma~\ref{lm:satake-duality} the image of $\hg$ in $\k[X_*(T)]$ is $\k[X_*(T)_-]$ and the image of
  $\hm$ in $\k[X_*(T)]$ is $\k[X^M_*(T)_-]$, where $X^M_*(T)_-$ is the monoid of antidominant coweights for $T$ in $M$.
  Finally note that $\k[X_*(T)_-] \INTO \k[X^M_*(T)_-]$ is a localisation map of integral domains: it suffices to invert any
  $\lambda \in X_*(T)_-$ such that 
  for $\alpha \in \Delta$, $\langle \lambda,\alpha\rangle = 0$ if and only if $\alpha \in \Delta_M$.
\end{proof}

Suppose now that $V$ is a weight for $G$ and $\sigma$ a smooth $M$-representation. By the Iwasawa decomposition we have
$(\opInd \sigma)|_K \cong \Ind_{\o P(\O)}^K \sigma$, so the natural map
\begin{align}\label{eq:6}
  \Hom_K(V,\opInd \sigma) &\to \Hom_{M(\O)} (\vonk,\sigma) \\
  \notag f &\mapsto \o f = \big(\o v \mapsto f(v)(1)\big)
\end{align}
is a bijection (here $v$ is any lift of $\o v$). Note that $\hg$ naturally acts on the left-hand side, and on the right-hand side
via $\psmg$.

\begin{lm}\label{lm:hecke-parab-ind}
  The natural map in \eqref{eq:6} is $\hg$-equivariant.
\end{lm}

\begin{proof}
  For $\vp \in \hg$ let $\vp_M = \psmg(\vp)$. We compute
  \begin{align*}
    (\o f \ast \vp_M)(\o v) &= \sum_{M(\O)\bs M} m^{-1} f\bigg( \sum_{\o N(\O)\bs \o N} \vp (\o n m) v\bigg)(1) \\
    &= \sum_{M(\O)\bs M} \sum_{\o N(\O)\bs \o N} f( \vp (\o n m) v)(m^{-1} \o n^{-1})
  \end{align*}
  and
  \begin{equation*}
    (\o{f \ast \vp})(\o v) = (f \ast \vp)(v)(1) = \sum_{K \bs G} f(\vp(g)v)(g^{-1}).
  \end{equation*}
  By the Iwasawa decomposition $G = K \o P$ these two expressions are equal.
\end{proof}

\subsection{Various lemmas}
\label{sec:various-lemmas}

Suppose that $\lambda \in X_*(T)$. Let $P_\lambda = M_\lambda N_\lambda$ denote the parabolic subgroup of $G$
defined by $\lambda$.  For the following proposition, see
\cite[Prop.~3.8]{bib:satake}.

\begin{prop}\label{prop:buildings-lemma}
  Let $t = \lambda(\varpi)$. Then $\red(K \cap K^t) = P_\lambda(k)$.
\end{prop}

Define $Z_M^- := Z_M \cap T^-$ and $Z_M^+ := Z_M \cap (T^-)^{-1}$. We also need
\begin{equation*}
  Z_M^{--} := \{ z \in Z_M : \ord_F(\alpha(z)) < 0 \quad \forall \alpha \in \Delta-\Delta_M \}.
\end{equation*}

\begin{lm}\label{lm:iwahori-decomp}
  Suppose $P = MN$ is a standard parabolic subgroup. We let $\oP^- = \o N \cap \oP$, $\oP^0 = M \cap \oP$, $\oP^+ = N \cap \oP$.
  Then $\oP = \oP^- \oP^0 \oP^+$ \(in any order\). Moreover, conjugation by $T^-$ contracts $\oP^-$ and expands $\oP^+$.

  If $h_0 \in Z_M^{--}$ then the sets $(\oP^+)^{h_0^n}$ for $n \ge 0$ form a neighbourhood basis of the identity in $N$.
\end{lm}

By interchanging positive and negative roots we obtain similarly that $\P = \P^- \P^0 \P^+$ in any order
with $T^-$ contracting $\P^-$ and expanding $\P^+$.

\begin{proof}
  Note that $\oP^- = \o N(\O)$, $\oP^0 = M(\O)$, and $\oP^+ = \ker(N(\O)\to N(k))$.
  For any $\O$-group scheme~$H$ let us write ${}^1 H(\O)$ for $\ker (H(\O) \to H(k))$.
  Assume for now that $K(1) = {}^1 \o U(\O) {}^1 T(\O) {}^1 U(\O)$. It is not hard to prove this by an elementary
  argument when $G = \GL_n$; we will justify it in general below.

  From our assumption, $K(1) \subset \o P(\O) {}^1 P(\O)$. From $\o P(\O) = \o N(\O) M(\O)$ and $P(\O) = M(\O)N(\O)$ it
  follows that $K(1) \subset \oP^- \oP^0 \oP^+$. Thus $\oP = \o P(\O) K(1) \subset \oP^- \oP^0 \oP^+$. The reverse inclusion
  is obvious.  Thus $\oP = \oP^- \oP^0 \oP^+$. Noting that $\oP^0$ normalises both $\oP^-$ and~$\oP^+$, and by using the
  inverse, we see that the order does not matter. (Moreover, it is a direct product decomposition because multiplication $\o
  N \times M \times N \to G$ is injective.)
  
  Let us show that $T^-$ contracts~$\oP^-$ and expands~$\oP^+$. Choose root homomorphisms $x_\alpha : \G_a \congto U_\alpha$
  \cite[II.1.2]{bib:Jan-reps}. Fix an ordering of the set of roots~$\Phi$.  Note that $\o N = \prod_{\Psi^-} U_\alpha$ and $N
  = \prod_{\Psi^+} U_\alpha$ (as $\O$-schemes), where $\Psi^\pm = \Phi^\pm-\Phi^\pm_M$. Thus $\oP^- = \prod_{\Psi^-} x_\alpha(\O)$ and $\oP^+ =
  \prod_{\Psi^+} x_\alpha(\varpi\O)$.  The claim follows since $t x_\alpha(u) t^{-1} = x_\alpha(\alpha(t) u)$ for $t \in
  T(F)$ and $u \in F$.

  Now fix an element $h_0 \in Z_M^{--}$. The multiplication map $\prod x_\alpha : \prod_{\Psi^+} \G_a \to N$ induces an
  isomorphism of $\O$-schemes, thus also a homeomorphism $\prod_{\Psi^+} F \to N(F)$. Write again $\oP^+ = \prod_{\Psi^+} x_\alpha(\varpi\O)$,
  so that $(\oP^+)^t = \prod_{\Psi^+} x_\alpha(\alpha(t)^{-1}\varpi\O)$ for any $t \in T(F)$. By the definitions,
  $\ord_F(\alpha(h_0^n)) \le -n$ for all $\alpha \in \Psi^+$ and all $n \ge 0$. Thus $(\oP^+)^{h_0^n}$ for $n \ge 0$ form a neighbourhood
  basis of the identity in~$N(F)$.
  
  Finally we justify that $K(1) = {}^1 \o U(\O) {}^1 T(\O) {}^1 U(\O)$ using Bruhat--Tits theory. Fraktur letters denote the
  $\O$-group schemes defined in \cite{bib:BT2}.  Since $G$ is a reductive group over~$\O$ with special fibre~$G_{/F}$, there
  is a (hyper-)special point~$x$ in the reduced building of~$G_{/F}$ such that $G \cong \fG_x^0$
  \cite[5.1.40]{bib:BT2}. 
  We identify these two group schemes. Let $f := f_x'$ be the corresponding ``quasi-concave'' function on the set of roots
  \cite[4.6.24]{bib:BT2}, so that $G = \fG_f^0$. Now this $\O$-group scheme has two split maximal tori, $T$ and~$\fT$.
  By~\cite[Exp.~XXVI, Prop.~6.16]{bib:SGA3} there is an element $\gamma \in G(\O)$ that conjugates $\fT$ to~$T$. In the
  generic fibre the two unipotent subgroups, $U_{/F}$ and~${}^\gamma (\fU^+_{/F})$ might correspond to two different choices
  of positive roots. By multiplying $\gamma$ by a suitable element of~$N(T)(\O)$ we can assure that $U$ and~${}^\gamma \fU^+$
  agree on the generic fibre. But then they are the same (see for example \cite[1.2.6]{bib:BT2}).  By~\cite[4.6.8]{bib:BT2}
  we have ${}^1 \fG_f^0(\O) = {}^1 \fU^-(\O) {}^1 \fT(\O) {}^1 \fU^+(\O)$. Conjugating by~$\gamma$ we obtain $K(1) = {}^1
  G(\O) = {}^1 \o U(\O) {}^1 T(\O) {}^1 U(\O)$.
\end{proof}

\begin{lm}\label{lm:support}
  Suppose that $P = MN$ and $Q = LN'$ are standard parabolics. Suppose that
  $V$ is both $M$ and $L$-regular. Then $p_{\o {N'}} (\kappa \vnk) \ne 0$ for $\kappa \in K$ implies $\kappa \in \oQ \P$.
  
  If $\Stab_W(\vuk)$ equals $W_M$ \(resp., $W_L$\), we may drop the assumption that $V$ be $L$-regular
  \(resp., $M$-regular\).
\end{lm}

\begin{proof}\mar{simpler proof? can show non-zero on some extremal weight?}
It is easy to reduce to the case that the derived subgroup of $G$ is simply connected,
just as in \cite[Lemma~2.5]{bib:satake}. (Note that the statements we want to prove just concern the finite group $G(k)$.)
In this case $V = F(\nu)$ for some $\nu \in X_q(T)$, moreover $V^{N(k)} = F(\nu)^N$.

We claim that for all weights $\mu$ of $F(\nu)^{N}$ and all $\alpha \in \Phi^+ - \Phi_M^+$ we have $\langle
\mu,\alpha\dual\rangle > 0$.  We know that $F(\nu)^N$ is the irreducible $M$-representation of highest weight $\nu$, so its
weights lie in the convex hull of $w\nu$ ($w \in W_M$). Note that any $W_M$ preserves $\Phi^+ - \Phi_M^+$. (Reduce to the
case of a simple reflection $s_\beta$ with $\beta \in \Phi_M$.  It preserves $\Phi_M$ and $\Phi^+-\{\beta\}$.) Therefore
$\langle w\nu,\alpha\dual\rangle = \langle \nu,w^{-1}\alpha\dual\rangle \ge 0$ for all $w \in W_M$.  If equality holds, then
$s_{w^{-1}\alpha}(\nu) = \nu$. Thus $s_{w^{-1}\alpha} \in W_M$ (as $V$ is $M$-regular), which implies that $\alpha \in
\Phi_M$, a contradiction. The claim follows.

By the rational Bruhat decomposition, $G(k) = \coprod_{W_L \bs W/W_M} \o Q(k) \dot\sigma P(k)$. It will thus suffice to show
that $p_{\o {N'}} (\dot\sigma \vnk) \ne 0$ for $\sigma \in W$ implies $\sigma \in W_LW_M$. By Lemma~\ref{lm:wts-invts} we see
that $V = V^{N'(k)} \oplus \ker p_{\o {N'}}$ and that these two subspaces are preserved by the torus $T$ acting on $V =
F(\nu)$, sharing no common weight.  (Note that $p_{\o {N'}}$ is identified with $F(\nu) \to F(\nu)_{\o {N'}}$. This
follows, for example, since the natural surjection $V_{\o {N'}(k)} \to F(\nu)_{\o {N'}}$ is an isomorphism, since the domain
is irreducible as $L(k)$-representation.)  So there is a weight~$\mu$ of $F(\nu)^N$ such that $\sigma\mu$ is a weight of
$F(\nu)_{\o{N'}} \cong F(\nu)^{N'}$. By the previous paragraph,
\begin{equation*}
  \langle \mu,\alpha\dual\rangle > 0 \quad \forall \alpha \in \Phi^+ - \Phi_M^+,
  \quad  \langle \sigma\mu,\beta\dual\rangle > 0 \quad \forall \beta \in \Phi^+ - \Phi_L^+.
\end{equation*}
It follows that $\sigma(\Phi^+ - \Phi_M^+) \subset \Phi^+ \cup \Phi_L^-$.

We claim that there is a $w \in W_M$ such that $\sigma w(\Phi^+) \subset \Phi^+ \cup \Phi_L^-$. Suppose there is a simple root
$\alpha$ of $M$ such that $\sigma(\alpha) < 0$. Then $s_\alpha$ preserves $\Phi^+ - \Phi_M^+$, so
$\sigma s_\alpha(\Phi^+ - \Phi_M^+) \subset \Phi^+ \cup \Phi_L^-$ while $\sigma s_\alpha$ maps one fewer simple root of $M$ to a negative
root. By induction we find a $w \in W_M$ such that $\sigma w(\Phi^+ - \Phi_M^+) \subset \Phi^+ \cup \Phi_L^-$
and $\sigma w$ maps all simple roots of $M$ to positive roots. This implies the claim.

Equivalently, $w^{-1}\sigma^{-1}(\Phi^- - \Phi_L^-) \subset \Phi^-$. The same argument
as in the previous paragraph shows that there is a $w' \in W_L$ such that $w^{-1}\sigma^{-1} w' (\Phi^-) \subset \Phi^-$.
This shows that $\sigma = w' w^{-1} \in W_L W_M$, which completes the proof.

To justify the final statement, suppose that $\Stab_W(\vuk) = W_M$. Then $\vnk = \vuk$ is one-dimensional. (One can see this
directly using the Bruhat decomposition; alternatively note that $\nu \in X^0_M(T)$.) Thus if $p_{\o {N'}} (\dot\sigma \vnk)
\ne 0$ then $\sigma\nu$ is a weight of $F(\nu)^{N'}$. Pick simple roots $\alpha_i \in \Delta$ with corresponding simple
reflections $s_i \in W$ such that
\begin{equation*}
  \nu \gneq s_1 \nu \gneq \dots \gneq s_r\cdots s_2s_1 \nu = \sigma\nu.
\end{equation*}
By Lemma~\ref{lm:wts-invts} it follows that $\nu - \sigma\nu \in \Z \Phi_L$.  So each $\alpha_i$ is in $\Phi_L$, hence
$\sigma \in s_r\cdots s_2s_1 \Stab_W(\nu) \subset W_L W_M$.  The argument in case $\Stab_W(\vuk) = W_L$ is similar, or
follows by duality.
\end{proof}

\begin{coroll}\label{cor:support}
  Suppose that $P = MN$ is a standard parabolic and that $\lambda \in X_*(T)_-$. Suppose that $V$ is $M$-regular, and suppose that
  either $\Stab_W(\vuk) = W_M$ or that $V$ is also $M_\lambda$-regular.
  Let $t = \lambda(\varpi)$.
  \begin{enumerate}
  \item If $T_\lambda(g)|_{\vnk} \ne 0$ then $g \in Kt \P$. If $p_{\o N} \circ T_\lambda(g) \ne 0$ then $g \in \oP t K$.
  \item We have $\smg(T_\lambda) = T_\lambda^M$ and $\psmg(T_\lambda) = T_\lambda^M$.
  \end{enumerate}
\end{coroll}

\mar{Alternatively could assume that $\Stab_W(\vuk) = W_\lambda$!}
This generalises \cite[Proposition~1.4]{bib:satake}. Note that if $\lambda(\varpi) \in Z_M^-$, then an $M$-regular
$V$ is also $M_\lambda$-regular.

\begin{proof}
  (i) Writing $g = \kappa' t \kappa$ with $\kappa'$, $\kappa \in K$, we see that $T_\lambda(t)\kappa|_{V^{N(k)}} \ne 0$, so
  $p_{N_\lambda}(\kappa V^{N(k)}) \ne 0$. Lemma~\ref{lm:support} shows that $\kappa \in \P_\lambda \cdot \P$. By Prop.~\ref{prop:buildings-lemma},
  $\red(\P_\lambda) = \red(K\cap K^t)$, so $\kappa \in (K \cap K^t) \P$ and $g = \kappa't\kappa \in Kt\P$.

  The other part is similar, or follows by duality.
  
  (ii) Take any $m \in M$ and suppose that $T_\lambda(mn)|_{\vnk} \ne 0$ (some $n \in N$) is a non-zero term contributing to $(\smg
  T_\lambda)(m)$. By part (i) and by Lemma~\ref{lm:iwahori-decomp} (twice), we have
  \begin{equation*}
    mn \in Kt\P \cap P = KtP(\O) \cap P = P(\O)tP(\O) = M(\O)tM(\O)N(\O). 
  \end{equation*}
  Thus $m \in M(\O)tM(\O)$ and $n \in N(\O)$. This shows that $\smg(T_\lambda)$ is supported on $M(\O)tM(\O)$ and that $(\smg
  T_\lambda)(t) = T_\lambda(t)|_{\vnk}$ is a linear projection (which is non-zero since $\vuk \subset \vnk$), so
  $\smg(T_\lambda) = T_\lambda^M$.  The other statement follows similarly or by using duality.
\end{proof}

\begin{coroll}\label{cor:smg-equals-psmg}
  Suppose that the derived subgroup of $G$ is simply connected.
  Let $P = MN$ be a standard parabolic subgroup.
  Under the natural identification of $\HH_M(\vnk)$ and $\HH_M(\vonk)$ we have $\smg = \psmg$.
\end{coroll}

\begin{proof}
  We first show that $\sg = \psg$ when $\dim V = 1$. We will make use of the ``Gelfand involution'' on $\hg$ (see~\cite{bib:satake}, end
  of Section~2.1). Let $\tau : G \to G$ denote a ``transpose'' involution as in \cite[II.1.16]{bib:Jan-reps}; it fixes $T$ pointwise
  and interchanges $U_\alpha$ and $U_{-\alpha}$ for all $\alpha \in \Phi$. Let ${}^\tau V$ denote $\Hom(V,\k)$ with $G(k)$ acting
  via $\tau$. As the derived subgroup of $G$ is simply connected, there is a $G(k)$-linear isomorphism $\varsigma : V \to {}^\tau V$.
  A linear map $\vp \in \End_\k V$ induces ${}^\tau \vp \in \End_\k V$ by dualising and applying $\varsigma$. We verified in~\cite{bib:satake}
  that ${}^\tau f({}^\tau g) = f(g)$ for all $f \in \hg$ and all $g \in G$.
  As $\dim V = 1$, we have $\End_\k V = \k$ and thus ${}^\tau \vp = \vp$. It follows that
  \begin{gather*}
    \sg(f)(t) = \sum_{U/U(\O)} f(tu) = \sum_{U/U(\O)} f({}^\tau (tu)) = \sum_{\o U(\O)\bs \o U} f(\o u t) = \psg(f)(t).
  \end{gather*}
  
  Next we show that $\sg = \psg$ in general. Let $M$ be the standard Levi with $\Stab_W (\vuk) = W_M$. Under the natural
  identification, we have by Cor.~\ref{cor:support} that $\smg(T_\lambda) = T_\lambda^M = \psmg(T_\lambda)$, for any $\lambda
  \in X_*(T)_-$.  Note that $\vnk$ is one-dimensional (for example, by the Bruhat decomposition), so by the above, $\sm =
  \psm$ for this weight.  By transitivity, $\sg = \psg$.

  Finally if $M$ is arbitrary we use that $\psg = \psm \circ \psmg$ and $\sg = \sm \circ \smg$. By the above we already know that
  $\sg = \psg$ and $\sm = \psm$. The claim follows by the injectivity of the latter map.
\end{proof}

\begin{lm}\label{lm:shrink-support}
  Suppose $\o P = M\o N$. Then $G = \oP T K$.

  Suppose the subset $X \subset G$ has finite image in $\oP \bs G$.
  Then there exists $h \in Z_M^-$ such that $h X \subset \oP T^- K$.
\end{lm}

See also \cite[Lemma~12]{bib:schneider-stuhler}.

\begin{proof}
  For the first claim, suppose $g \in G$. By the Cartan decomposition, there is a $\lambda \in X_*(T)_-$ such
  that $g \in K\lambda(\varpi)K$. By Prop.~\ref{prop:buildings-lemma}, we have
  \begin{equation*}
    \oP \bs K\lambda(\varpi)K/K \cong \o P(k) \bs G(k) / P_{-\lambda}(k).
  \end{equation*}
  The rational Bruhat decomposition shows that $g \in \oP \dot w \lambda(\varpi) K = \oP ({}^w\lambda(\varpi)) K$,
  for some $w \in W$. Alternatively one could use the Bruhat--Tits decomposition $G = \oB N(T) \oB$.

  For the second claim, say $X \subset \bigcup \oP t_i k_i$ (finite union). We now use that $(G,\oB,N(T))$ is a generalised
  Tits system \cite[\S 1]{bib:Iwahori}. The relevant axioms are as follows: (i) $H := N(T) \cap \oB$ is normal in $N(T)$,
  (ii) $N(T)/H$ is a semidirect product of a subgroup $\Omega$ and a normal subgroup $W'$, (iii) there exists a generating
  set $S = \{w_i : i \in I\}$ of $W'$ each of whose elements has order~2, (iv) for any $\sigma \in \Omega W'$ and any $i \in
  I$, $\sigma \oB w_i \subset \oB \sigma w_i \oB \cup \oB \sigma \oB$, 
  (v) any element of $\Omega$ normalises $\oB$.  Note that $H = T(\O)$ in our case. We remark that
  it is in fact not hard to justify using Bruhat--Tits theory that $(G,\oB,N(T))$ is a generalised Tits system, namely by applying
  \cite[Thm.~6.5]{bib:BT1}; so $S$ consists of the reflections in the walls of the chamber $C$ corresponding to $\oB$
  and $\Omega$ is the stabiliser of~$C$, as a subset of the apartment of~$T$.

  For $G = \GL_n$ it is not hard to write down $\Omega$ and $S$ explicitly \cite[\S 2]{bib:Iwahori}. 

  It follows easily from the above axioms that for all $n' \in N(T)$ there are $n_1$, \ldots, $n_r \in N(T)$ such that
  \begin{equation*}
    n \oB n' \subset \bigcup \oB n n_j \oB \quad \forall n \in N(T).
  \end{equation*}
  (Just write the image of $n'$ in $N(T)/H$ as a product of elements in $S \cup \Omega$ and induct.)
  Thus there are finitely many $n_{ij} \in N(T)$ such that $h \oB t_i \subset \bigcup_j \oB h n_{ij} \oB$ for
  all $i$ and for all $h \in Z_M^-$. It follows using Lemma~\ref{lm:iwahori-decomp} that
  \begin{equation*}
    h \oP t_i k_i = \oP^- \oP^0 h \oP^+ t_i k_i \subset \bigcup_j \oP h n_{ij} K,
  \end{equation*}
  as $\oP^+ \subset \oB \subset \oP$. Writing $n_{ij} = t_{ij} \dot w_{ij}$ with $t_{ij} \in T$ and $w_{ij} \in W$ we see that
  the right-hand side equals $\bigcup_j \oP h t_{ij} K$. Since $M(\O) \subset \oP$ we may replace $t_{ij}$ with any $W_M$-conjugate
  without changing its double coset. In this way we can ensure that $t_{ij}$ is antidominant as element of $M$.
  Then it is possible to find an $h \in Z_M^-$ such that $h t_{ij} \in T^-$ for all $i$, $j$.
\end{proof}

\subsection{Compatibilities between Hecke actions}
\label{sec:comp-betw-hecke}

Let $\o V$ be a weight for $M$.
Consider the following subspaces of Hecke algebras:
\begin{gather*}
  \HH'_M(\o V) := \{ \vp : \supp(\vp) \subset M(\O) Z_M^- M(\O) \} \subset \HH_M(\o V), \\
  \HH'_\P(\o V) := \{ \vp : \supp(\vp) \subset \P Z_M^- \P \} \subset \HH_\P(\o V), \\
  \HH'_\oP(\o V) := \{ \vp : \supp(\vp) \subset \oP Z_M^- \oP \} \subset \HH_\oP(\o V).
\end{gather*}
The following lemma shows that each of them is a (commutative) subalgebra and that we can naturally identify them.

\begin{lm}\label{lm:hecke-compat-m-p}
  The subspaces $\HH'_M(\o V)$, $\HH'_\P(\o V)$, $\HH'_\oP(\o V)$ are subalgebras. The map $\vp \mapsto \vp|_M$ from $\HH'_\P(\o V)$ to $\HH'_M(\o V)$
  \(resp., from $\HH'_\oP(\o V)$ to $\HH'_M(\o V)$\)
  is an algebra isomorphism.
\end{lm}

See also Vign\'eras \cite[\S A.7]{bib:vigneras2} for the subalgebras $\HH'_\P(\o V)$ and $\HH'_\oP(\o V)$.

\begin{proof}
  Clearly $\HH'_M(\o V)$ is a subalgebra. Once we show that the two maps are bijective and compatible with the algebra
  structures, it will follow that the other two subspaces are subalgebras.  We concentrate on the first map. The verification
  for the second map is similar or may be deduced by duality.

  The map is clearly injective. To check surjectivity, note that the $T_h^M$ for $h \in Z_M^-$ span $\HH'_M(\o V)$.  For $h
  \in Z_M^-$ define $E_h : G \to \End(\o V)$ in $\HH'_\P(\o V)$ with support $\P h \P$ and such that $E_h(h) = \id_{\o V}$.
  Let us check that it is well defined. If $\p_1 h = h \p_2$ for $\p_i \in \P$, write $\p_i = \p_i^- \p_i^0 \p_i^+$ according
  to the Iwahori decomposition $\P = \P^-\P^0\P^+$ (Lemma~\ref{lm:iwahori-decomp}). It follows
  that $\p_1^- (\p_1^0 h) (\p_1^+)^h = {}^h \p_2^- (h\p_2^0) \p_2^+$ in $\o N M N$, so $\p_1^0 = \p_2^0$.
  Thus $\p_1 E_h(h) = \p_1^0 = \p_2^0 = E_h(h) \p_2 \in \End(\o V)$ and $E_h$ is well defined. As $\P h\P = \P^- (\P^0 h)\P^+$, we see
  that $E_h|_M =  T_h^M$.

  To check that the map is compatible with the algebra structures, by~\eqref{eq:2} it suffices to show that $E_{h_1} \ast
  E_{h_2} = E_{h_1h_2}$ ($h_i \in Z_M^-$). We have
  \begin{equation*}
      \P h_1 \P h_2 \P = \P (h_1\P^- \P^0) (\P^+h_2) \P = \P h_1h_2 \P
  \end{equation*}
  by Lemma~\ref{lm:iwahori-decomp} and
  \begin{equation*}
    (E_{h_1} \ast E_{h_2})(h_1h_2) = \sum_{\P/(\P \cap {}^{h_1} \P)} E_{h_1}(\p h_1) E_{h_2}(h_1^{-1} \p^{-1} h_1 h_2).
  \end{equation*}
  Since $\P \cap {}^{h_1} \P = ({}^{h_1} \P^-) \P^0 \P^+$, we may replace the index set by $\P^-/{}^{h_1} \P^-$.
  To obtain a non-zero term for $\p \in \P^-$ we also need $(h_1^{-1} \p^{-1} h_1) h_2 \in \P h_2 \P = \P^-(\P^0 h_2) \P^+$.
  As the product map $\o N \times M \times N \to G$ is injective, $\p \in {}^{h_1} \P^-$. Thus
  $(E_{h_1} \ast E_{h_2})(h_1h_2) = \id_{\o V}$.
\end{proof}

Suppose that $V$ is an $M$-regular weight. Consider the subspace
\begin{equation*}
  \HH'_G(V) := \{ \vp : \supp(\vp) \subset K Z_M^- K \} \subset \hg.
\end{equation*}
The following lemma shows that it is a subalgebra and that we can identify it with $\HH'_\P(\vnk)$ and with
$\HH'_\oP(\vonk)$. Moreover, the resulting identification $\HH'_\P(\vnk) \cong \HH'_G(V) \cong \HH'_\oP(\vonk)$ is the one
obtained in Lemma~\ref{lm:hecke-compat-m-p}. (Recall that $\vnk \congto \vonk$ is a weight for $M$.)

\begin{lm}\label{lm:i_P}
  Assume that $V$ is $M$-regular.
  The subspace $\HH'_G(V)$ is a subalgebra. There is an algebra isomorphism $i^{\P} : \HH'_G(V) \to \HH'_\P(\vnk)$, which is
  characterised as follows. For any $\vp \in \HH'_G(V)$ the map $i^\P(\vp)$ is supported on $\P Z_M^-
  \P$ and for $g \in \P Z_M^- \P$, $i^\P(\vp)(g) \in \End(\vnk)$ is the restriction of $\vp(g) \in \End(V)$ to $\vnk$.

  Similarly there is an algebra isomorphism $i_{\oP} : \HH'_G(V) \to \HH'_\oP(\vonk)$ such that for any $\vp \in \HH'_G(V)$ the map
  $i_\oP(\vp)$ is supported on $\oP Z_M^- \oP$ and for $g \in \oP Z_M^- \oP$, $i_\oP(\vp)(g) \in \End(\vonk)$ is induced by
  $\vp(g) \in \End(V)$.

  The resulting identification of $\HH'_\P(\vnk)$ with $\HH'_\oP(\vonk)$ is the one obtained in
  Lemma~\ref{lm:hecke-compat-m-p}. The resulting identification of $\HH'_G(V)$ with $\HH'_M(\vnk)\cong \HH'_M(\vonk)$ is
  given by $\smg = \psmg$.
\end{lm}

\begin{proof}
  Let $\vp \in \HH'_G(V)$.  We verify that $\vp(g)$ induces an endomorphism of $\vnk$ for $g \in \P Z_M^- \P$ by checking it
  for $g = h \in Z_M^-$.  Take $n \in N(\O)$ and note that $n\vp(h) = \vp(h) n^h$, where $n^h \in N(\O)$ since $h^{-1}$ contracts
  $N(\O) = \P^+$ (Lemma~\ref{lm:iwahori-decomp}).  Thus $i^\P(\vp)$ is well defined.
  
  We now use Lemma~\ref{lm:hecke-compat-m-p} and its proof. Note that $\HH'_G(V)$ is spanned by the $T_\lambda$ with
  $\lambda(\varpi) \in Z_M^-$. Fix such a~$\lambda$ and let $h := \lambda(\varpi)$. By definition, $T_\lambda(h) \in \End(V)$
  induces the identity map on~$V^{N_{-\lambda}(k)}$. But $N_{-\lambda} \subset N$ since $\lambda(\varpi) \in Z_M^-$, so 
  $T_\lambda(h)$ induces the identity map on $\vnk$. As a consequence, $i^\P(T_\lambda) = E_h$. So $T_\lambda$ in $\HH'_G(V)$
  is mapped to $E_h|_M = T_h^M = T_\lambda^M$ in $\HH'_M(\o V)$.  By Cor.~\ref{cor:support}, the composite map $\HH'_G(V) \to
  \HH'_M(\o V)$ is given by~$\smg$. As the $T_\lambda^M$ with $\lambda(\varpi) \in Z_M^-$ span $\HH'_M(\o V)$, the map is
  surjective. But $\smg$ is injective on~$\hg$, so $\HH'_G(V)$ is the inverse image under~$\smg$ of~$\HH'_M(\o V)$, in
  particular it is a subalgebra. Moreover $i^\P$ is an algebra homomorphism as $\smg$ is.

  The algebra isomorphism $i_\oP$ is obtained by duality from $i^\P$ (interchanging positive and negative roots and $V$ with
  $V^*$).
  
   The two identifications coincide since 
   for $m \in M(\O) Z_M^-$, the endomorphisms of $\vnk$ and $\vonk$ induced by $\vp(m)$ are identified
   via $\vnk \congto \vonk$.
\end{proof}

After this subsection we will no longer use the primed notation.  If $\o V$ is a weight for $M$, we will denote by $\HH(\o
V)$ the isomorphic Hecke algebras of Lemma~\ref{lm:hecke-compat-m-p}.  Thus $\HH(\o V)$ can be thought of as subalgebra of
$\HH_M(\o V)$, $\HH_\P(\o V)$, and $\HH_\oP(\o V)$.  Similarly, suppose that $V$ is an $M$-regular weight for $G$. Then $\vnk
\congto \vonk$ is a weight for $M$ and we denote it by $\o V$.  We will denote by $\HH(V)$ the isomorphic Hecke algebras of
Lemma~\ref{lm:i_P}. Thus $\HH(V)$ can be thought of as subalgebra of $\HH_G(V)$, $\HH_M(\o V)$, $\HH_\P(\o V)$, and
$\HH_\oP(\o V)$.  Note that $\HH(V) \cong \HH(\o V)$. There can be no confusion since $M$-regular weights for $G$
are in natural bijection with weights for $M$ (by Lemma~\ref{lm:lift-to-m-reg-wt}).

We will even write $\HH$ when it is clear from the context what $V$ (or $\o V$) is.

In the following proposition we will use natural maps $\ind_\P^G \vnk \to \ind_K^G V \to \ind_\oP^G V_{\o N(k)}$.
They are obtained by Frobenius reciprocity from $\vnk \to V \subset \ind_K^G V$, respectively
$V \to \Ind_\oP^K \vonk \subset \ind_\oP^G V_{\o N(k)}$. Alternatively, as Hecke operators they are supported on $K$
and map the identity of $K$ to the natural map $\vnk \to V$, respectively $V \to \vonk$.

\begin{prop}\label{prop:comp-betw-hecke}
  Assume that $V$ is $M$-regular. We have the following diagram of Hecke operators.
  \ifkuvio
  \[ \cellwidth=45pt \Diag
  \ind_\P^G \vnk & \rEpi & \ind_K^G V & \rInto & \ind_\oP^G V_{\o N(k)} \\
  \dTo_{i^{\P}(\vp)} &\ldDashto & \dTo^{\vp} \dy{-3mm} &\ldDashto & \dTo^{i_{\oP}(\vp)} \dy{-3mm} \\
  \ind_\P^G \vnk & \rEpi & \ind_K^G V & \rInto & \ind_\oP^G V_{\o N(k)} \\
  \endDiag \]
  \fi
  For all $\vp \in \HH \subset \hg$ the two squares commutes. If moreover $\supp(\vp) \subset K Z_M^{--} K$ then there are diagonal
  arrows making the whole diagram commute.
\end{prop}

\begin{proof}
First note that the horizontal map on the top left is surjective: this is because the map is
obtained by compactly inducing the $K$-linear surjection $\ind_\P^K \vnk \to V$ from $K$ to $G$.
Similarly, the horizontal map on the top right is injective.

To check commutativity, we first deal with the left half of the diagram. 
Without loss of generality, $\vp$ is supported on a single double coset $K h K$ with $h = \lambda(\varpi) \in Z_M^-$.  By
Frobenius reciprocity it suffices to check that the two maps around the left square agree on $\vnk \subset \ind_\P^G \vnk$.
A vector $v \in \vnk$ is mapped to $[1,v] \in \ind_K^G V$ which in turn is mapped to $\sum_{K\bs K h K} [g^{-1}, \vp(g)v] \in
\ind_K^G V$ under $\vp$.  If $\vp(g)v \ne 0$ then $g \in Kh\P$ by Cor.~\ref{cor:support}.  We verify that the natural map
$\P\bs \P h\P \to K\bs Kh\P$ is a bijection. It is enough to show that $\P \cap K^{h} = \P \cap \P^{h}$ or equivalently that
${}^{h} \P \cap K = {}^{h} \P \cap \P$. Now note that ${}^{h} \P \cap K = ({}^{h} \P^-) \P^0 \P^+ \subset \P$ by
Lemma~\ref{lm:iwahori-decomp}.

Thus we see that
\begin{equation}\label{eq:1}
  \sum_{K\bs KhK} [g^{-1}, \vp(g)v] = \sum_{\P\bs \P h\P} [g^{-1}, i^\P(\vp)(g)v],
\end{equation}
which shows that the left square commutes.

Now suppose that $\supp(\vp) = KhK$ for some $h = \lambda(\varpi) \in Z_M^{--}$. In this case $P_{-\lambda} = P$. We will
define the diagonal map $\vp_\delta \in \HH_{K,\P}(V,\vnk)$ as follows (see~\S\ref{sec:cpt-ind-generalities} for the notation). We let $\vp_\delta$
agree with $\vp$ on $\P hK$ and we let it vanish outside.  To see that it is well defined, note that $\vp(h)$ maps $V$ to
$V^{N_{-\lambda}(k)} = \vnk$ and that $\P$ preserves $\vnk$.  We check that the top left triangle commutes. For $v \in
\vnk$ the arrow to the right maps it to $[1,v]$ as before and $\vp_\delta$ further maps it to $\sum_{\P\bs \P h K} [g^{-1}, \vp(g)v]$.
As noted above, if $\vp(g)v \ne 0$ then $g \in Kh\P$. But $\P hK \cap Kh\P = \P h\P$ since $K \cap {}^h K \subset \P_{-\lambda} = \P$
by Prop.~\ref{prop:buildings-lemma}, so we are done as in~\eqref{eq:1}. Since the map at the top is surjective, the bottom triangle also commutes.

We now dualise the left half of the diagram, in the following sense. We think of maps between compact inductions
as Hecke operators and then apply the duality $\vp \mapsto \vp'$ considered in \S\ref{sec:satake-variants}. (Strictly
speaking, we also use this duality for maps between different compact inductions as in \S\ref{sec:cpt-ind-generalities}.)
We obtain
\ifkuvio
\[ \cellwidth=45pt \Diag
\ind_K^G \vd & \rTo & \ind_\P^G (\vd)_{N(k)}  \\
\dTo_{\vp'} &\ldDashto & \dTo^{i^{\P}(\vp)'} \dy{-3mm} \\
\ind_K^G \vd & \rTo & \ind_\P^G (\vd)_{N(k)} \\
\endDiag \]
\fi
Since the natural maps $\vnk \to V$ and $V^* \to (V^*)_{N(k)}$ are dual, the top and bottom maps are the natural ones.  By
construction of $i_\P$, $i^\P(\vp)' = i_\P(\vp')$. By replacing $V$ by $V^*$ and
by interchanging positive and negative roots, we obtain the right half of the diagram.
\end{proof}

\begin{coroll}\label{cor:parahoric-and-spherical-inductions}
  Suppose that $V$ is $M$-regular and that $\chi : \HH_M(\vnk) \to \k$ is an algebra homomorphism.
  Then the maps of Prop.~\ref{prop:comp-betw-hecke} induce isomorphisms
  \begin{equation*}
    \ind_\P^G \vnk \otimes_{\HH,\chi} \k \congto \ind_K^G V  \otimes_{\HH,\chi} \k\congto \ind_\oP^G V_{\o N(k)} \otimes_{\HH,\chi} \k,
  \end{equation*}
  where $\HH = \HH(V)$ as above.
\end{coroll}

\begin{proof}
  Pick any $\lambda \in X_*(T)_-$ such that $h := \lambda(\varpi) \in Z_M^{--}$. Since $T^M_h T^M_{h^{-1}} = 1$, we have
  $\chi(T_h^M) \ne 0$. Let $\mm = \ker(\chi)$, an ideal of $\HH$.  We will use that $\HH$ is commutative.
  
  The first map is clearly surjective. Let $\vp = T_\lambda \in \HH$; it is identified with $T_h^M$ so $\chi(\vp) \ne 0$.
  Recall that we denoted the diagonal map corresponding to $\vp$ on the left-hand side of the diagram in
  Prop.~\ref{prop:comp-betw-hecke} by $\vp_\delta$. We now show that $\vp_\delta$ is $\HH$-linear.  Let $\xi : \ind_\P^G \vnk
  \to \ind_K^G V$ denote the natural map, and suppose that $\psi \in \HH \subset \hg$. Then
  \begin{equation*}
    \vp_\delta \circ \psi \circ \xi = i^\P(\vp) \circ i^\P(\psi) = i^\P(\psi) \circ i^\P(\vp) = i^\P(\psi) \circ \vp_\delta \circ \xi,
  \end{equation*}
  and the claim follows, since $\xi$ is surjective.
  
  If $f \in \pind \vnk$ and $\xi(f) \otimes 1 = 0$ then $\xi(f) \in \mm (\kind V)$. Then
  \begin{equation*}
    \vp(f) = \vp_\delta(\xi(f)) \in \vp_\delta(\mm (\kind V)) = \mm\, \vp_\delta(\kind V),
  \end{equation*}
  so $\vp(f) \otimes 1 = 0$. Thus $f \otimes 1 = \chi(\vp)^{-1}(\vp(f) \otimes 1) = 0$, so the first map is injective.

  Suppose $y \otimes 1 \in \opind\vonk \otimes_{\HH,\chi} \k$. Then $y \otimes 1 = \chi(\vp)^{-1} (\vp(y) \otimes 1)$ comes
  from $\kind V$, by the commuting triangle on the bottom right. This proves that the second map is surjective.

  Suppose $x \otimes 1 \in \kind V \otimes_{\HH,\chi} \k$ maps to zero. Let $\eta : \ind_K^G V \to \ind_\oP^G V_{\o N(k)}$ denote
  the natural map. Then $\eta(x) \in \mm (\opind \vonk)$, so 
  \begin{equation*}
    \eta\vp(x) = \vp \eta(x) \in \mm (\vp \opind \vonk) \subset \mm\, \eta(\kind V). 
  \end{equation*}
  As $\eta$ is injective, $\vp(x) \in \mm(\kind V)$. Finally $x \otimes 1 = \chi(\vp)^{-1} (\vp(x) \otimes 1) = 0$.
  This proves that the second map is injective.
\end{proof}

\section{Parabolic inductions and compact inductions}
\label{sec:parab-ind-comp-ind}

The following theorem is inspired by work of Barthel--Livn\'e for $\GL_2$ (see \cite[Thm.~25]{bib:BL-general}). 
One aspect of the proof, namely the comparison of parahoric and parabolic inductions, crucially
use ideas of Schneider--Stuhler~\cite{bib:schneider-stuhler} and Vign\'eras~\cite{bib:vigneras2}. See also the comment
after Cor.~\ref{cor:parab-ind-cpt-ind}.

\begin{thm}\label{thm:cptind-parabind}
  Let $P = M N$ be a standard parabolic in $G$ and suppose that $V$ is an $M$-regular weight for $G$.
  Then for any algebra homomorphism $\chi : \HH_{M} (V_{\o N(k)}) \to \bar k$, there is a natural isomorphism
  of smooth $G$-representations,
  \begin{equation*}
    \ind_K^G V \otimes_{\hg,\chi} \k \congto
    \Ind_{\o P}^G \Big\{ \ind_{M(\O)}^M V_{\o N(k)} \otimes_{\HH_{M} (V_{\o N(k)}),\chi} \k \Big\}.
  \end{equation*}
\end{thm}

Note that $\chi$ becomes a character of $\hg$ by composing with the partial Satake homomorphism $\psmg : \hg \to \hm$.

\begin{proof} We begin by defining maps
\begin{equation}\label{eq:4}
  \ind_K^G V \xrightarrow{\;\eta\;} \ind_\oP^G \vonk \xrightarrow{\;\zeta\;} \Ind_{\o P}^G (\ind_{M(\O)}^M \vonk).
\end{equation}

The map $\eta$ is the natural one that we already used in Prop.~\ref{prop:comp-betw-hecke}. Explicitly,
$\eta$ sends $f \in \ind_K^G V$ to $g \mapsto p_{\o N}(f(g))$ in $\ind_\oP^G \vonk$.

We define $\zeta$ by the following formula, for $f \in \opind \vonk$:
\begin{equation*}
  \zeta(f)(g) = \sum_{\o P(\O)\bs \o P} \o p^{-1} [1, f(\o p g)].
\end{equation*}
Note that $\o P$ acts via $M$.
By using that $\o P(\O) = M(\O) \o N(\O)$ we see that the term in the sum only depends on $\o p \in \o P(\O)\bs \o P$.  We
check that the sum only involves finitely non-zero terms. Without loss of generality $f$ is supported on a single coset $\oP
\gamma$. Since $\o P \cap \oP = \o P (\O)$, we see that the sum involves at most one term. Clearly $\zeta(f)$ is $\o
P$-equivariant and $\zeta$ is $G$-equivariant. It follows that $\zeta(f)$ is smooth, so $\zeta(f) \in \opInd (\mind \vonk)$.
Therefore $\zeta$ is well defined.

Let $\HH = \HH(V)$ as in \S\ref{sec:comp-betw-hecke}. It acts on the first term of~\eqref{eq:4} via $\HH_G(V)$, on the second
term via $\HH_\oP(\vonk)$ and on the third term via $\HH_M(\vonk)$.

\emph{Step 1.} Check that $\eta$ and $\zeta$ are $\HH$-equivariant. We already checked this for $\eta$ in Prop.~\ref{prop:comp-betw-hecke}.

We introduce a useful shorthand for certain functions in $\opInd \sigma$ for any smooth $M$-representation $\sigma$. For any
compact open subset $\Omega \subset N$ and any $x \in \sigma$ we write $[\Omega, x]$ for the function that is supported on
$\o P \Omega^{-1}$ and sends all $\nu \in \Omega^{-1}$ to $x$. This function is locally constant, since $N \into \o P \bs G$
has open image and $\Omega$ is compact open. Note that for $m \in M$ and $n \in N$,
\begin{equation*}
  m [\Omega,x] = [{}^m \Omega,mx],\quad n [\Omega,x] = [n \Omega,x].
\end{equation*}
We claim that $\zeta([1,\o v]) = [\oP^+,[1,\o v]]$. Clearly the left-hand side is supported on $\o P\oP = \o P\oP^+$. It is then
an easy computation to check that both sides agree on $\oP^+$.

To check that $\zeta$ is $\HH$-equivariant, we can immediately reduce to the case when $\vp = T_\lambda$ for some $\lambda
\in X_*(T)_-$ with $h := \lambda(\varpi) \in Z_M^-$. Then $\psmg(T_\lambda) = T_h^M$ by Cor.~\ref{cor:support}. Suppose $\o v \in
\vonk$. On the one hand, in $\ind_{M(\O)}^M \vonk$ we have $T_h^M([1,\o v]) = [h^{-1},\o v]$ by~\eqref{eq:3}. Thus
$(\psmg(\vp)\circ \zeta)([1,\o v]) = [\oP^+,[h^{-1},\o v]]$. On the other hand,
\begin{align*}
  (\zeta \circ i_\oP(\vp))([1,\o v]) &= \sum_{\oP\bs\oP h\oP} g^{-1} [\oP^+, [1,\vp(g)\o v]] \\
  &= \sum_{(\oP \cap \oP^h)\bs\oP} \o\p^{-1}h^{-1} [\oP^+, [1,T_\lambda(h)\o\p\o v]].
\end{align*}
Since $(\oP^+)^h \bs \oP^+ \to (\oP \cap \oP^h)\bs\oP$ is a bijection by Lemma~\ref{lm:iwahori-decomp} and taking into account
that $\oP^+$ fixes $\o v$ and that $T_\lambda(h)$ is trivial on $\vonk$ (see the proof of Lemma~\ref{lm:i_P}), the sum above simplifies to
\begin{equation*}
  \sum_{(\oP^+)^h \bs \oP^+} [\o\p^{-1} (\oP^+)^h, [h^{-1},\o v]] = [\oP^+,[h^{-1},\o v]].
\end{equation*}

\emph{Step 2.}\mar{We cite this later! So don't change number\dots}
Check that $\zeta \circ \eta$ is $\hg$-equivariant (via partial Satake). Let
$\theta : V \to \Ind_{\o P}^G (\ind_{M(\O)}^M \vonk)$ be the composition of $V \to \kind V$ with $\zeta \circ \eta$.
It is $K$-linear. Note that by definition of $\zeta$ and since $\eta([1,v])$ is supported on $K$,
\begin{equation*}
  \theta(v)(1) = \zeta(\eta([1,v]))(1) = [1,\eta([1,v])(1)] = [1,p_{\o N}(v)].
\end{equation*}
Thus, in the notation of \eqref{eq:6}, $\o \theta \in \Hom_{M(\O)} (\vonk,\mind \vonk) = \hm$ is the natural map.
By Lemma~\ref{lm:hecke-parab-ind}, $\o{\theta \ast \vp} = \o\theta \ast \vp_M$, where $\vp_M = \psmg \vp$. Since
$\hm$ is commutative and since the map in \eqref{eq:6} is natural in $\sigma$, we get that
$\o{\theta \ast \vp} = \vp_M \circ \o\theta = \o{\vp_M \circ \theta}$. But the map in \eqref{eq:6} is injective, so
$\theta \ast \vp = \vp_M \circ \theta$, as required.

\emph{Step 3.}
Check that $\eta \otimes_{\HH,\chi} \k$ is an isomorphism.
This is the content of Cor.~\ref{cor:parahoric-and-spherical-inductions}.

\emph{Step 4.}
Check that $\zeta \otimes_{\HH,\chi} \k$ is surjective.

Fix any non-zero $\o v\in \vonk$.
We claim that $\zeta([1,\o v]) \otimes 1 = [\oP^+,[1,\o v]] \otimes 1$ generates $\opInd (\mind \vonk) \otimes_{\HH,\chi} \k$.
Pick $h_0 \in Z_M^{--}$ and let $x_0 = [1,\o v] \in \mind \vonk$. Since
\begin{equation*}
  [\oP^+,h_0^n x_0] \otimes 1 = [\oP^+,T_{h_0^{-n}}^M (x_0)] \otimes 1 = [\oP^+, x_0] \otimes \chi(T_{h_0^{-n}}^M),
\end{equation*}
it is enough to show that $[\oP^+,h_0^n x_0]$ ($n \in \Z$) generate $\opInd (\mind \vonk)$.

Note that the smooth $M$-representation $\sigma := \mind \vonk$ is generated by $x_0$. (This is the only property of $\sigma$ that we will use.)
We want to show that any $f \in \opInd \sigma$ is contained in the $G$-representation generated by the $[\oP^+,h_0^n x_0]$ for $n \in \Z$.
By writing $\o P\bs G$ as a (finite) disjoint union of compact open subsets, each of which is contained
in a $G$-translate of $\o P \bs \o P N$, we can reduce to the case that $\supp(f) \subset \o P N$. Since $f|_N$ is locally constant
and compactly supported, we can moreover assume that $f = [\Omega,x]$, for some compact open subset $\Omega \subset N$ and
$x \in \sigma$. We can write $x = \sum_i \lambda_i m_i x_0$ as finite linear combination with $\lambda_i \in \k$, $m_i \in M$, so
\begin{equation*}
  f = \sum_i \lambda_i m_i [\Omega^{m_i}, x_0].
\end{equation*}
Thus we may assume that $f = [\Omega,x_0]$, for some compact open $\Omega \subset N$. By Lemma~\ref{lm:iwahori-decomp} there is an $n \gg 0$
such that $\Omega = \coprod_j \nu_j \cdot (\oP^+)^{h_0^n}$, a finite disjoint union with $\nu_j \in N$. Therefore
\begin{equation*}
  f = \sum_j \nu_j [(\oP^+)^{h_0^n}, x_0] = \sum_j \nu_j h_0^{-n} [\oP^+, h_0^nx_0],
\end{equation*}
which completes the argument.

\emph{Step 5.}
Check that $\zeta \otimes_{\HH,\chi} \k$ is injective.

Let $f' \in \opind \vonk$ such that $\zeta(f') \otimes 1 = 0$. We need to show that $f' \otimes 1 = 0$.
Note that for $h \in Z_M^-$ we have $\supp(E_h f') \subset \oP h \supp(f')$ since $E_h$ is supported on $\oP h\oP$.
(Recall that $E_h \in \HH$ was defined in the proof of Lemma~\ref{lm:hecke-compat-m-p}. It is identified with $T_h^M$.)
Thus by Lemma~\ref{lm:shrink-support} there is an $h \in Z_M^-$ such that $\supp(E_h f') \subset \oP T^- K$. Since
$E_h f' \otimes 1 = \chi(E_h)(f' \otimes 1)$ and $\chi(E_h) = \chi(T_h^M)\ne 0$, we may assume that $\supp(f') \subset \oP T^- K$.

Lemma~\ref{lm:iwahori-decomp} shows that for $t \in T^-$, ${}^t (\oP^-) \subset \oP^-$ and $(\oP^+)^t \subset \oP^+$.  Thus
for $h \in Z_M^-$ and $t \in T^-$, $\oP h \oP t K = \oP (h \oP^- \oP^0) (\oP^+ t) K \subset \oP ht K$. We can write $f' =
\sum_{i=1}^r f'_i$ such that $\supp(f'_i) = \oP t'_i k_i$ ($t'_i \in T^-$, $k_i \in K$).  Since any element of $Z_M$ can be
written in the form $h' h^{-1}$ for $h'$, $h \in Z_M^-$, we can find $h_i \in Z_M^-$ and $\tau_i \in T(\O)$ such that $\tau_i h_i
t'_i = \tau_j h_j t'_j$ whenever $T(\O) Z_M t'_i = T(\O) Z_M t'_j$. Let $f_i = \chi(E_{h_i})^{-1} E_{h_i} f'_i$, $f =
\sum_i f_i$, and $t_i = \tau_i h_i t_i' \in T^-$.  Then $f' \otimes 1 = f \otimes 1$, moreover $\supp(f_i) \subset \oP t_i K$ 
(as $T(\O) \subset \oP$) and
\begin{equation}
  \label{eq:7}
  \text{$t_i = t_j$ whenever $T(\O) Z_M t_i = T(\O) Z_M t_j$.}
\end{equation}
We can now write $f_i = \sum f_{ij}$ (finite sum) such that $\supp (f_{ij}) = \oP t_i k_{ij}$ for some $k_{ij} \in K$.
By combining those $f_{ij}$ that have identical support, we may assume moreover:
\begin{equation}
  \label{eq:9}
  \text{the sets $\oP t_i k_{ij}$ are pairwise disjoint.}
\end{equation}

We now show that $f = 0$. In fact we show that $f|_{\oP \o P} = 0$, but all our conditions on $f$ are invariant under $K$,
so $f = 0$ since $\oP \o P K = G$. (When we replace $f$ by $kf$ for some $k \in K$, then $t_i$ is unchanged and $k_{ij}$ is replaced by $k_{ij}k^{-1}$.
Also note that $\zeta(kf) \otimes 1 = k\zeta(f) \otimes 1 = 0$.)

All we will use is that the image of $\zeta(f)(1)$ in $\mind \vonk \otimes_{\HH, \chi} \k$ is zero. Let us show that this
latter space is naturally isomorphic to $\mzind \vonk$. We let $Z_M$ act on $\vonk$ by declaring that $h \o v =
\chi(T_{h^{-1}}^M) \o v$ for $h \in Z_M$ and $\o v \in \vonk$.  This is compatible with the $M(\O)$-action: for $h_0 \in
Z_M(\O)$, $T_{h_0^{-1}}^M = \omega(h_0) T_1^M = \omega(h_0)$ by~\eqref{eq:2} and~(\ref{eq:5}), where $\omega : Z_M(k) \to
k\s$ is the central character of $\vonk$, so $\chi(T_{h_0^{-1}}^M) = \omega(h_0)$. The map
\begin{equation*}
  \mind \vonk \to \mzind \vonk
\end{equation*}
is the obvious one induced by Frobenius reciprocity. It sends $[m,\o v]$ to $[m,\o v]_{Z_M}$, where the subscript is used to distinguish between the
two induced representations. In particular it is surjective. For $h \in Z_M$ it sends $T_h^M[m,\o v] = h^{-1}[m,\o v]$
to $[h^{-1}m,\o v]_{Z_M} = [m,h^{-1} \o v]_{Z_M} = \chi(T_h^M) [m,\o v]_{Z_M}$. (See~\eqref{eq:3} for the first identity.)
The induced map
\begin{equation}\label{eq:18}
  \mind \vonk \otimes_{\HH,\chi} \k \onto \mzind \vonk
\end{equation}
is injective: we can lift any element in the kernel to one of the form $\sum [m_i,\o v_i] \in \mind \vonk$, where the $m_i$ lie in distinct
$M(\O)Z_M$-cosets. By considering its image on the right-hand side we see that all $\o v_i$ are zero.

Since $\zeta(f)(1) \otimes 1$ depends only on $f|_{\oP\o P}$, we will assume from now on that $f$ is supported on $\oP \o P =
\oP^+ \o P$.  Thus $k_{ij} \in (\oP^+)^{t_i} \cdot \o P \cap K = (\oP^+)^{t_i} \cdot \o P(\O)$ for all $i$. Thus we may assume
that $k_{ij} \in \o P(\O) = \oP^- \oP^0$ without changing $\oP t_i k_{ij}$. Since $t_i \in T^-$ shrinks $\oP^-$ we may even assume
that $k_{ij} \in \oP^0 = M(\O)$.  Let us write $f_{ij} = [k_{ij}^{-1}t_i^{-1},\o v_{ij}]$. As $k_{ij}^{-1}t_i^{-1} \in M$, we see that
$\zeta(f_{ij})(1) \otimes 1 = [k_{ij}^{-1}t_i^{-1},\o v_{ij}]_{Z_M}$.  By showing that $M(\O)Z_M t_i k_{ij}$ are pairwise disjoint we
will deduce from $\zeta(f)(1) \otimes 1 = 0$ that $\o v_{ij} = 0$ for all $i$, which will complete the argument. Thus suppose that
$M(\O)h t_i k_{ij} = M(\O) t_r k_{rs}$ for some $h \in Z_M$. Note that $Z_M T^-$ is the antidominant part of~$T$ with respect
to~$M$.  Thus by the Cartan decomposition for $M$ we have $T(\O) h t_i = T(\O) t_r$, so by~\eqref{eq:7} it follows that $t_i
= t_r$ and hence that $h \in Z_M(\O)$. Therefore $M(\O)t_i k_{ij} = M(\O)t_r k_{rs}$ which implies $\oP t_i k_{ij} = \oP t_r k_{rs}$, so
by~\eqref{eq:9} we see that $(i,j) = (r,s)$ and we are done.

\emph{Step 6.}
We showed that $(\zeta \circ \eta) \otimes_{\HH,\chi} \k$ is an $\hg$-linear isomorphism (as $\hg$ commutative).
By tensoring over $\hg$ with $\chi$ we see that $(\zeta \circ \eta) \otimes_{\hg,\chi} \k$ is an isomorphism.
The target of that isomorphism is isomorphic to $\Ind_{\o P}^G (\ind_{M(\O)}^M \vonk) \otimes_{\hm,\chi} \k$,
since $\psmg : \hg \INTO \hm$ is a localisation map.
Since $\hm$ is noetherian, we can pick a finite generating set $\vp_1$, \dots, $\vp_n$ of the ideal $\ker \chi$.
We have an exact sequence
\begin{equation*}
  \sigma^{\oplus n} \xrightarrow{\,\sum \vp_i\,} \sigma \to \sigma \otimes_{\hm,\chi} \k \to 0.
\end{equation*}
As $\Ind_{\o P}^G$ is exact, we deduce that $\Ind_{\o P}^G (\ind_{M(\O)}^M \vonk) \otimes_{\hm,\chi} \k$ is isomorphic to 
$\Ind_{\o P}^G (\ind_{M(\O)}^M V_{\o N(k)} \otimes_{\hm,\chi} \k)$.
This completes the proof.
\end{proof}

We record the following corollary to the proof.

\begin{coroll}\label{cor:parab-ind-cpt-ind}
  Suppose that $\o V$ is a weight for $M$ and that $\chi_M : Z_M \to \k\s$ such that $\chi_M|_{Z_M(\O)}$ is the central
  character of $\o V$. Then there is an algebra homomorphism $\chi : \HH \to \k$ such that $\chi(T_h^M) = \chi_M(h)^{-1}$ 
  for all $h \in Z_M^-$ and we have
  \begin{equation*}
    \opind \o V\otimes_{\HH,\chi} \k \congto \opInd (\mzind \o V).
  \end{equation*}
\end{coroll}

On the right-hand side we let $h \in Z_M$ act on $\o V$ by $\chi_M(h)$.

\begin{proof}
  This isomorphism was obtained in the above proof, in case $\o V$ is of the form $\vonk$, where $V$ is $M$-regular, and $\chi$ is the restriction
  of an algebra homomorphism $\HH_M(\vonk) \to \k$. But it did not matter that $\o V$ was of that form. (It is anyway, by Lemma~\ref{lm:lift-to-m-reg-wt}.)
  Moreover, by Prop.~\ref{prop:param-hecke-evals} and
  Cor.~\ref{cor:param-hecke-evals} (or directly) the pair $(M,\chi_M)$ gives rise to $\chi : \HH_M(\o V) \to \k$ such that
  $\chi(T_h^M) = \chi_M(h)^{-1}$ for $h \in Z_M^-$. 
\end{proof}

If $G = \GL_n$ and $P = B$ is the Borel, this recovers the results of Schneider--Stuhler \cite[Prop.~11]{bib:schneider-stuhler} ($\chi = 1$)
and Vign\'eras \cite[Thm.~4.10]{bib:vigneras2} ($\chi$ arbitrary). In those cases the right-hand side simplifies
to the principal series $\obInd (\chi_T)$.

\section{Hecke eigenvalues and supersingularity}
\label{sec:hecke-eigenv}

The following proposition allows one to compare Hecke eigenvalues between different weights. It is analogous to the classical
result parameterising unramified Hecke eigenvalues by unramified characters of the torus \cite[Cor.~4.2]{bib:Cartier}.
It will be convenient to define a \emph{standard Levi} to be the unique Levi subgroup containing $T$ of a standard parabolic subgroup.

\begin{prop}\label{prop:param-hecke-evals}
  There is a natural bijection between $\chi \in \Hom\kalg(\hg,\k)$ and pairs $(M,\chi_M)$, where $M$ is a standard Levi and
  $\chi_M : Z_M \to \k\s$ a character such that $\chi_M|_{Z_M(\O)}$ is the central character of $\vonk$.

  Given such a pair $(M,\chi_M)$
  the corresponding set of eigenvalues is the composite of algebra homomorphisms
  \begin{equation*}
     \chi: \hg \xrightarrow{\psg}\HH_T^-(\vouk) \xrightarrow{\chi'} \k,
  \end{equation*}
  where $\chi'(\vp) = \sum_{Z_M(\O)\bs Z_M} \vp(z) \chi_M(z)^{-1}$.
\end{prop}

We will say in the following that $\chi \in \Hom\kalg(\hg,\k)$ \emph{is parameterised by the pair $(M,\chi_M)$} if they
correspond under the bijection in the proposition. Suppose now that $(M,\chi_M)$ consists of a standard Levi $M$ and
an arbitrary smooth character $\chi_M : Z_M \to \k\s$. Then there may be more than one weight $V$ such that
the central character of $\vonk$ equals $\chi_M|_{Z_M(\O)}$; for each one we obtain a corresponding 
algebra homomorphism $\hg \to \k$. In \S\ref{sec:classification} we will see that if $G = \GL_n$ and $\pi$ is an irreducible admissible
$G$-representation, then all Hecke eigenvalues in all weights of $\pi$ are identified in this manner.

We will see below that $M$ is in fact the smallest standard Levi such that $\chi$ factors through $\psmg : \hg \to \hm$ (as an algebra
homomorphism).
We have the following immediate consequence.

\begin{coroll}\label{cor:param-hecke-evals}
  In the situation of Prop.~\ref{prop:param-hecke-evals}, we have for $\lambda \in X_*(T)_-$,
  \begin{equation*}
    \chi'(\tau_\lambda) = 
    \begin{cases}
      \chi_M(\lambda(\varpi))^{-1} & \text{if $\lambda(\varpi) \in Z_M$,} \\
      0 & \text{otherwise.}
    \end{cases}
  \end{equation*}
\end{coroll}

Given $\chi \in \Hom\kalg(\hg,\k)$ we can consider its Satake transform $\chi' \in \Hom\kalg(\HH_T^-(\vouk),\k)$.
We will say that $\chi'$ vanishes on an open subset $X \subset T^-$ if $\chi'$ vanishes on all elements of
$\HH_T^-(\vouk)$ that are supported on $X$. Then $\chi'$ has a well-defined \emph{support}, namely the complement in $T^-$ of the biggest open
subset on which it vanishes. It has the following alternative description: under the isomorphism $T^-/T(\O) \congto X_*(T)_-$ it corresponds
to a subset of $X_*(T)_-$. And this is precisely the subset of $X_*(T)_-$ on which $\chi' : \k [X_*(T)_-] \to \k$ is non-vanishing.

\begin{lm}\label{lm:factor-through-hm}
  Suppose $P = MN$ is a standard parabolic.
  A $\k$-algebra homomorphism $\chi : \hg \to \k$ factors through $\psmg$ if and only if $\supp \chi' \supset Z_M^- T(\O)$.
\end{lm}

\begin{proof}
  As we noted in the proof of Prop.~\ref{prop:partial-satake}, $\hm$ is the localisation of $\hg$ at any element
  $\vp$ such that the support of $\psg(\vp)$ is a single $T(\O)$-coset in $Z_M^{--}$. If $\supp \chi' \supset Z_M^- T(\O)$ then
  $\chi(\vp) = \chi'(\psg(\vp)) \ne 0$ and so $\chi$ factors through $\hm$.
  
  Conversely, if $\chi$ factors through $\psmg$, then $\chi'$ extends to the subalgebra of $\HH_T(\vouk)$ consisting of those elements whose
  support in $T$ is antidominant for $M$. Since $Z_M T(\O)$ is a subgroup of $T$ that is contained in the part of $T$ that is antidominant for $M$,
  it follows that $\supp \chi' \supset Z_M^- T(\O)$. 
\end{proof}

\begin{proof}[Proof of Proposition~\ref{prop:param-hecke-evals}]
  By the proof of Cor.~1.5 in \cite{bib:satake}, $\supp \chi'$ is of the form $Z_{M}^- T(\O)$ for some
  standard parabolic $P = MN$. We will show that algebra homomorphisms $\chi$ with $\supp \chi' = Z_M^- T(\O)$ biject with
  characters $\chi_M : Z_M \to \k\s$ such that $\chi_M|_{Z_M(\O)}$ is the central character of $\vonk$.
  By Lemma~\ref{lm:factor-through-hm}, $\chi$ factors through an algebra homomorphism $\wt\chi:\hm\to\k$.
  It is not hard to verify that $\wt\chi'$ has support $Z_M T(\O)$.
  Thus by replacing $(M,\vonk,\wt\chi)$ by $(G,V,\chi)$, we are reduced to the case $\supp \chi' = Z T(\O)$.

  Let us write $\HH_T^- (\vouk) = \HH \oplus \I$, where $\HH$ (resp., $\I$) consists of those elements $\vp$ whose
  support is contained in $ZT(\O)$ (resp., disjoint from $ZT(\O)$). Clearly $\HH$ is a subalgebra and $\I$ is an ideal,
  so the $\chi'$ with support $ZT(\O)$ biject with $\Hom\kalg(\HH,\k)$. By restricting functions to $Z$, $\HH$ is isomorphic to
  \begin{equation*}
    \{ \vp : Z \to \k : \vp(z_0 z) = \omega(z_0) \vp(z)\ \forall z_0 \in Z(\O), z \in Z;\ \supp \vp \text{\ cpt.} \}
  \end{equation*}
  (as algebra under convolution), where $\omega$ is the central character of $V$. The linear dual of this space consists of all
  functions $f : Z \to \k$ such that $f(z_0z) = \omega(z_0)^{-1} f(z)$ under the pairing $\langle \vp, f\rangle =
  \sum_{Z(\O)\bs Z} \vp(z) f(z)$. A simple argument shows that the linear map $\HH \to \k$ induced by $f$ is an algebra
  homomorphism if and only if $f$ is a homomorphism $Z \to \k\s$. Finally we take the inverse of this homomorphism.
\end{proof}

\begin{lm}\label{lm:hecke-evals-supersing}
  Suppose that $\pi$ is a smooth $G$-representation that has a central character~$\omega_\pi$. Suppose that $V$ is a weight and that $\chi : \hg \to \k$ is a
  set of Hecke eigenvalues on $\Hom_K(V,\pi)$. If $\chi$ is parameterised by the pair $(M,\chi_M)$, then $\chi_M|_Z = \omega_\pi$.
\end{lm}

\begin{proof}
  It is enough to show that $\chi_M(z) = \omega_\pi(z)$ for all $z$ of the form $\lambda(\varpi) \in Z$. First, note that
  $\psg(T_\lambda) = \tau_\lambda$: this follows, for example, directly from the definition of $\psg$ or from~\eqref{eq:10}. By
  Cor.~\ref{cor:param-hecke-evals}, $\chi(T_\lambda) = \chi_G(z)^{-1}$. Finally, consider the induced $G$-linear map $\kind V
  \to \pi$. It follows from~\eqref{eq:3} that $\chi(T_\lambda) = \omega_\pi(z)^{-1}$.
\end{proof}

\begin{lm}\label{lm:hecke-evals-parab-ind}
  Suppose that $P = MN$ is a standard parabolic and that $\sigma$ is an admissible $M$-representation.
  Suppose that $V$ is a weight for $G$. Then the Hecke eigenvalues of $V$ in $\opInd \sigma$ and of $\vonk$ in $\sigma$ are parameterised
  by the same set of pairs $(L,\chi_L)$. In particular, $L \subset M$ in each case.
\end{lm}

\begin{proof}
  By Lemma~\ref{lm:hecke-parab-ind} the possible $\hg$-eigenvalues in $\Hom_K(V,\opInd \sigma)$ are obtained from the
  possible $\hm$-eigenvalues in $\Hom_{M(\O)} (\vonk,\sigma)$ by composing with $\psmg$. By construction,
  the pair associated to $\chi : \hm \to \k$ is the same as the pair associated to $\chi \circ \psmg$.
\end{proof}

\begin{lm}\label{lm:hecke-evals-twist}
  Suppose that $\pi$ is a smooth $G$-representation and that $f \in \Hom_K(V, \pi)$ is an $\hg$-eigenvector. Suppose that
  $\eta : G \to \k\s$ is a smooth character.
  If the Hecke eigenvalues of $f$ are parameterised by $(M,\chi_M)$, then the Hecke eigenvalues of $f \otimes \eta$
  are parameterised by $(M,\chi_M \eta|_{Z_M})$.
\end{lm}

\begin{proof}
  Since $\kind (V\otimes \eta) \cong \kind V \otimes \eta$, we have a natural isomorphism $\HH_G(V) \congto \HH_G(V \otimes
  \eta)$, $\vp \mapsto \vp_\eta$, under which the Hecke eigenvalues of $f$ and $f \otimes \eta$ are identified. A calculation
  shows that $\vp_\eta(g) = \eta(g)\vp(g)$ when we canonically identify $\End V$ and $\End(V \otimes \eta)$, which implies
  that $\psg(\vp_\eta)(t) = \eta(t)\psg(\vp)(t)$. (Note that $\eta$ is trivial on $\o U$, as it kills all pro-$p$
  subgroups.) The claim now follows from the formula for $\chi'$ in Prop.~\ref{prop:param-hecke-evals}.
\end{proof}

Finally we define what it means for an irreducible admissible representation $\pi$ to be supersingular. Roughly speaking, 
all Hecke eigenvalues occurring in~$\pi$ should be as trivial as possible. We fixed a hyperspecial maximal compact subgroup $K$
of $G(F)$ at the beginning, but our definition should be independent of this choice and we should take into account all of them.
(We thank M.-F.~Vign\'eras for this observation.) When $G = \GL_n$ all hyperspecial maximal compact subgroups of $G(F)$ are
conjugate and we can just work with our fixed choice of $K$.

Suppose the reductive group scheme $G'_{/\O}$ is another reductive integral structure of $G_{/F}$. Fix a maximal split torus
$T'$ of $G'$ and a Borel subgroup $B'$ containing $T'$. Let $K' = G'(\O)$. We claim that if $K$ and $K'$ are $G(F)$-conjugate,
then even $(K,T,B)$ and $(K',T',B')$ are $G(F)$-conjugate.
Without loss of generality we may assume that $K = K'$. By \cite[3.8.1]{bib:Tits}, \cite[II.5.1.40]{bib:BT2} reductive integral
structures of $G_{/F}$ are naturally in bijection with hyperspecial points of the reduced building. In this bijection
$G$ corresponds to the unique fixed point of $G(\O)$ in the reduced building. Thus 
$G = G'$.  Now by the same argument that was used towards the end of the proof of Lemma~\ref{lm:iwahori-decomp} there is an
element of $G(\O)$ that conjugates $(T,B)$ to $(T',B')$.

A \emph{$K'$-weight} is an irreducible representation~$V'$ of~$G'(k)$, or equivalently of~$K'$. (Hence a $K$-weight is a
weight.)  We denote the corresponding Hecke algebra by $\HH_{G,K'}(V')$ to emphasise the dependence on~$K'$.

\begin{df}\label{df:supersingular}
  Let $\pi$ be an irreducible admissible representation. We say that $\pi$ is \emph{supersingular} if for all triples $(K',
  T', B')$ as above, for all $K'$-weights $V'$ and for all Hecke eigenvalues $\chi'$ on $\Hom_{K'}(V',\pi)$ the following
  equivalent conditions hold.
  \begin{enumerate}
  \item $\chi'$ is parameterised by the pair $(G,\omega_\pi)$.
  \item $\chi'$ is parameterised by a pair $(G,\chi_G)$ for some $\chi_G : Z \to \k\s$.
  \item $\chi'$ does not factor through ${}'\SS^{M'}_G : \HH_{G,K'}(V') \to \HH_{M',M'(\O)}(V'_{\o{N'}(k)})$ for any proper parabolic $P'
    = M'N'$ containing $B'$.
  \end{enumerate}
\end{df}

Note that the first two conditions are equivalent by Lemma~\ref{lm:hecke-evals-supersing} and the last two conditions are equivalent
by Lemma~\ref{lm:factor-through-hm} and the proof of Prop.~\ref{prop:param-hecke-evals}.

It is easy to see that in this definition, the triple $(K',T',B')$ only matters up to $G(F)$-conjugacy.
Thus we only need to let $K'$ run through a set of representatives for the finitely many conjugacy classes of hyperspecial maximal compact subgroups in
$G(F)$ (and choose compatible $T'$, $B'$ for each).

Suppose now that $G = \GL_n$. By the above, the supersingularity condition has to be checked only for our fixed choice
$(K,T,B)$. Even better, as a corollary to our main results we will see that it has to be checked only for one weight~$V$ and for
one~$\chi$. (See Cor.~\ref{cor:uniqueness}.) Moreover supersingularity can be characterised using parabolic inductions. (See
Cor.~\ref{cor:supercuspidal}.)

\section{Computing the Satake transform}
\label{sec:compute-satake}

In this section we determine the inverse of the mod~$p$ Satake transform explicitly. We deduce this from the corresponding
result over the complex numbers, the Lusztig--Kato formula. The final result has a very simple shape: since $q = 0$ in $\k$,
only the constant terms of the intervening Kazhdan--Lusztig polynomials matter.

If $M$ is a standard Levi, we will denote by $\ge_M$ the usual partial order on $X_*(T)$ with respect to $M$, i.e., $\lambda \ge_M \mu$
means that $\lambda -\mu$ is a non-negative integral linear combination of the simple coroots of $M$. We also write $\ge$ for $\ge_G$.

\begin{prop}\label{prop:lusztig-kato}
  Suppose that the derived subgroup of $G$ is simply connected.
  Let $V$ be a weight.
  Let $M$ be the standard Levi subgroup such that $\Stab_W(\vuk) = W_M$. Then for all $\mu \in X_*(T)_-$,
  \begin{equation*}
    \tau_\mu = \sum_{\substack{\lambda \in X_*(T)_- \\ \lambda \ge_M \mu}} \SS_G(T_\lambda).
  \end{equation*}
\end{prop}

The assumption on $G$ should be unnecessary.
\mar{Need integral $z$-extension to remove it! Good example: $\det : \PGL_{q-1}(k) \to k\s$.}

\begin{proof}
  We can easily reduce to the case when $M = G$ by Cor.~\ref{cor:support}(ii). Just note that if $\lambda \in X_*(T)$ is antidominant for $M$
  and $\lambda \ge_M \mu$, then $\lambda \in X_*(T)_-$. (The point is that $\langle \alpha\dual,\beta\rangle \le 0$ for all $\alpha$, $\beta \in \Delta$.)

  Next we will reduce to the case when $V$ is the trivial weight, by showing that the coefficients of $\SS_G(T_\lambda)$ in
  the basis $(\tau_\mu)_\mu$ do not depend on $V$.  Since the derived subgroup of $G$ is simply connected, we can write $V =
  F(\nu)$ for some $q$-restricted weight $\nu$. We have $\nu \in X^0(T)$ since $\Stab_W(\nu) = W$.  The dual Weyl module
  $H^0_\O(\nu)$ (see \cite[\S II.8]{bib:Jan-reps}) of $G_{/\O}$ of highest weight $\nu$ is free of rank one (for example by
  the Weyl character formula \cite[Prop.~II.5.10]{bib:Jan-reps}). In particular there is a character $\wt\nu: G \to F\s$ that
  agrees with $\nu$ on $T$ and whose reduction modulo $\varpi$ is the one-dimensional representation $V$. (Note that
  $\wt\nu(K) \subset \O\s$, as $K = G(\O)$ is compact.
  Alternatively we could obtain $\wt\nu$ using the isogeny of $F$-tori $T \to G/G'$, where $G'$ is the derived subgroup of $G$.)
  Let $\lambda$, $\mu \in X_*(T)_-$ and put $t := \lambda(\varpi)$, $t' := \mu(\varpi)$. For $g = k_1 t k_2 \in KtK$ we see
  that $T_\lambda(g) = k_1 k_2 \in \End_\k(V) = \k$, i.e., $T_\lambda(g) = \overline{\wt\nu(g)\wt\nu(t)^{-1}} \in \k\s$.  Now
  we follow a standard argument (e.g., \cite[\S 3]{bib:Gross_Satake}). By the Iwasawa decomposition we can write $KtK =
  \coprod_{i=1}^r t_iu_i K$ with $t_i = \lambda_i(\varpi)$, $\lambda_i \in X_*(T)$, $u_i \in U$.  Then $t_iu_iK \cap t'U \ne
  \varnothing$ if and only if $\lambda_i = \mu$, in which case the intersection equals $t' u_i U(\O)$. Thus the coefficient
  of $\tau_\mu$ in $\SS_G (T_\lambda)$ equals
  \begin{equation*}
    \sum_{i:\lambda_i = \mu} T_\lambda(t'u_i) = \#\{ i : \lambda_i = \mu \} \cdot \overline{\wt\nu(t')\wt\nu(t)^{-1}} = \#\{ i : \lambda_i = \mu \},
  \end{equation*}
  where the last equality follows since $\mu \ge_{\R} \lambda$ if $KtK \cap t'U \ne \varnothing$ \cite[Lemma~3.6]{bib:satake}
  and thus $\wt\nu(t't^{-1}) = \varpi^{\langle \nu,\mu-\lambda\rangle} = 1$.
  It is therefore indeed independent of $\nu \in X^0(T)$.
  \mar{alternatively: $\wt\nu$ induces smooth char. $G \to F\s \to \O\s \to k\s$ extending $\nu$ on $G(k)$. then
    use twisting argument as in Lemma~\ref{lm:hecke-evals-twist}.}
  
  We now prove the proposition in case $V$ is trivial.
  Recall the Lusztig--Kato formula \cite[Thm.~7.8.1]{bib:HKP} (see also \cite[\S 4]{bib:Gross_Satake}). This is an identity in $\Z[q^{1/2},q^{-1/2}][X_*(T)]$ which,
  when the variable $q^{1/2}$ is specialised to a complex square root of $q = \# k$, gives the following in $\C[X_*(T)]$:
  \begin{equation*}
    \ch V_\mu = \sum_{\lambda \le \mu} q^{-\langle \mu,\rho\rangle} P_{w_\lambda,w_\mu}(q) 1_{K\lambda(\varpi)K}\dual.
  \end{equation*}
  Here $\lambda$, $\mu$ denote any \emph{dominant} coweights, $V_\mu$ is the irreducible complex representation of highest
  weight $\mu$ of the dual group, $\rho$ is the half-sum of all positive roots, $w_\lambda$ is the element $\lambda\cdot w_0$ in the
  extended affine Weyl group $\wt W = X_*(T) \rtimes W$, where $w_0$ is the longest Weyl element,
  and $1_{K\lambda(\varpi)K}\dual$ denotes the classical (i.e.,
  normalised) Satake transform of the characteristic function of $K\lambda(\varpi)K$. Besides, for any elements $w \le w'$ in
  $\wt W$, $P_{w,w'}(q) \in 1+q\Z[q]$ denotes the corresponding Kazhdan--Lusztig polynomial. We remark that $\lambda \le \mu$
  if and only if $w_\lambda \le w_\mu$ in $\wt W$.

  For any $\mu' \in X_*(T)$ take the coefficient of $\mu'$ in the above formula and rescale:
  \begin{equation*}
    q^{\langle\mu-w_0\mu',\rho\rangle} \dim V_\mu(\mu') = \sum_{\lambda \le \mu}  P_{w_\lambda,w_\mu}(q) \sum_{U/U(\O)} 1_{K\lambda(\varpi)K}(\mu'(\varpi)u).
  \end{equation*}
  Here $V_\mu(\mu')$ denotes the $\mu'$-weight space in $V_\mu$. Also note that $\delta_B^{1/2}(\mu'(\varpi)) = q^{\langle w_0\mu',\rho\rangle}$
  for the modulus character of $B$.

  Consider the left-hand side. We have $V_\mu(\mu') = 0$ unless $w_0\mu \le \mu' \le \mu$. If these inequalities hold, then
  $\langle \mu-w_0\mu',\rho\rangle \ge 0$ with equality if and only if $\mu = w_0\mu'$. Therefore we may reduce both sides modulo $p$
  and obtain
  \begin{equation*}
    \tau_{w_0\mu} = \sum_{\lambda \le \mu} \SS_G(T_{w_0\lambda}),
  \end{equation*}
  noting again that the Kazhdan--Lusztig polynomials have constant coefficient~1. Finally we interchange $\lambda$ with $w_0\lambda$ and
  $\mu$ with $w_0\mu$.
\end{proof}

\section{Maps between compact inductions}
\label{sec:maps-betw-cpt-ind}

\subsection{Generalities}
\label{sec:cpt-ind-generalities}

Suppose that we are given compact open subgroups $H_i$ of $G$ and finite-dimensional smooth $H_i$-representations
$V_i$ ($i = 1$, $2$). \mar{in fact it seems it's enough that the $V_i$ are finitely generated as $H_i$-reps.}We define
\begin{equation*}
  \HH_{H_1,H_2}(V_1,V_2) := \Hom_G (\ind_{H_1}^G V_1, \ind_{H_2}^G V_2).
\end{equation*}
If moreover $V_3$ is a finite-dimensional smooth representation of a compact open subgroup $H_3$, we have a natural bilinear map
given by composition
\begin{equation}\label{eq:13}
  \HH_{H_2,H_3}(V_2,V_3) \times \HH_{H_1,H_2}(V_1,V_2) \to \HH_{H_1,H_3}(V_1,V_3).
\end{equation}
In particular $\HH_{H_1,H_2}(V_1,V_2)$ is a Hecke bimodule, with $\HH_{H_1}(V_1)$ acting on the
right and $\HH_{H_2}(V_2)$ acting on the left.

By Frobenius reciprocity, $\HH_{H_1,H_2}(V_1,V_2)$ is isomorphic to $\Hom_{H_1}(V_1,\ind_{H_2}^G V_2)$ and thus,
by thinking of it inside the space of functions on $G \times V_1 \to V_2$,
to
\begin{align*}
  \{ \vp : G \to \Hom_{\k}(V_1,V_2) : {}&\text{$\supp \vp$ compact}, \\ &\vp(h_2gh_1) = h_2\circ \vp(g)\circ h_1 \ \forall h_i \in H_i, g \in G \}.
\end{align*}
In this language the composition (\ref{eq:13}) is given by convolution:
$(\vp' \ast \vp)(g) = \sum_{G/H_2} \vp'(gx) \vp(x^{-1})$.

If $H_1 = H_2 = K$ and the $V_i$ are weights, we just write $\HH_G(V_1,V_2)$ for $\HH_{K,K}(V_1,V_2)$.

\begin{prop}\label{prop:maps-betw-cpt-ind}
  Suppose that $V_1$, $V_2$ are two weights. Then $\HH_G(V_1,V_2)$ is non-zero if and only if $V_1^{U(k)} \cong
  V_2^{U(k)}$ as $T(k)$-representations.

  If this is satisfied, the Hecke algebras $\HH_G(V_1)$ and $\HH_G(V_2)$ can naturally
  be identified via $\SS_G$. Under this identification the actions of the two Hecke algebras on $\HH_G(V_1,V_2)$ agree.
  If moreover the centre of $G$ is connected and the derived subgroup of $G$ is simply connected
  then $\HH_G(V_1,V_2)$ is a free module of rank one under $\HH_G(V_1) \cong \HH_G(V_2)$.
\end{prop}

This follows immediately from the following proposition and its proof.

\begin{prop}\label{prop:maps-betw-cpt-ind-satake}
  Suppose that $V_1$, $V_2$ are two weights. The Satake transform
  \begin{align*}
    \SS_G : \HH_G(V_1,V_2) &\to \HH_T(V_1^{U(k)},V_2^{U(k)}) \\
    \vp &\mapsto \left(t \mapsto \sum_{u \in U/U(\O)} \vp(tu)\Big|_{V_1^{U(k)}}\right)
  \end{align*}
  is injective. The transform is compatible with compositions:
  \ifkuvio
  \begin{equation}\label{eq:19}
    \cellwidth=45pt \Diag
    \HH_{G}(V_2,V_3) & \times & \HH_{G}(V_1,V_2) &\rTo &\HH_{G}(V_1,V_3) \\
    \dTo_{\SS_G}  && \dTo_{\SS_G} && \dTo_{\SS_G} \dy{-3mm} \\
    \HH_T(V_2^{U(k)},V_3^{U(k)}) & \;\times\; & \HH_T(V_1^{U(k)},V_2^{U(k)}) & \rTo & \HH_T(V_1^{U(k)},V_3^{U(k)}) \\ 
    \endDiag 
  \end{equation}
  \else
  \begin{equation}\label{eq:19}\end{equation}
  \fi
  If $V_1^{U(k)} \cong V_2^{U(k)}$ as $T(k)$-representation, $G$ has connected centre, and the derived subgroup of $G$ is
  simply connected, the image of $\SS_G$ is a free \mar{maybe need to assume $G'$ simply conn. in this last statement?}
  module of rank one under the compatible actions of $\HH_T^-(V_1^{U(k)}) \cong \HH_T^-(V_2^{U(k)})$.
\end{prop}

\begin{proof}
  That the transform is well defined and compatible with compositions is formal using the Iwasawa decomposition
  and follows exactly the same steps as in the case when $V_1 = V_2$ (see~\cite{bib:satake}).
  
  Suppose $\lambda \in X_*(T)_-$. By Step~1 of the proof of \cite[Thm.~1.2]{bib:satake}, we see that the vector space of $\vp
  \in \HH_G(V_1,V_2)$ that are supported on $K \lambda(\varpi) K$ is one-dimensional if $V_1^{N_{-\lambda}(k)} \cong
  V_2^{N_{-\lambda}(k)}$ as $M_\lambda(k)$-representations and zero otherwise. Since $N_{-\lambda}(k) \subset U(k)$ it
  follows that $\HH_G(V_1,V_2) \ne 0$ implies $V_1^{U(k)} \cong V_2^{U(k)}$ as $T(k)$-representation.  To see the converse,
  suppose $V_1^{U(k)} \cong V_2^{U(k)}$ and choose any $\lambda \in X_*(T)_-$ such that $\langle \lambda,\alpha\rangle < 0$ for
  all $\alpha \in \Delta$.  Then $N_{-\lambda} = U$, so that there exists a non-zero Hecke operator supported on $K \lambda(\varpi) K$.

  We pick $\vp \ne 0$ in $\HH_G(V_1,V_2)$ and show that $\SS_G(\vp)$ is non-zero. Suppose first that $\vp$ is supported on $K
  \lambda(\varpi) K$ for some $\lambda \in X_*(T)_-$. For $\mu \in X_*(T)$ the same argument as in Step~3 of
  \cite[Thm.~1.2]{bib:satake} shows that $\SS_G(\vp)(\mu(\varpi)) \ne 0$ implies $\mu \ge_{\R} \lambda$ and that
  $\SS_G(\vp)(\lambda(\varpi)) \ne 0$. For general $\vp$, picking a minimal $\lambda \in
  X_*(T)_-$ for $\ge_\R$ such that $\vp$ is non-zero on $K \lambda(\varpi)K$, we see that
  $\SS_G(\vp)(\lambda(\varpi)) \ne 0$.  Thus $\SS_G$ is injective. (This gives an alternative proof that $\HH_G(V_1,V_2) = 0$
  unless $V_1^{U(k)} \cong V_2^{U(k)}$.)

  Let us now assume that $V_1^{U(k)} \cong V_2^{U(k)}$ and that $G$ has connected centre and simply connected derived subgroup.
  We first show that there is a $\lambda_0 \in X_*(T)_-$ such that there is a non-zero $\vp \in \HH_G(V_1,V_2)$ supported
  on $K\lambda(\varpi)K$ if and only if $\lambda-\lambda_0 \in X_*(T)_-$.
  We can write $V_i \cong F(\nu_i)$ with $\nu_i$ $q$-restricted. As $V_1^{U(k)} \cong V_2^{U(k)}$, we have $\nu_1-\nu_2 \in (q-1)X^*(T)$.
  There is a non-zero $\vp \in \HH_G(V_1,V_2)$ supported
  on $K\lambda(\varpi)K$ if and only if $V_1^{N_{-\lambda}(k)} \cong V_2^{N_{-\lambda}(k)}$ if and only if
  $\langle\nu_1-\nu_2,\alpha\rangle = 0$ for all roots $\alpha$ of $M_\lambda$, i.e., for all $\alpha \in \Phi$ such
  that $\langle \lambda,\alpha\rangle = 0$ (the last step uses \cite[Lemma 2.5]{bib:satake} and Prop.~1.3 in the appendix of \cite{bib:thesis}).
  Since the centre of $G$ is connected we may choose $\lambda_0 \in X_*(T)_-$ such that for all simple roots $\alpha$,
  $\langle \lambda_0,\alpha\rangle = 0$ if $\langle \nu_1-\nu_2,\alpha\rangle = 0$ and $\langle \lambda_0,\alpha\rangle = -1$ otherwise.
  This clearly satisfies the desired property.

  We can canonically identify $\HH_T(V_1^{U(k)})$, $\HH_T(V_2^{U(k)})$ and we denote them simply by $\HH_T$. By
  Thm.~\ref{thm:satake} and~\eqref{eq:19} we see that $\HH_T^-$ preserves the image of $\SS_G$.  Choose $\vp_0 \in
  \HH_G(V_1,V_2)$ that is non-zero and supported on $K \lambda_0(\varpi)K$.  We show below that $\HH_T^-\SS_G(\vp_{0}) =
  \im(\SS_G)$. Since $\HH_T \cong \k[X_*(T)]$ is an integral domain, this completes the proof.

  \emph{Claim:} For all $\lambda \in \lambda_0 + X_*(T)_-$ there is a $\psi_\lambda \in \HH_T^-\SS_G(\vp_{0})$ such that for
  all $\mu \in \lambda_0 + X_*(T)_-$, $\psi_\lambda(\mu(\varpi)) \ne 0$ if and only if $\mu = \lambda$.
  
  We prove the claim by induction with respect to $\ge_\R$, noting that $\{\mu \in X_*(T)_-: \mu \ge_\R \lambda\}$ is finite.
  Consider $\tau_{\lambda-\lambda_0} \SS_G(\vp_0) \in \HH_T^-\SS_G(\vp_{0})$. By the above this is zero on $\mu(\varpi)$
  unless $\mu \ge_\R \lambda$ and non-zero on $\lambda(\varpi)$. If there is no $\mu >_\R \lambda$ in $\lambda_0 + X_*(T)_-$ we are done.
  Otherwise by induction we can subtract off multiples of $\psi_\mu$ for such $\mu$ to find a $\psi_\lambda$.

  Finally given $\psi \in \im(\SS_G)$, by the claim we can subtract off multiples of the $\psi_\lambda \in \HH_T^-\SS_G(\vp_{0})$ to assume without loss
  of generality that $\psi(\lambda(\varpi)) = 0$ for all $\lambda \in \lambda_0 + X_*(T)_-$. Then the above argument for the injectivity
  of $\SS_G$ shows that $\psi = 0$.
\end{proof}

\begin{coroll}
  Suppose $V$ is a weight. Then $\kind V$ is a torsion-free $\HH_G(V)$-module.
\end{coroll}

\begin{proof}
  If $\kind V$ has torsion, there is a non-zero $\vp : \kind V \to \kind V$ in $\hg$ that has non-zero kernel. As a non-zero smooth
  $G$-representation, the kernel has to contain a weight $V'$. By Frobenius reciprocity we get a non-zero map $\kind V'
  \to \kind V$ whose composite with $\vp$ is zero. More generally suppose that a composite of non-zero maps
  $\vp_1 : \kind V_1 \to \kind V_2$ and $\vp_2 : \kind V_2 \to \kind V_3$ is zero. By Prop.~\ref{prop:maps-betw-cpt-ind}, the $T(k)$-representations $V_1^{U(k)}$,
  $V_2^{U(k)}$, and $V_3^{U(k)}$ are isomorphic, and we identify them (non-canonically). Then for all $i$, $j$ the bimodules
  $\HH_T(V_i^{U(k)},V_j^{U(k)})$ are all naturally identified with the integral domain $\HH_T(V_1^{U(k)})$ such that
  moreover for any triple $(i,j,k)$ the bimodule multiplication corresponds to the ring multiplication.
  Since $\SS_G(\vp_2) \ast \SS_G(\vp_1) = 0$ and $\SS_G$ is injective, one of $\vp_1$, $\vp_2$ has to be zero.
\end{proof}

\subsection{The minuscule case}
\label{sec:minuscule-case}

Let $V$, $V'$ denote distinct weights such that $V^{U(k)} \cong (V')^{U(k)}$. As we just saw, in this
case we can identify $\hg$ and $\HH_G(V')$, and we will denote them simply by $\HH_G$.
The goal of this subsection is to find an explicit criterion when
\begin{equation}
\kind V \otimes_{\HH_G,\chi} \k \cong \kind V' \otimes_{\HH_G,\chi} \k,\label{eq:20}
\end{equation}
in case $V$ and $V'$ differ only minimally (in the sense that $V^{N(k)} \cong (V')^{N(k)}$ as $M(k)$-representations for some
maximal parabolic $P = MN$) and provided there exists a corresponding minuscule fundamental coweight.
We remark that when $G = \GL_n$, any simple root admits a minuscule fundamental coweight.

An isomorphism as in \eqref{eq:20} is extremely useful because it allows us to ``change the weight'' in a smooth $G$-representation $\pi$.
The point is that $\Hom_G(\kind V \otimes_{\HH_G,\chi} \k,\pi) \ne 0$ is equivalent to saying that $V$ occurs in $\pi$ with
Hecke eigenvalues $\chi$. So if~\eqref{eq:20} holds and $V$ occurs in $\pi$ with Hecke eigenvalues $\chi$, then $V'$ occurs
in $\pi$ with Hecke eigenvalues $\chi$.
The results in this subsection will play a key role in the proofs of Theorems~\ref{thm:irred-parab-ind-gln}
and~\ref{thm:classification}.

Suppose that $G$ has simply connected derived subgroup. \mar{can remove this by $z$-ext.?} This implies that simple coroots possess fundamental weights. 
Fix a simple root $\alpha$. We denote by $\omega_\alpha$ a
fundamental weight for $\alpha\dual$, i.e., $\langle \beta\dual,\omega_\alpha\rangle = \delta_{\alpha\beta}$ for all simple roots
$\beta$. We make the following assumption.
\begin{center}
  \emph{There exists a minuscule fundamental coweight $-\lambda$ associated to $\alpha$.}
\end{center}
Recall that this means that $\langle -\lambda,\beta\rangle = \delta_{\alpha\beta}$ for $\beta \in \Delta$ and
$\langle -\lambda,\gamma\rangle \in \{0, 1\}$ for all $\gamma \in \Phi^+$. In particular, $\lambda \in X_*(T)_-$.

Suppose that $\nu$ is a $q$-restricted weight satisfying $\langle \nu,\alpha\dual\rangle = 0$. We let
$\nu' := \nu + (q-1)\omega_\alpha$ (again $q$-restricted) and define weights $V := F(\nu)$, $V' := F(\nu')$.
Note that $V^{U(k)} \cong (V')^{U(k)}$ as $T(k)$-representations. There exist Hecke operators
$\vp^+_\lambda \in \HH_G(V,V')$ and $\vp^-_\lambda \in \HH_G(V',V)$ whose support is $K\lambda(\varpi)K$. (In fact
$\lambda$ is precisely a possible $\lambda_0$ in the proof of Prop.~\ref{prop:maps-betw-cpt-ind-satake}.)

\begin{prop}\label{prop:minu-case}
  With the above notation, $\SS_G(\vp^-_\lambda \ast \vp^+_\lambda) \in \HH_T^-$ has support
  $\lambda(\varpi)^2 T(\O) \cup \lambda(\varpi)^2 \alpha\dual(\varpi) T(\O)$. The values at
  $\lambda(\varpi)^2$ and $\lambda(\varpi)^2 \alpha\dual(\varpi)$ add to zero.
\end{prop}

\begin{proof}
  We first compute $\vp^-_\lambda \ast \vp^+_\lambda \in \HH_G(V)$. Let $t = \lambda(\varpi)$.
  
  \begin{sublm}\label{sublm:minuscule}
    As $-\lambda$ is minuscule, we have $K \cap {}^t K = \red^{-1}(P_{-\lambda}(k))$.
  \end{sublm}

  Note that this is elementary when $G = \GL_n$.

  \begin{proof}[Proof of Sublemma~\ref{sublm:minuscule}]
    By Prop.~\ref{prop:buildings-lemma} we know that the left-hand side is contained in the right-hand side.
    We can rewrite the right-hand side by Lemma~\ref{lm:iwahori-decomp} as
    \begin{equation*}
      N_{-\lambda}(\O) M_\lambda(\O) \ker(N_\lambda(\O) \to N_\lambda(k)).
    \end{equation*}
    Since $t \in Z(M_\lambda)$, we see that $M_\lambda(\O) \subset K \cap {}^t K$.
    Lemma~\ref{lm:iwahori-decomp} shows that $t^{-1}$ contracts $N_{-\lambda}(\O)$ so that
    $N_{-\lambda}(\O) \subset K \cap {}^t K$.

    Fix an order on the set of roots $\alpha$ such that $\langle \lambda,\alpha\rangle > 0$. Multiplication induces an
    isomorphism $\prod_{\langle \lambda,\alpha\rangle > 0} U_\alpha \congto N_{\lambda}$ \cite[II.1.7(1)]{bib:Jan-reps}.
    Moreover we choose root homomorphisms $x_\alpha : \G_a \congto U_\alpha$ \cite[II.1.2]{bib:Jan-reps}.  An element of
    $N_\lambda(\O)$ can then be expressed uniquely as $\prod_{\langle \lambda,\alpha\rangle > 0} x_\alpha(u_\alpha)$ with
    $u_\alpha \in \O$. If its reduction is trivial in $N_\lambda(k)$ then $u_\alpha \in \varpi \O$ for all $\alpha$. Then
    $t^{-1} \prod x_\alpha(u_\alpha) t = \prod x_\alpha(\alpha(t^{-1}) u_\alpha)$, and this is contained in $K$ as $\langle
    \lambda,\alpha\rangle = 1$ whenever $\langle \lambda,\alpha\rangle > 0$ (as $-\lambda$ minuscule).  It follows that $\ker(N_\lambda(\O) \to
    N_\lambda(k)) \subset K \cap {}^t K$.
  \end{proof}

  For each $w \in W$ choose a representative $\dot w \in N(T)(\O)$.
  It follows from the sublemma that
  \begin{equation*}
    K/(K \cap {}^t K) \cong G(k)/P_{-\lambda}(k) \cong \coprod_{W/W_{\lambda}} U(k) \dot w P_{-\lambda}(k)/P_{-\lambda}(k),
  \end{equation*}
  where we used the rational Bruhat decomposition in the last step and where $W_{\lambda}$ denotes $W_{M_{\lambda}} = \Stab_W(\lambda)$.
  Therefore $KtK = \bigcup_W U(\O) \dot w t K$.

  We claim that the support of $\vp^-_\lambda \ast \vp^+_\lambda$ is $Kt^2K$. Let $\lambda' \in X_*(T)_-$ and let $t' =
  \lambda'(\varpi)$.  When we compute $(\vp^-_\lambda \ast \vp^+_\lambda)(t')$, by the above all terms are of the form
  $\vp^-_\lambda(u \dot w t) \vp^+_\lambda (t^{-1}\dot w^{-1} u^{-1} t')$, where $u \in U(\O)$ and $w \in W$.  As $t' \in
  T^-$, $(t')^{-1} u^{-1} t' \in U(\O)$, so if the term is non-zero then so is $\vp^-_\lambda({}^w t) \vp^+_\lambda ({}^w
  t^{-1} t')$ (where we also pulled a $\dot w$ from the left to the right).  By the refined Cartan decomposition we need that
  $-w\lambda + \lambda' = w'\lambda$ for some $w' \in W$. Moreover $\vp^-_\lambda({}^w t) \vp^+_\lambda ({}^{w'} t) \ne 0$,
  which implies that $\vp^-_\lambda(t)\dot w^{-1} \dot w' \vp^+_\lambda (t) \ne 0$.  This in turn implies that
  $p_{N_\lambda}(\dot w^{-1} \dot w' (V')^{N_{-\lambda}(k)}) \ne 0$. The proof of Lemma~\ref{lm:support} then shows that
  $w^{-1}w' \in W_\lambda$. (Note that $\Stab_W(\nu') \subset W_\lambda$ since the stabiliser is generated by simple
  reflections and since $\langle \nu',\alpha\dual\rangle = q-1 > 0$.) Thus $w'\lambda = w\lambda$.
  As $-w\lambda + \lambda' = w'\lambda$, we obtain $\lambda' = 2w\lambda$. Taking into account that $\lambda'$
  and $\lambda$ are antidominant coweights, we see that $w\lambda = \lambda$ and that $\lambda' = 2\lambda$.
  The latter equation shows that the support of $\vp^-_\lambda \ast \vp^+_\lambda$ is contained in $Kt^2K$.
  But if we take $\lambda' = 2\lambda$, the former equation shows that $u \dot w t = t (t^{-1}ut)\dot w \in tK$,
  so that only the trivial term $\vp^-_\lambda(t) \vp^+_\lambda(t)$ contributes to $(\vp^-_\lambda \ast \vp^+_\lambda)(t^2)$,
  but that term is clearly non-zero.

  To complete the proof, we will show that $\SS_G(T_{2\lambda}) = \tau_{2\lambda} - \tau_{2\lambda+\alpha\dual}$. 
  First note that $2\lambda + \alpha\dual \in X_*(T)_-$ since for $\beta \in \Delta$ we have
  \begin{equation*}
    \langle 2\lambda + \alpha\dual, \beta\rangle = 
    \begin{cases}
      0 & \text{if $\beta = \alpha$,} \\
      \langle \alpha\dual, \beta\rangle \le 0 & \text{if $\beta \ne \alpha$.}
    \end{cases}
  \end{equation*}
  Let $M$ be the standard Levi with $W_M = \Stab_W(\nu)$. Note that $\alpha \in \Delta_M$. 
  By Prop.~\ref{prop:lusztig-kato}, it suffices to show the following:

  \emph{Claim:} Suppose $\mu \in X_*(T)_-$. Then $\mu \ge_M 2\lambda$ if and only if $\mu = 2\lambda$ or $\mu \ge_M 2\lambda+\alpha\dual$.
  
  It suffices to prove the claim with $\ge_M$ replaced by $\ge$. Suppose $\mu \ge 2\lambda$, so $\mu-2\lambda = \sum_i
  \beta_i\dual$, where the $\beta_i$ are simple roots. Let $\alpha_0 \in X^*(T)$ be the sum of the longest roots of all irreducible
  components of the root system (it need not itself be a root). Then $\alpha_0 = \sum_{\beta \in \Delta} n_\beta \beta$ with
  $n_\beta \ge 1$. Moreover $\alpha_0$ is dominant and $\langle \lambda,\alpha_0\rangle = -n_\alpha = -1$ (as $-\lambda$ is minuscule).
  As $\mu$ is antidominant, we find
  \begin{equation*}
    \langle \mu,\alpha\rangle \ge \sum_{\beta\in\Delta} n_\beta \langle \mu,\beta\rangle = \langle \mu,\alpha_0\rangle \ge
    \langle 2\lambda,\alpha_0\rangle = -2.
  \end{equation*}
  If $\langle \mu,\alpha\rangle > -2$ then $\sum_i \langle \beta_i\dual,\alpha\rangle > 0$, so $\alpha$ has to occur among the simple roots
  $\beta_i$ and thus $\mu \ge 2\lambda + \alpha\dual$. If $\langle \mu,\alpha\rangle = -2$, then $\langle \mu,\beta\rangle = 0$ for
  all $\beta \in \Delta-\{\alpha\}$. Then $\mu-2\lambda = 0$ since it is orthogonal to all roots and contained in $\Z \Phi\dual$.
\end{proof}

\begin{qu}
  Suppose that the centre of $G$ is connected so that every simple root $\alpha$ possesses a fundamental coweight $\lambda$.
  Is the result of Prop.~\ref{prop:minu-case} true when $\lambda$ is not assumed to be minuscule? It seems that this is the case 
  when $G = \GSp_4$.
\end{qu}

\begin{coroll}\label{cor:minu-case-criterion}
  We keep the above notation.
  Suppose $\chi \in \Hom\kalg(\HH_G,\k)$ is parameterised by $(M,\chi_M)$. Assume that $\alpha \not\in \Delta_M$. Then
  \begin{equation}\label{eq:14}
    \kind V \otimes_{\HH_G,\chi} \k \cong \kind V' \otimes_{\HH_G,\chi} \k,
  \end{equation}
  provided either $\alpha\dual(\varpi) \not\in Z_M$ or $\chi_M(\alpha\dual(\varpi)) \ne 1$.
\end{coroll}

We can reformulate the conditions in a way that does not use~$\varpi$. The first is equivalent to
$\langle \alpha\dual,\beta\rangle \ne 0$ for some $\beta \in \Delta_M$, the second is equivalent to
$\chi_M \circ \alpha\dual : F\s \to k\s$ being non-trivial.

\begin{proof}
  By the above, $\vp^+_\lambda$, $\vp^-_\lambda$ induce $G$-linear maps between the two representations in~(\ref{eq:14}).  If
  we can show that $\chi(\vp^-_\lambda \ast \vp^+_\lambda) = \chi(\vp^+_\lambda \ast \vp^-_\lambda)$ is non-zero, we will be
  done.  Let $\chi'$ be the Satake transform of $\chi$. By Prop.~\ref{prop:minu-case} it suffices to show that
  $\chi'(\tau_{2\lambda}) \ne \chi'(\tau_{2\lambda+\alpha\dual})$. Recall that $\supp \chi' = Z_M^- T(\O)$.

  Since $\langle \lambda,\beta\rangle = 0$ for all simple roots $\beta \in \Delta_M$ (as $\alpha \not\in\Delta_M$),
  we have $\lambda(\varpi)^2 \in Z_M^-$; moreover $\lambda(\varpi)^2\alpha\dual(\varpi) \in Z_M^-$ if and only if
  $\alpha\dual(\varpi) \in Z_M$. We now use Cor.~\ref{cor:param-hecke-evals}.
  If $\alpha\dual(\varpi) \not\in Z_M$ then $\chi'(\tau_{2\lambda}) \ne 0 = \chi'(\tau_{2\lambda+\alpha\dual})$.
  Otherwise both are non-zero and $\chi'(\tau_{2\lambda})\chi'(\tau_{2\lambda+\alpha\dual})^{-1} = \chi_M(\alpha\dual(\varpi)) \ne 1$.

  In the latter case let us verify that $\chi_M \circ \alpha\dual$ is unramified, i.e., trivial on $k\s$. We have
  that $\chi_M|_{Z_M(k)}$ is the central character of $F(\nu)_{\o N(k)}$, hence it equals $\nu|_{Z_M(k)}$.
  Thus if $\langle \alpha\dual,\beta\rangle = 0$ for all $\beta \in \Delta_M$ then $\alpha\dual(x) \in Z_M(k)$ for $x \in k\s$,
  so $\chi_M(\alpha\dual(x)) = \nu(\alpha\dual(x)) = 1$. (Recall that $\langle \nu, \alpha\dual\rangle = 0$.)
\end{proof}

\begin{rk}\mar{this shows that the trivial and Steinberg weights occur
  simultaneously in supersingulars of $\GL_2(F)$}
  Equation~(\ref{eq:14}) also holds if $\alpha \in \Delta_M$ and $\langle \alpha\dual,\beta\rangle = 0$
  for all $\beta \in \Delta_M-\{\alpha\}$.
\end{rk}

We showed in~\cite{bib:satake} (Step~4 of the proof of Thm.~1.2) that the image of $\SS_G$ is supported on the ``almost antidominant'' part of the torus, namely the set
of those $t \in T$ such that $\ord_F(\alpha(t)) \le 1$ for all $\alpha \in \Delta$.
(The proof still goes through when $V_1 \not\cong V_2$.)
We now show that the support of the image of $\SS_G$ may extend outside $T^-$ if $V_1 \not\cong V_2$.

\begin{prop}
  With the above notation, $\SS_G(\vp^+_\lambda) \in \HH_T$ has support
  $\lambda(\varpi) T(\O) \cup \lambda(\varpi) \alpha\dual(\varpi) T(\O)$. 
\end{prop}

Note that $\langle \lambda + \alpha\dual,\alpha\rangle = 1$.
\mar{I think more generally, for $V$, $V'$, the Satake transform of map between them (in some given support) only depends
on $\Stab_W(\vuk)$ and $\Stab_W((V')^{U(k)})$. Point: reduce to $V$ being one-dimensional using Cor.~\ref{cor:support}
(as in \ref{prop:lusztig-kato}), then use twisting as in \ref{lm:hecke-evals-twist}/\ref{prop:lusztig-kato}. Once
$V$ is trivial, the other weight has to be one of the $V_P$, so just down to shape of Levi\dots}

\begin{proof}
  By Prop.~\ref{prop:minu-case} it suffices to show that $\SS_G(\vp^-_\lambda)$ has support $\lambda(\varpi)T(\O)$.
  Let $M = M_\lambda$. We first claim that $\smg(\vp^-_\lambda) \in \HH_M((V')^{N(k)}, V^{N(k)})$ has support $M(\O) \lambda(\varpi) M(\O)$.
  Since $\langle \nu',\alpha\dual\rangle > 0$ and $-\lambda$ is a fundamental coweight for $\alpha$, it follows
  that $V'$ is $M_\lambda = M$-regular. Then Cor.~\ref{cor:support}(ii) establishes the claim. (Note that the proof still goes through
  and that it only depends on the ``source'' weight $V'$.)
  
  Since $\nu' - \nu = (q-1)\omega_\alpha$ and $\omega_\alpha$ is orthogonal to all simple coroots of~$M$, we see that
  $(V')^{N(k)} \cong V^{N(k)}$. By fixing an isomorphism between them we are reduced to consider the image of $T_\lambda^M$
  under $\SS_M : \HH_M(V^{N(k)}) \to \HH_T(\vuk)$. But this equals $\tau_\lambda$, as $\lambda(\varpi) \in Z_M$. (Just as in
  the proof of Lemma~\ref{lm:hecke-evals-supersing}.)
\end{proof}

\begin{ex}
  Consider $G = \GL_2$ with diagonal torus $T$ and upper-triangular Borel $B$.  A pair of integers $(a,b)$ denotes the
  element of $X^*(T)$ that sends $\big(\begin{smallmatrix} t_1 \\ & t_2 \end{smallmatrix}\big)$ to $t_1^a t_2^b$.  Let $V =
  F(0,0)$, the trivial weight, and $V' = F(q-1,0) \cong \Sym^{q-1} \k^2$, the Steinberg weight.  Then we are in the situation
  considered above, with $\alpha = (1,-1)$. In particular there are non-zero $\HH_G[G]$-linear maps $\kind V \to \kind V'$
  and $\kind V' \to \kind V$ whose support is $K\big(\begin{smallmatrix} 1 \\ &\varpi \end{smallmatrix}\big)K$. A variant of
  the latter map was used in~\cite[Lemma~1.5.5]{bib:Kisin}.
\end{ex}

\section{Generalised Steinberg representations}
\label{sec:gen-steinberg}

In this section we will determine the Jordan--H\"older factors of $\Ind_{\o B}^G 1$ (completing the work of Gro\ss e-Kl\"onne) as well as
their weights and Hecke eigenvalues. As usual, $P = MN$ and $Q$ denote standard parabolic subgroups of $G$ throughout.

Recall that the \emph{generalised Steinberg representations} are defined as follows:
\begin{equation}\label{eq:26}
  \Sp P = \frac{\Ind_{\o P}^G 1}{\sum_{Q \supsetneq P}\Ind_{\o Q}^G 1}.
\end{equation}

\begin{thm}\label{thm:steinb-irred}
  For any standard parabolic subgroup $P$, the generalised Steinberg representation $\Sp P$ is irreducible and admissible.
\end{thm}

For the Steinberg representation $\Sp B$ this was proved by Vign\'eras \cite[\S 4]{bib:vigneras}. The general case was proved
by Gro\ss e-Kl\"onne \cite[Cor.~4.3]{bib:grosse-kloenne}, under the assumption that $G$ is of type A, B, C, or D. We recall Cor.~4.4
of~\cite{bib:grosse-kloenne}, which now holds without any restriction on the root system.

\begin{coroll}\label{cor:steinb-irred-constit}
  The generalised Steinberg representations $\Sp P$ are pairwise non-isomorphic. They form the irreducible constituents of $\Ind_{\o B}^G 1$,
  each occurring with multiplicity one. In particular, $\Ind_{\o B}^G 1$ is of finite length $2^{\# \Delta}$. 
\end{coroll}

We deduce Thm.~\ref{thm:steinb-irred} from the work of Gro\ss e-Kl\"onne with the help of Thm.~\ref{thm:cptind-parabind}.
More precisely we show that $\Sp P$ contains a unique weight, that that weight lifts to $\Ind_{\o P}^G 1$, and finally that
the lifted weight generates $\Ind_{\o P}^G 1$.

\begin{prop}\label{prop:steinb-wts}
  For any standard parabolic subgroup $P$, the representation $\Sp P$ contains a unique weight $V_P$. It occurs with multiplicity
  one.  The unique set of Hecke eigenvalues in weight $V_P$ is parameterised by the pair $(T,1)$.
\end{prop}

We will see in the proof that $V_P$ is the unique $M$-regular weight such that $(V_P)^{N(k)} \cong (V_P)_{\o N(k)}$ is trivial.
It follows that $\Stab_W (V_P^{U(k)}) = W_M$, so that the weights $V_P$ are pairwise distinct.

\begin{proof}[Proof of Theorem~\ref{thm:steinb-irred} and Proposition~\ref{prop:steinb-wts}]
  Let $\oSp P$ denote the generalised Steinberg representation of $G(k)$, i.e.,
  \begin{equation*}
    \oSp P = \frac{\Ind_{\o P(k)}^{G(k)} 1}{\sum_{Q \supsetneq P}\Ind_{\o Q(k)}^{G(k)} 1}.
  \end{equation*}
  Then by~\cite[\S 3]{bib:grosse-kloenne}, \mar{improve ref. once Elmar has revised}
  there is a natural $K$-linear embedding $\oSp P \INTO \Sp P$ such that
  \begin{equation}\label{eq:15}
    (\oSp P)^{B(k)} = (\Sp P)^I = (\Sp P)^{I(1)}.
  \end{equation}
  Since $I(1)$ is a pro-$p$ group, we see in particular that $\Sp P$ is admissible.
  
  \emph{Step 1.} We show that $\soc_K(\Sp P)$ is irreducible. Let $V \subset \Sp P$ be any weight.  By~\eqref{eq:15} we have
  $V^{U(k)} \subset (\oSp P)^{B(k)}$, so $\vuk = V^{B(k)}$ is a one-dimensional subspace of $(\oSp P)^{B(k)}$, which is
  stable under the action of the Hecke algebra $\k[B(k)\backslash G(k)/B(k)]$. In \cite[Prop.~3.4]{bib:grosse-kloenne} it is
  shown that any non-zero $\k[B(k)\backslash G(k)/B(k)]$-submodule of $(\oSp P)^{B(k)}$ contains the element $g_P$, which is by definition
  the image of $1_{\o P(k)B(k)} \in \Ind_{\o P(k)}^{G(k)} 1$ in $\oSp P$.  This shows that $V$ is generated by $g_P$ as
  $G(k)$-module.
  
  \emph{Step 2.} By Lemma~\ref{lm:lift-to-m-reg-wt}, there is a unique $M$-regular weight $V_P$ such that $(V_P)_{\o N(k)}$
  is trivial.  Let $\chi : \HH_M(1) \to \k$ denote the Hecke eigenvalues of the trivial weight in the trivial
  $M$-representation. It is parameterised by the pair $(T,1)$, for example by Lemmas~\ref{lm:hecke-evals-parab-ind}
  and~\ref{lm:hecke-evals-supersing} since the trivial $M$-representation is contained in the trivial principal series for
  $M$.  From Thm.~\ref{thm:cptind-parabind} we get a surjective $G$-linear map $\kind V_P \otimes_{\HH_G(V_P),\chi} \k \onto
  \opInd 1$.  Thus the image of $V_P$ generates $\opInd 1$ as $G$-representation; a fortiori, its image in $\Sp P$ is
  non-zero and generates. By Step~1, $V_P$ is the unique weight of $\Sp P$.
  
  \emph{Step 3.} Suppose that $\pi \subset \Sp P$ is a non-zero subrepresentation. Then $\pi$ contains a weight. By the previous
  two steps, we know that this is $V_P$ and that $V_P$ generates $\Sp P$. Thus $\pi = \Sp P$.

  The final statement of Prop.~\ref{prop:steinb-wts} follows from the proof of Lemma~\ref{lm:hecke-evals-parab-ind}. That is,
  $\chi \circ \psmg$ is parameterised by $(T,1)$ because $\chi$ is.
\end{proof}

The following proposition will be useful later.

\begin{prop}\label{prop:opind-spq}
  Suppose that $P = MN$ is a standard parabolic and that $Q$ is a standard parabolic of $M$. Then the constituents of $\opInd
  \Sp Q$ are given by $\Sp {P'}$, where $P'$ runs through all standard parabolic subgroups of $G$ such that $P' \cap M = Q$.
  They all occur with multiplicity one.
\end{prop}

\begin{proof}
  The last statement follows from Cor.~\ref{cor:steinb-irred-constit}: $\opInd \Sp Q$ is a subquotient of
  $\obInd 1$ since $\Sp Q$ is a subquotient of $\Ind_{B \cap M}^M 1$.
  
  For each standard parabolic subgroup $P'$ of~$G$, $P'\cap M$ is a standard parabolic subgroup of $M$.  
  In the parameterisation of standard parabolic subgroups by subsets of $\Delta$ the map $P' \mapsto P' \cap M$ becomes
  $\Delta' \mapsto \Delta' \cap \Delta_M$. Thus there is a smallest standard parabolic subgroup $Q'_0$ of $G$ such
  that $Q'_0 \cap M = Q$ and we have $P' \supset Q'_0$ if and only if $P' \cap M \supset Q$.
  In fact, $Q'_0 = QN$.

  By exactness and transitivity of parabolic induction,
  \begin{equation*}
    \opInd \Sp Q \cong \frac{\opInd \Ind_{\o Q}^M 1}{\sum_{R \supsetneq Q} \opInd \Ind_{\o R}^M 1} \cong 
    \frac{\Ind_{\o {Q'_0}}^G 1}{\sum_{R \supsetneq Q} \Ind_{\o {R'_0}}^G 1}.
  \end{equation*}
  Using~\eqref{eq:26} again it follows by induction that the constituents of $\opInd 1$ are given by $\Sp {P'}$ for $P' \supset P$.
  The claim now follows, recalling that we have already established that all constituents occur with multiplicity one.
\end{proof}

\section{Irreducibility of parabolic inductions}
\label{sec:irreducibility}

The goal of this section is to construct many irreducible admissible $G$-representations by parabolic induction.
The most precise results are obtained for $\GL_n$.

We remark that parabolic induction preserves admissibility. Suppose that $P = MN$ is a standard parabolic and that $\sigma$
is an admissible $M$-representation. It suffices to show that $(\opInd \sigma)^{K(1)}$ is finite-dimensional.  This follows
from the admissibility of $\sigma$ since
\begin{equation*}
  (\opInd \sigma)^{K(1)} = (\Ind_{\o P(\O)}^K \sigma)^{K(1)} = \Ind_{\o P(k)}^{G(k)} (\sigma^{M(1)}),
\end{equation*}
where $M(1)$ is the kernel of $M(\O) \to M(k)$.

\subsection{The case of $\GL_n$}
\label{sec:parab-ind-gln}

When $G = \GL_n$ we will always let $T$ be the diagonal torus and $B$ the upper-triangular Borel subgroup.

\begin{thm}\label{thm:irred-parab-ind-gln}
  Suppose that $G = \GL_n$. Let $P$ be the standard parabolic with Levi $\prod_i \GL_{n_i}$, where $\sum_i n_i = n$.
  Suppose that $\sigma_i$ is an irreducible admissible supersingular representation of $\GL_{n_i}(F)$ for $i = 1$, \ldots, $r$.
  Then $\Ind_{\o P}^G (\sigma_1 \otimes \cdots \otimes \sigma_r)$ is irreducible if and only if
  there is no $i$ such that $n_i = n_{i+1} = 1$ and $\sigma_i \cong \sigma_{i+1}$.
\end{thm}

When $n_1 = \cdots = n_r = 1$, this was proved by Ollivier \cite{bib:ollivier}. Her method relied on a detailed knowledge of
the $I(1)$-Hecke algebra and the study of intertwinings of certain induced $I(1)$-Hecke modules.  Henniart (unpublished)
found a different proof in that case, which utilises the structure of a principal series as $B$-representation and which works, under some
conditions, for quasi-split groups.  By contrast, our strategy is to show that every non-zero subrepresentation of an induced
representation satisfying the criterion in Thm.~\ref{thm:irred-parab-ind-gln} has to contain a sufficiently regular
weight (by ``changing the weight'' as in \S\ref{sec:minuscule-case}), and that such a weight has to generate the whole representation
(by using the surjectivity of the map in Thm.~\ref{thm:cptind-parabind}).

\begin{lm}\label{lm:product-decomp}
  Suppose that $G = \prod_{i=1}^r G_i$ is a product of split reductive groups.
  If $\sigma_i$ is an irreducible admissible $G_i$-representation for $i = 1$, \ldots, $r$, then
  $\sigma := \bigotimes_i \sigma_i$ is an irreducible admissible $G$-representation. Conversely, every
  irreducible admissible $G$-representation is of this form.
  
  The weights of $\sigma$ are of the form $V := \bigotimes_i V_i$
  where $V_i$ is an arbitrary weight of $\sigma_i$. The Hecke eigenvalues of
  $V$ in $\sigma$ are parameterised by the pairs $(\prod_i M_i,\prod \chi_i)$, where
  each $(M_i,\chi_i)$ ranges over the pairs parameterising Hecke eigenvalues of $V_i$ in $\sigma_i$.
\end{lm}

\begin{proof}
  The first part follows just as in the case of finite groups (using a proof that does not use
  character theory). Similarly, the weights
  for $G$ are of the form $\bigotimes V_i$ where each $V_i$ is a weight for $G_i$.
  Next note that $\bigotimes \Hom_{G_i(\O)}(V_i,\sigma_i) \cong \Hom_{G(\O)}(V,\sigma) $. Finally
  we have a natural isomorphism $\bigotimes \HH_{G_i}(V_i) \cong \HH_G(V)$ which is induced by the map
  $\otimes \vp_i \mapsto \vp$ with $\vp(g_1,\dots,g_r) = \bigotimes \vp_i(g_i)$.
  (To see that it is a bijection, note that it sends $\bigotimes T_{\lambda_i}$ to $T_{\prod \lambda_i}$.)
  One also verifies that this isomorphism is compatible with the Satake transform (for $G$ and the $G_i$).
  The lemma then follows easily from Prop.~\ref{prop:param-hecke-evals}.
\end{proof}

\begin{proof}[Proof of Theorem~\ref{thm:irred-parab-ind-gln}]
  Let $M_i = \GL_{n_i}$. We have $P = MN$ with $M = \prod_i \GL_{n_i}$. We let $\sigma = \sigma_1 \otimes \cdots \otimes
  \sigma_r$; it is an irreducible admissible representation of $M$.  Note first that the criterion is clearly necessary since
  otherwise $\sigma$ extends to a representation of a larger parabolic subgroup.

  To show that the criterion is sufficient, assume that $\pi \subset \Ind_{\o P}^G \sigma$ is a non-zero subrepresentation.  
  Let $V$ be a weight of $\pi$. We claim that the action of $\hg$ on $\Hom_K(V,\pi)$ factors through $\psmg : \hg \to \hm$,
  which is a localisation map by Prop.~\ref{prop:partial-satake}. Consider the natural maps of finite-dimensional
  $\hg$-modules,
  \begin{equation}\label{eq:22}
    \Hom_K(V,\pi) \INTO \Hom_K(V,\opInd \sigma) \congto \Hom_{M(\O)} (\vonk, \sigma),
  \end{equation}
  where $\hg$ acts on the final term via $\psmg$ (see Lemma~\ref{lm:hecke-parab-ind}). Any vector space automorphism of the
  right-hand side that preserves the left-hand side obviously acts invertibly on it. The claim follows since $\psmg$ is a
  localisation map. The maps in~\eqref{eq:22} are now $\hm$-linear with the right-hand side carrying its natural
  $\hm$-action.
  
  By the above we can pick Hecke eigenvalues $\chi : \hm \to \k$ and a corresponding Hecke eigenvector $f : V \to \pi$ that
  occur in the finite-dimensional $\hm$-module $\Hom_K(V,\pi)$. By~\eqref{eq:22} we can think of $f$ also as $\hm$-eigenvector
  with eigenvalues $\chi$ in $\Hom_{M(\O)} (\vonk, \sigma)$
  and we get a surjection
  \begin{equation*}
  \mind \vonk \otimes_{\hm, \chi} \k \onto \sigma.  
  \end{equation*}
  
  \emph{Case 1: $V$ is $M$-regular.} Then Thm.~\ref{thm:cptind-parabind} applies. Since $\opInd$ is an
  exact functor, we obtain
  \begin{equation}\label{eq:23}
    \kind V \otimes_{\hg,\chi} \k \onto \opInd \sigma.
  \end{equation}
  This map is naturally induced from $V\!\xrightarrow{f} \pi \subset \opInd \sigma$. (See the computation of $\o \theta$ in Step~2 of the
  proof of Thm.~\ref{thm:cptind-parabind}.) As $V$ generates the left-hand side of~\eqref{eq:23}, it follows that $\opInd \sigma$ is
  generated by $f(V)$ as $G$-representation.  In other words, $\pi = \opInd \sigma$.
  
  \emph{Case 2: $V$ is not $M$-regular.} We show by induction that $\pi$ has to contain an $M$-regular weight (in fact the
  one corresponding to $\vnk$ by Lemma~\ref{lm:lift-to-m-reg-wt}). Then Case~1 implies that $\pi = \opInd \sigma$, which
  concludes the proof that $\opInd \sigma$ is irreducible.
  
  Since the derived subgroup of $G$ is simply connected, we can write $V \cong F(\nu)$. As $\Stab_W(\nu) \not\subset W_M$ and
  since the left-hand side is generated by simple reflections, there is a simple root $\alpha \in \Delta-\Delta_M$ such that
  $s_\alpha(\nu) = \nu$. Since the centre of $G$ is connected, there is a fundamental coweight $-\lambda$ associated to
  $\alpha$, which is minuscule as $G = \GL_n$. Just as in \S\ref{sec:minuscule-case} we now put $\nu' = \nu +
  (q-1)\omega_\alpha$ and $V' = F(\nu')$. The map $f$ above gives rise to a non-zero map $\kind V \otimes_{\HH_G,\chi} \k \to
  \pi$ (we also write $\chi$ for $\chi\circ \psmg$). If we can show that Cor.~\ref{cor:minu-case-criterion} applies, there is
  a non-zero map $\kind V' \otimes_{\HH_G,\chi} \k \to \pi$, in particular $V'$ is a weight of $\pi$.  Since
  $\Stab_W(\nu')$ is strictly smaller than $\Stab_W(\nu)$, by induction we eventually find that $\pi$ has to contain an
  $M$-regular weight.
  
  Finally we show that the criterion in Thm.~\ref{thm:irred-parab-ind-gln} implies that the assumptions in
  Cor.~\ref{cor:minu-case-criterion} are satisfied. By Lemma~\ref{lm:product-decomp} and Def.~\ref{df:supersingular}, the
  Hecke eigenvalues $\chi$ are parameterised by $(M,\prod \chi_i)$, where $\chi_i$ is the central character of $\sigma_i$.
  We have $\alpha \not\in \Delta_M$.  Suppose first that there exists $\beta \in \Delta_M$ that is adjacent to $\alpha$ in
  the Dynkin diagram. In this case $\langle \alpha\dual,\beta\rangle \ne 0$, so $\alpha\dual(\varpi) \not \in Z_M$ and we are
  done. See the left part of Fig.~\ref{fig:levis}. Otherwise, we have $n_i = n_{i+1} = 1$ for the two Levi blocks
  $(M_i,\sigma_i)$ and $(M_{i+1},\sigma_{i+1})$ that $\alpha$ ``separates''. (We work with the diagonal torus and the
  upper triangular Borel.)  See the right part of Fig.~\ref{fig:levis}. In particular $\sigma_i = \chi_i$, $\sigma_{i+1} = \chi_{i+1}$. 
  Therefore on $F\s$ we have $\chi_M \circ \alpha\dual = \chi_i \chi_{i+1}^{-1} \ne 1$ by assumption, so we are done. 
\end{proof}

\begin{figure}[t]
  \centering
  \input{levi.pstex_t}
  \caption{Final part of the proof of Thm.~\ref{thm:irred-parab-ind-gln}: $\alpha$ and $\alpha\dual$.}\label{fig:levis}
\end{figure}

\begin{thm}\label{thm:finite-length}
  Suppose that $G = \GL_n$. Let $P$ be the standard parabolic with Levi $\prod_i \GL_{n_i}$, where $\sum_i n_i = n$.
  Suppose that $\sigma_i$ is an irreducible admissible supersingular representation of $\GL_{n_i}(F)$ for $i = 1$, \ldots,
  $r$. Then $\Ind_{\o P}^G (\sigma_1 \otimes \cdots \otimes \sigma_r)$ has finite length.
\end{thm}

By the transitivity of parabolic induction, we can rewrite such a representation in the form
$\Ind_{\o Q}^G (\tau_1 \otimes \cdots \otimes \tau_s)$ with $Q \supset P$ the standard parabolic with
Levi $\prod_i \GL_{m_i}$, such that for all $i$ either
\begin{itemize}
\item $\tau_i$ is supersingular and $m_i > 1$, or
\item $\tau_i \cong \Ind_{\o {B_i}}^{\GL_{m_i}} \eta_i$ for some $\eta_i : F\s \to \k\s$
\end{itemize}
and such that $\eta_i \ne \eta_{i+1}$ whenever both $\tau_i$ and $\tau_{i+1}$ fall into the second case.  (Here and in the following, we
often write $\eta$ when we really mean $\eta \circ \det$. Moreover $B_i$ denotes the Borel in $\GL_{m_i}$.)

By Cor.~\ref{cor:steinb-irred-constit}, the irreducible constituents of $\Ind_{\o {B_i}}^{\GL_{m_i}} \eta_i$ are of the form
$\Sp{Q_i} \otimes \eta_i$, where $Q_i$ ranges over the standard parabolics of $\GL_{m_i}$.

We have thus reduced the proof of Thm.~\ref{thm:finite-length} to the following generalisation of Thm.~\ref{thm:irred-parab-ind-gln}.
It shows that $\Ind_{\o Q}^G (\tau_1 \otimes \cdots \otimes \tau_s)$ has length $2^{\sum_i d_i}$, where
$d_i = 0$ or $m_i-1$ depending on whether $\tau_i$ falls into the first or into the second case.

\begin{thm}\label{thm:irred-parab-ind-gln-2}
  Suppose that $G = \GL_n$. Let $P$ be the standard parabolic with Levi $\prod_i \GL_{n_i}$, where $\sum_i n_i = n$.
  Suppose that $\sigma_i$ is an irreducible admissible representation of $\GL_{n_i}(F)$ for $i = 1$, \ldots,
  $r$ such that for all $i$ either
  \begin{itemize}
  \item $\sigma_i$ is supersingular and $n_i > 1$, or
  \item $\sigma_i \cong \Sp{Q_i} \otimes \eta_i$ for some $\eta_i$ and some standard parabolic $Q_i \subset \GL_{n_i}$.
  \end{itemize}
  Assume that $\eta_i \ne \eta_{i+1}$ whenever both $\sigma_i$ and $\sigma_{i+1}$ fall into the second case.
  Then $\Ind_{\o P}^G (\sigma_1 \otimes \cdots \otimes \sigma_r)$ is irreducible.
\end{thm}

\begin{proof}
  The proof differs from the one of Thm.~\ref{thm:irred-parab-ind-gln} only in the final part. This time the Hecke eigenvalues are parameterised
  by the pair $(M' = \prod M_i',\prod \chi_i)$, where $M_i' = M_i$ and $\chi_i$ is the central character of $\sigma_i$ if $\sigma_i$ falls
  into the first case, whereas $M_i'$ is the torus of $M_i$ and $\chi_i = \eta_i \circ \det$ if $\sigma_i$ falls
  into the second case. It follows that $\alpha\dual(\varpi) \not\in Z_{M'}$ unless the blocks adjacent to $\alpha$, $(M_i,\sigma_i)$ and
  $(M_{i+1},\sigma_{i+1})$, both fall into the second case and the argument goes through since $\eta_i \eta_{i+1}^{-1} \ne 1$.
\end{proof}

\begin{thm}\label{thm:decomp-parab-ind-gln}
  Suppose that the final condition on the $\eta_i$ is dropped in Thm.~\ref{thm:irred-parab-ind-gln-2}.  Then $\Ind_{\o P}^G
  (\sigma_1 \otimes \cdots \otimes \sigma_r)$ is of finite length with explicit irreducible constituents, each occurring with multiplicity one.
\end{thm}

In the proof we describe the irreducible constituents; in particular we show that the length of $\Ind_{\o P}^G (\sigma_1
\otimes \cdots \otimes \sigma_r)$ equals $2^\delta$, where $\delta$ is the number of times the condition $\eta_i \ne
\eta_{i+1}$ fails.

\begin{proof}
  Note that $\bigotimes_{i=1}^s \Sp {Q_i} \cong \Sp {Q_1 \times \cdots \times Q_s}$ as representation of $\prod_{i=1}^s
  \GL_{n_i}$. (For example, directly from~\eqref{eq:26}.)
  We may thus collect consecutive $\sigma_i$ with the same twisting character $\eta_i$ and use transitivity of
  parabolic induction to rewrite our representation as $\Ind_{\o Q}^G (\tau_1 \otimes \cdots \otimes \tau_s)$, where $Q \supset
  P$ is the standard parabolic with Levi $\prod_i \GL_{m_i}$, such that for all $i$ either
  \begin{itemize}
  \item $\tau_i$ is supersingular and $m_i > 1$, or
  \item $\tau_i \cong \Ind_{\o R_i}^{\GL_{m_i}} \Sp {S_i} \otimes \eta'_i$ for some $\eta'_i : F\s \to \k\s$
  \end{itemize}
  and such that $\eta'_i \ne \eta'_{i+1}$ whenever both $\tau_i$ and $\tau_{i+1}$ fall into the second case.  Here $R_i = L_i
  N_i'$ is a standard parabolic of $\GL_{m_i}$ and $S_i$ a standard parabolic of $L_i$.  (Note that $L_i$ is a product of
  consecutive $\GL_{n_j}$ and $S_i$ a product of consecutive $Q_j$.)  By exactness of parabolic induction, by
  Prop.~\ref{prop:opind-spq} and Thm.~\ref{thm:irred-parab-ind-gln-2} it follows that the irreducible constituents are
  obtained by replacing each $\tau_i$ that falls into the second case by $\Sp {S_i'} \otimes \eta'_i$ for any standard
  parabolic $S_i'$ of $\GL_{m_i}$ such that $S_i' \cap L_i = S_i$.
  
  The multiplicity one assertion is true because these irreducible representations do not share any common weight. (This
  holds for the generalised Steinberg representations by the results of Section~\ref{sec:gen-steinberg}.  Now
  use~\eqref{eq:6}.) Alternatively this follows using ordinary parts; see Cor.~\ref{cor:uniqueness}.
\end{proof}

\subsection{General results}
\label{sec:parab-ind-general}

We show that parabolic inductions of irreducible admissible representations are usually irreducible, even for general $G$.
As we do not have the results of \S\ref{sec:minuscule-case} available in general, we need to put stronger hypotheses on the
weights that are allowed to occur. On the other hand, the representation of the Levi does not have to be supersingular.

\begin{thm}\label{thm:irred-parab-ind}
  Let $P = MN$ be a standard parabolic and suppose that $\sigma$ is an irreducible admissible $M$-representation satisfying:
  \begin{enumerate}
  \item[($*$)] For all simple roots $\alpha \in \Delta - \Delta_M$, the restriction to $k\s$ via the coroot $\alpha\dual$
    of the $T(k)$-representation $(\soc_{M(\O)} \sigma)_{(\o U\cap M)(k)}$ does not contain the trivial representation.
  \end{enumerate}
  Then $\opInd \sigma$ is irreducible.
\end{thm}

\begin{proof}
  This follows by the proof of Thm.~\ref{thm:irred-parab-ind-gln}, since we now show that ($*$) implies that all
  weights $V$ of $\opInd \sigma$ are $M$-regular. First assume that the derived subgroup of $G$ is simply
  connected.  In this case, $V \cong F(\nu)$ for some $q$-restricted weight $\nu$ and $\vonk\cong F^M(\nu)$. 
  Assumption ($*$) shows that for all simple roots $\alpha \in \Delta - \Delta_M$, $\langle \nu,\alpha\dual\rangle \not\equiv
  0 \pmod{q-1}$ (as $\vonk$ is a direct summand of $\soc_{M(\O)} \sigma$ and $(\vonk)_{(\o U\cap M)(k)}$ is isomorphic to
  $\nu|_{T(k)}$). In particular, $\langle \nu,\alpha\dual\rangle > 0$ for such $\alpha$. Since $\Stab_W(\nu)$ is generated by simple
  reflections, we see that $\Stab_W(\nu) \subset W_M$, so $V$ is indeed $M$-regular.
  
  If $G$ is general, pick a weight $V$ of $\opInd \sigma$. The restriction of the $T(k)$-representation $\vouk$
  to $k\s$ via $\alpha\dual$ is non-trivial for $\alpha \in \Delta-\Delta_M$ by ($*$).
  Now take a $z$-extension $\wt G \onto G$ of the special fibre, just as in \cite[Lemma~2.5]{bib:satake}. 
  Since coroots are compatible in $z$-extensions, we get the analogous statement for
  $V$ as $\wt G(k)$-representation. By the previous paragraph, $V$ is $\wt M$-regular and therefore $M$-regular.
\end{proof}

\begin{rk}
  Condition ($*$) is best possible in the sense that it is equivalent to the condition that all weights of $\opInd \sigma$
  are $M$-regular.
\end{rk}

Here is a simple application to principal series representations.

\begin{coroll}
  Suppose that $\chi : T \to \k\s$ is a smooth character. Then the principal series representation $\obInd \chi$ is
  irreducible provided that $\chi \circ \alpha\dual|_{k\s} \ne 1$ for all simple roots~$\alpha$.
  \(Note that $\chi|_{T(\O)}$ factors through $T(k)$.\)
\end{coroll}

\begin{proof}
  This follows from the theorem using~\eqref{eq:6}.
\end{proof}

\section{Classification results}
\label{sec:classification}

\subsection{The case of $\GL_n$}
\label{sec:class-case-gl_n}

In this subsection, $G = \GL_n$. The goal is to show that all irreducible admissible representations are of the form 
as in Thm.~\ref{thm:irred-parab-ind-gln-2}. Our strategy is roughly speaking as follows. If $\pi$ is irreducible admissible,
it contains a weight $V$ and Hecke eigenvalues $\chi$ in that weight, so $\kind V \otimes_{\hg,\chi} \k \onto \pi$.
If there is a standard parabolic $P = MN \ne G$ such that $V$ is $M$-regular and $\chi$ factors through $\psmg$, then we can
apply Thm.~\ref{thm:cptind-parabind} to express $\pi$ as a quotient of a (non-trivial) parabolic induction. We can then
use Emerton's ordinary parts functor (a right adjoint to parabolic induction) to write $\pi$ as a quotient of $\opInd \sigma$
with $\sigma$ irreducible admissible, and we can induct. If the process gets stuck, we try to ``change the weight'' (Cor.~\ref{cor:minu-case-criterion}).
The one crucial ingredient that is still needed is that we can detect the trivial representation by its weight and Hecke eigenvalues;
see Prop.~\ref{prop:triv-wt-hecke-evals} below.

We recall that Emerton's functor $\Ord_P$ \cite{bib:emerton-ordinary1} sends smooth $G$-representations to smooth $M$-representations.
It is left exact and preserves admissibility. It is a right adjoint to $\opInd$ between the full subcategories of admissible representations.

\begin{prop}\label{prop:triv-wt-hecke-evals}
  Suppose $\pi$ is a smooth $G$-representation such that $\pi$ contains the trivial
  weight with Hecke eigenvalues parameterised by $(T,1)$. Then either $\pi$ contains the trivial $G$-representation
  or there exists a proper standard parabolic $P$ such that $\Ord_P \pi \ne 0$.
\end{prop}

In fact we will see in Cor.~\ref{cor:characterise-steinberg} below that if $\pi$ is moreover irreducible and admissible, then
$\pi$ is the trivial $G$-representation.
(Note also that $\Ord_P 1 = 0$ for any proper parabolic $P$.)

\begin{proof}
  The first assumption means that $\pi^K \ne 0$.  We use the more classical notation $[KgK]$ for a Hecke operator in the
  unramified Hecke algebra $\HH_G(1)$.  The action on $\pi^K$ is the natural left action: if $KgK = \coprod g_\alpha K$, then $[KgK]v
  := \sum g_\alpha v \in \pi^K$. Note that $[KgK] = 1_{Kg^{-1}K} \in \HH_G(1)$.
  Let $T_i := [Kt_iK]$, where $t_i = \diag(\varpi,\ldots,\varpi,1,\ldots,1)$ ($i$ copies of $\varpi$).
  
  Pick $v \in \pi^K$, a Hecke eigenvector with Hecke eigenvalues parameterised by $(T,1)$. This means that $T_iv = v$
  for all $1 \le i \le n$. (This follows from Cor.~\ref{cor:param-hecke-evals}. Note that $\ge_\R$ is the same as $\ge$ on $X_*(T)$ for
  $\GL_n$ and so by~\eqref{eq:10}, $\SS_G(T_\lambda) = \tau_\lambda$ whenever $-\lambda$ is minuscule.)
  From $i = n$ or from Lemma~\ref{lm:hecke-evals-supersing} we see that $\pi$ has trivial central character.

  Suppose that $\Ord_P \pi = 0$ for all proper standard parabolics $P$.

  \emph{Step 1.} We show that $U_i := [I t_i I] \in \HH_I(1)$ acts nilpotently on $v$ for $i = 1$, \dots, $n-1$. (Recall that $I$ denotes the Iwahori subgroup.)

  Let $P$ be the standard parabolic subgroup with Levi $\GL_i \times \GL_{n-i}$ and let $\P$ be the corresponding parahoric subgroup.
  Suppose that the Hecke operator $[\P t_i \P]$ has a non-zero eigenvalue on $\pi^{\P}$. Then a corresponding eigenvector
  is even stable under $\HH(1) \subset \HH_\P(1)$ because $\HH(1)$ is generated by $E_{h} = [\P h^{-1} \P]$ for $h = t_i^{-1}$ and $h = t_n^{-1}$
  (in the notation of Lemma~\ref{lm:hecke-compat-m-p}). Then there is a unique unramified character $\chi_M : Z_M \to \k\s$ that sends
  $h^{-1}$ to the (non-zero!) eigenvalue of $E_h$ for $h = t_i^{-1}$ and $h = t_n^{-1}$. Since $E_h$ is identified with $T_h^M \in \HH_M(1)$,
  Prop.~\ref{prop:weights-of-ord-p} applies with $\o V = 1$ and we see that $\Ord_P\pi\ne 0$. Contradiction.
  
  Thus $[\P t_i \P]$ is nilpotent on $\pi^\P$. But $[\P t_i \P] = [I t_i I]$ on $\pi^{\P} \subset \pi^I$, so $U_i$ acts
  nilpotently on $v$. (To see this, note that $\P t_i \P/\P$ naturally bijects with $N(\O)/{}^{t_i} N(\O)$ by
  Lemma~\ref{lm:iwahori-decomp}. Similarly $I t_i I/I$ naturally bijects with $U(\O)/{}^{t_i} U(\O)$. Now use
  that $U = N \rtimes (M \cap U)$.)

  \emph{Step 2. Setup and some lemmas.}

  Define the following elements in the Iwahori Hecke algebra $\HH_I(1)$:
  \begin{equation*}
    S_1 = \bigg[I\left(\begin{smallmatrix} &1\\
      1\\&&\varddots \\ &&&1\end{smallmatrix}\right)I\bigg], \dots, S_{n-1} = \bigg[I\left(\begin{smallmatrix} 1\\ &\varddots
      \\&&&1\\&&1\end{smallmatrix}\right)I\bigg], \Pi = \bigg[I\left(\begin{smallmatrix} &1\\ &&\varddots\\
      &&&1\\\varpi\end{smallmatrix}\right)I\bigg],
  \end{equation*}
  where $S_i$ is defined by the $i$-th simple reflection. It is straightforward to verify the following relations.
  \begin{alignat*}{2}
    S_i^2 &= -S_i  &&\text{for all $i$}, \\
    S_iS_j &= S_jS_i  &&\text{whenever $|i-j| > 1$}, \\
    S_kS_{k+1}S_k &= S_{k+1}S_kS_{k+1} \quad && \text{for all $k < n-1$}, \\
    S_k \Pi &= \Pi S_{k+1} && \text{for all $k < n-1$}. \\
    \noalign{\noindent We also have}
    \Pi^n &= 1\ \text{on}\ \pi^I.
  \end{alignat*}
  The quadratic relations take a simple form, since we work in characteristic~$p$.  For the last property note that the
  matrix defining $\Pi$ normalises $I$. (We remark that the above is a presentation of $\HH_I(1)$, but we will not need that fact.
  It can be deduced from Cor.~4 in~\cite{bib:Vig-ss} for $R = \fpb$ and $\gamma$ consisting of the trivial character.)

  Note that
  \begin{equation*}
    T_1 v = \sum \left(\begin{smallmatrix} \varpi &a_2&a_3&\cdots\\ & 1\\&&1\\&&&\varddots\end{smallmatrix}\right)v
    + \sum \left(\begin{smallmatrix} 1 \\ & \varpi&a_3&\cdots\\&&1\\&&&\varddots\end{smallmatrix}\right)v +
    \cdots,
  \end{equation*}
  where in each sum, the $a_i$ run over representatives of $\O$ modulo $\varpi$.

  It is not hard to express this in terms of the Iwahori Hecke action. First, in each of these sums we can interchange the
  first column with the column containing $\varpi$ (since $v \in \pi^K$). For $1 \le i \le n$ let $x_i \in G$ be the product
  of the matrices defining $S_{i(i+1)\cdots(n-1)}\Pi$. It will suffice to show that $S_{i(i+1)\cdots(n-1)}\Pi = [I x_i I]$ and that
  $I x_i I$ is the disjoint union of the cosets $u(a_{i+1},\dots,a_n) x_i I$, where $u(a_{i+1},\dots,a_n) \in U(\O)$ is
  zero above the diagonal except for the $i$-th row which ends in $(a_{i+1},\dots,a_n)$, and the $a_j$ run
  through representatives of $\O$ modulo $\varpi$. The second claim is a simple calculation. The first claim follows
  since both sides contain $q^{n-i}$ one-sided cosets.
  Since $T_1v = v$, we get:
  \begin{equation}
    v = \sum_{i=1}^n S_{i(i+1)\cdots(n-1)} \Pi v, \label{eq:16}
  \end{equation}
  where we abbreviate $S_{i(i+1)\cdots (j-1)} := S_i S_{i+1}\cdots S_{j-1}$ for any $1 \le i \le j \le n$.
  (This equals 1 when $i = j$.)

  \begin{sublm}\label{lm:commutation}
    Suppose $i \le k \le \ell \le j$. Then
    \begin{equation*}
      S_{i(i+1)\cdots j}S_{k(k+1)\cdots (\ell-1)} = S_{(k+1)(k+2)\cdots \ell}S_{i(i+1)\cdots j}
    \end{equation*}
  \end{sublm}

  \begin{proof}
    It suffices to consider the case $\ell-1 = k$. In that case, $i \le k \le j-1$. By the braid relations,
    \begin{align*}
      S_{i(i+1)\cdots j}S_{k} &=  S_{i(i+1)\cdots (k-1)k(k+1)}S_{k}S_{(k+2)\cdots j} \\
      &=  S_{i(i+1)\cdots (k-1)}S_{(k+1)}S_{k(k+1)(k+2)\cdots j}\\
      &=  S_{(k+1)}S_{i(i+1)\cdots j}.
    \end{align*}
  \end{proof}

  \begin{sublm}\label{lm:kill-v}
    We have $S_{n-1}\Pi^2 v = 0$.
  \end{sublm}

  \begin{proof}
    We first note that for all $i$ we have $S_i v = 0$. This is immediate from the fact that $v \in \pi^K$
    and that the double coset defining $S_i$ is a disjoint union of $q = \# k$ one-sided cosets.

    Then we get $S_{n-1}\Pi^2 v = \Pi^2 S_1 v = 0$ (since $\Pi^2 v = \Pi^{-(n-2)} v$).
  \end{proof}

  \begin{sublm}\label{lm:U_i} For all $i$, we have
    $U_i = (S_{i\cdots (n-1)} \Pi)^i$.
  \end{sublm}

  \begin{proof}
    If we multiply out the double cosets on the right-hand side, we certainly find the matrix
    $t_i$. (We just multiply
    the defining matrices of all these double cosets.)
    To see that the product of the double cosets yields only the double coset of this diagonal matrix,
    it's enough to count the number of one-sided cosets. On the left-hand side the number of cosets
    is $q^{i(n-i)}$. Since each $S_i$ contains $q$ cosets and $\Pi$ contains just one, the result follows.
  \end{proof}

  \emph{Step 3. We show that $U_1 v = S_{12\cdots(n-1)} \Pi v = 0$.}

  Sublemma~\ref{lm:U_i} tells us that $U_1 = S_{12\cdots(n-1)} \Pi$. We now see that
  \begin{align*}
    (S_{12\cdots(n-1)} \Pi)^2 v &= S_{12\cdots(n-1)} \Pi \left(v - \sum_{i = 2}^n S_{i\cdots(n-1)}\Pi v\right) \\
    \noalign{\noi by (\ref{eq:16}),}
    &= S_{12\cdots(n-1)} \Pi v - \sum_{i = 2}^n S_{12\cdots(n-1)}S_{(i-1)\cdots(n-2)}\Pi^2 v \\
    \noalign{\noi by pushing $\Pi$ to the right,}
    &= S_{12\cdots(n-1)} \Pi v - \sum_{i = 2}^n S_{i\cdots(n-1)}S_{12\cdots(n-1)}\Pi^2 v \\
    \noalign{\noi by Sublemma~\ref{lm:commutation},}
    &= S_{12\cdots(n-1)} \Pi v,
  \end{align*}
  by Sublemma~\ref{lm:kill-v}. In other words, $U_1^2 v = U_1v$.
  So by Step~1, we have $U_1v = 0$.

  \emph{Step 4. We show that $U_2 v = S_{23\cdots(n-1)} \Pi v = 0$.}

  By Step~3, the first term in (\ref{eq:16}) vanishes and we get
  \begin{equation}
    \label{eq:17}
    v = \sum_{i=2}^n S_{i(i+1)\cdots(n-1)} \Pi v.
  \end{equation}

  We repeat the calculation we did in the previous step, using \eqref{eq:17} instead of \eqref{eq:16}; we find:
  \begin{align*}
    (S_{23\cdots(n-1)} \Pi)^2 v = S_{23\cdots(n-1)} \Pi v.
  \end{align*}

  It follows from Sublemma~\ref{lm:U_i} that
  \begin{equation}
      U_2 v = (S_{23\cdots(n-1)} \Pi)^2 v = S_{23\cdots(n-1)} \Pi v
  \end{equation}
  and moreover that $U_2^2 v = U_2 v$. By Step~1 we see that $U_2 v = 0$.

  \medskip
  $\vdots$
  \medskip

  \emph{Step $n+2$. We show that $v \in \pi^G$.}

  From Step $i+2$, we know that $U_i v = S_{i(i+1)\cdots(n-1)} \Pi v = 0$ ($1 \le i \le n-1$).
  Equation~(\ref{eq:16}) thus simplifies to $v = \Pi v$.
  Since $G$ is generated by $K$ and $\left(\begin{smallmatrix} &1\\ &&\varddots\\
      &&&1\\\varpi\end{smallmatrix}\right) \in N_G(I)$ (for example by the Cartan decomposition, generating $t_i$
  first as in the proof of Sublemma~\ref{lm:U_i}), the claim follows.
  This completes the proof.
\end{proof}

\begin{thm}\label{thm:classification}
  Let $\pi$ be any irreducible admissible $G$-representation.  Then there exists a standard parabolic $P$ with Levi $\prod_{i=1}^r
  \GL_{n_i}$ and irreducible admissible representations $\sigma_i$ of $\GL_{n_i}(F)$ such that $\pi
  \cong \Ind_{\o P}^G (\sigma_1 \otimes \cdots \otimes \sigma_r)$ and such that for all $i$ either
  \begin{itemize}
  \item $\sigma_i$ is supersingular and $n_i > 1$, or
  \item $\sigma_i \cong \Sp{Q_i} \otimes \eta_i$ for some $\eta_i$ and some standard parabolic $Q_i \subset \GL_{n_i}$.
  \end{itemize}
  Moreover $\eta_i \ne \eta_{i+1}$ whenever both $\sigma_i$ and $\sigma_{i+1}$ fall into the second case.
\end{thm}

Together with Thm.~\ref{thm:irred-parab-ind-gln-2}, this achieves the classification of irreducible admissible
$\GL_n(F)$-representations.

We first prove a useful lemma using ordinary parts.

\begin{lm}\label{lm:quotient-of-parab-ind}
  Suppose that $Q = LN'$ is a standard parabolic subgroup. Suppose there is a $G$-invariant surjection $\oqInd \tau \onto
  \pi$ for some smooth $L$-representation $\tau$ having a central character and some irreducible admissible
  $G$-representation~$\pi$. Then there exists an irreducible admissible $L$-representation $\sigma$ and a $G$-invariant
  surjection $\oqInd \sigma \onto \pi$.
\end{lm}

\begin{proof}
  As $\tau$ is locally $Z_L$-finite and $\pi$ is smooth, the natural map
  \begin{equation*}
    \Hom_G(\oqInd \tau,\pi) \to \Hom_L(\tau,\Ord_Q \pi)
  \end{equation*}
  is injective: this follows from the injectivity of the first map in~\cite[(4.4.7)]{bib:emerton-ordinary1} (a
  geometric fact) and the isomorphism in \cite[Cor.~4.2.8]{bib:emerton-ordinary1}.  In particular, $\Ord_Q \pi \ne 0$. As
  $\Ord_Q \pi$ is admissible (\cite[Thm.~3.3.3]{bib:emerton-ordinary1}), it has an irreducible admissible
  subrepresentation $\sigma \INTO \Ord_Q \pi$ (the dual of an admissible $L$-representation is a finitely-generated module
  over the noetherian ring $\k[[L(\O)]]$, see~\cite[Lemma~2.2.11]{bib:emerton-ordinary1}). By the
  adjunction~\cite[Thm.~4.4.6]{bib:emerton-ordinary1} we get a $G$-linear map $\oqInd \sigma \onto \pi$.
\end{proof}

\begin{proof}[Proof of Theorem~\ref{thm:classification}]
  We argue by induction on $n$.  For $n = 1$, there is nothing to show. For $n > 1$, we may assume that there is at least one
  weight $V$ that occurs in $\pi$ with non-supersingular Hecke eigenvalues $\chi$ parameterised by $(M,\chi_M)$
  (otherwise $\pi$ is supersingular and we are done). This means that $M \ne G$.

  Let $\Delta_V = \{\alpha \in \Delta : s_\alpha \in \Stab_W(\vuk) \}$, where $s_\alpha \in W$ denotes the simple reflection
  corresponding to the root $\alpha$. Thus $V$ is $M$-regular if and only if $\Delta_V \subset \Delta_M$.

  \emph{Case 1: $\Delta_V \cup \Delta_M \ne
    \Delta$.} 
  Pick $\alpha \in \Delta-(\Delta_V \cup \Delta_M)$. Let $Q = LN'$ be the standard, maximal parabolic defined by $\alpha$, so
  $\Delta_L = \Delta-\{\alpha\} \supset \Delta_M$.  Then $V$ is $L$-regular and $\chi$ factors through $\HH_L(V_{\o {N'}(k)})$.
  Thm.~\ref{thm:cptind-parabind} gives rises to
  \begin{equation*}
    \oqInd (\lind V_{\o {N'}(k)} \otimes_{\HH_L(V_{\o {N'}(k)}),\chi} \k) \onto \pi.
  \end{equation*}
  Note that the $L$-representation that is being induced on the left-hand side has a central character (namely $h
  \mapsto \chi(T_{h^{-1}}^L)$). Thus by Lemma~\ref{lm:quotient-of-parab-ind} there is an 
  irreducible admissible $L$-representation $\sigma$ and a $G$-linear surjection $\oqInd \sigma \onto \pi$.
  We can now decompose $\sigma$ as a tensor product of irreducible admissible representations of the Levi blocks
  (Lemma~\ref{lm:product-decomp}).  By induction, each of them is of the desired form. By transitivity of parabolic induction
  we see that $\pi$ is a quotient of a parabolic induction of the desired form, except that consecutive $\eta_i$ might be
  equal. But by Thm.~\ref{thm:decomp-parab-ind-gln} we know that all Jordan--H\"older factors of this parabolic induction are of
  the desired form.

  \emph{Case 2: $\Delta_V \cup \Delta_M = \Delta$.} Suppose first that there is an $\alpha \in \Delta - \Delta_M$ such that
  $\alpha\dual(\varpi) \not\in Z_M$ or $\chi_M(\alpha\dual(\varpi)) \ne 1$. Then by Cor.~\ref{cor:minu-case-criterion}, there is
  a weight $V'$ that occurs in $\pi$ with the same Hecke eigenvalues $\chi$ and such that $\Delta_{V'} = \Delta_V - \{\alpha\}$.
  We are thus reduced to Case~1.

  Otherwise, for all $\alpha \in \Delta - \Delta_M$ we have $\alpha\dual(\varpi) \in Z_M$ and $\chi_M(\alpha\dual(\varpi)) =
  1$.  The first condition implies that $M = T$ and then the second condition implies that $\chi_M = \eta \circ \det$, for
  some smooth character $\eta$. By twisting we may assume that $\chi_M = 1$ (see Lemma~\ref{lm:hecke-evals-twist}).  Thus $\vouk$
  is trivial. Since we are in Case~2, we also have $\Delta_V = \Delta$. Putting together these two facts and using
  Lemma~\ref{lm:wts-invts}, we see that $V$ is the trivial weight. By Prop.~\ref{prop:triv-wt-hecke-evals},
  either $\pi$ is trivial (and we are done) or $\Ord_Q \pi \ne 0$ for some proper standard parabolic $Q$ (and we induct
  as in Case~1).
\end{proof}

\begin{coroll}\label{cor:uniqueness}
  Suppose that $\pi$ is an irreducible admissible $G$-representation.

  \begin{enumerate}
  \item All Hecke eigenvalues of $\pi$ are parameterised by the same pair $(M',\chi_{M'})$.
  \item There is a unique datum $(P,(\sigma_i)_{i=1}^r)$ as in Thm.~\ref{thm:classification} such that $\pi \cong \opInd
    (\sigma_1 \otimes \cdots \otimes \sigma_r)$.  In other words, there are no non-trivial intertwinings between the
    representations in Thm.~\ref{thm:classification}.
  \end{enumerate}
\end{coroll}

\begin{proof}
  By Thm.~\ref{thm:classification} we know that $\pi \cong \opInd
    (\sigma_1 \otimes \cdots \otimes \sigma_r)$ with $(P = MN,(\sigma_i)_{i=1}^r)$ satisfying the conditions in 
    that theorem.

  As in Thm.~\ref{thm:irred-parab-ind-gln-2} we write $M = \prod M_i$ as product of Levi blocks.
  By Lemmas~\ref{lm:hecke-evals-parab-ind} and~\ref{lm:product-decomp} we see that the Hecke eigenvalues
  of $\pi$ are parameterised by $(\prod M_i',\prod \chi_i)$, where $(M_i',\chi_i)$ runs over the Hecke eigenvalues of
  weights of $\sigma_i$. If $\sigma_i$ is supersingular, we have $M_i' = M_i$ and $\chi_i$ is the central
  character of $\sigma_i$ (see Def.~\ref{df:supersingular}).
  If $\sigma_i \cong \Sp {Q_i} \otimes \eta_i$, then $M_i'$ is the torus in $M_i$ and $\chi_i = \eta_i \circ \det$
  (use Prop.~\ref{prop:steinb-wts} and Lemma~\ref{lm:hecke-evals-twist}). Part (i) follows.

  Since consecutive $\eta_i$ are distinct, these common Hecke eigenvalues determine $M$. Part (ii) follows since
  $\Ord_P(\opInd \sigma) \cong \sigma$ \cite[Prop.~4.3.4]{bib:emerton-ordinary1}.
\end{proof}

We can now prove a converse to Prop.~\ref{prop:steinb-wts}.

\begin{coroll}\label{cor:characterise-steinberg}
  Suppose that $\pi$ is an irreducible admissible $G$-representation. Suppose that $Q$ is a standard parabolic
  and that the weight $V_Q$ occurs in $\pi$ with Hecke eigenvalues parameterised by $(T,1)$. Then
  $\pi \cong \Sp Q$.
\end{coroll}

\begin{proof}
  If $\pi \cong \opInd (\sigma_1 \otimes \cdots \otimes \sigma_r)$ as in Thm.~\ref{thm:classification}, then the analysis of
  Hecke eigenvalues in the proof of Cor.~\ref{cor:uniqueness} shows that $P = G$, and that $\sigma_1$ is a
  generalised Steinberg representation (note that $r = 1$). Then Prop.~\ref{prop:steinb-wts} show that
  $\pi \cong \Sp Q$.
\end{proof}

\begin{df}\label{df:supercusp}
  Suppose that $\pi$ is an irreducible admissible $G$-representation. We say that $\pi$ is \emph{supercuspidal} if it does
  not occur as subquotient in $\opInd \sigma$, where $P = MN$ is any proper standard parabolic and $\sigma$ any irreducible
  admissible $M$-representation.
\end{df}

\begin{coroll}\label{cor:supercuspidal}
  Suppose $\pi$ is an irreducible admissible $G$-representation.
  Suppose $P = MN$ is a standard parabolic and $\sigma$ an irreducible admissible $M$-representation.
  \begin{enumerate}
  \item $\opInd \sigma$ is of finite length, and all constituents occur with multiplicity one. All Hecke eigenvalues of all constituents are parameterised
    by the same pair $(M',\chi_{M'})$.
  \item $\pi$ is supersingular if and only if $\pi$ is supercuspidal.
  \end{enumerate}
\end{coroll}

\begin{proof}
  Part (i) follows from Thm.~\ref{thm:decomp-parab-ind-gln} and the analysis of Hecke eigenvalues in the proof of Cor.~\ref{cor:uniqueness}.
  (In terms of the classification in Thm.~\ref{thm:classification}, only the parabolics $Q_i$ differ among the constituents. But $Q_i$ plays
  no role for the Hecke eigenvalues.) It is clear that the Hecke eigenvalues of $\sigma$ are parameterised by the same pair
  $(M',\chi_{M'})$. In particular, $M' \subset M$.
  
  If $\pi$ is supercuspidal, it follows from Thm.~\ref{thm:classification} and the definition of the generalised Steinberg
  representations that $\pi$ is supersingular.  Conversely, suppose that $\pi$ occurs in $\opInd \sigma$ for some proper
  parabolic $P$.  By the above, the Hecke eigenvalues of $\pi$ are parameterised by $(M',\chi_{M'})$ with $M' \subset M
  \ne G$.  Thus $\pi$ cannot be supersingular.
\end{proof}

\subsection{Some general results}
\label{sec:class-general-results}

We show for general $G$ that any irreducible admissible $G$-representation is parabolically induced from a supersingular
representation, provided that it does not contain certain weights at the boundary.
First we need the following result on the weights of ordinary parts.

\begin{prop}\label{prop:weights-of-ord-p}
  Suppose that $\o V$ is a weight for $M$.  Suppose $\chi_M : Z_M \to \k\s$ is a homomorphism such that
  $\chi_M|_{Z_M(\O)}$ is the central character of $\o V$.  Then there is an algebra homomorphism $\chi : \HH_M(\o V) \to \k$ such
  that $\chi(T_h^M) = \chi_M(h)^{-1}$ for all $h \in Z_M$ and for any smooth $G$-representation $\pi$ we have a natural
  isomorphism
  \begin{equation}\label{eq:12}
    \Hom_{M(\O)} (\o V, (\Ord_P \pi)^{Z_M = \chi_M}) \cong \Hom_G (\pind \o V \otimes_{\HH(\o V),\chi} \k,\pi).
  \end{equation}
\end{prop}

\begin{proof}
  We let $h \in Z_M$ act on $\o V$ by $\chi_M(h)$. As this agrees with the usual action on $Z_M(\O)$, the weight $\o V$ becomes an $M(\O)Z_M$-representation.
  Thus by the definition of $\Ord_P \pi$ \cite[Def.~3.1.9]{bib:emerton-ordinary1}, the left-hand side of \eqref{eq:12} is naturally isomorphic to
  \begin{equation*}
    \Hom_{M(\O)Z_M} (\o V, \Hom_{Z_M^+}(Z_M,\pi^{N(\O)})) \cong \Hom_{M(\O)Z_M^+} (\o V, \pi^{N(\O)}).
  \end{equation*}
  (Even though we let $Z_M$ denote the \emph{connected} centre, as opposed to \cite{bib:emerton-ordinary1}, it is straightforward to verify
  that this is irrelevant for the definition of $\Ord_P$.)
  We verify that any $M(\O)Z_M^+$-linear map $f : \o V\to \pi^{N(\O)}$ factors through $\pi^{\P(1)}$, where
  $\P(1) = \ker(\P\to M(k))$. By the smoothness of $\pi$ and by Lemma~\ref{lm:iwahori-decomp}, there is an $h\in Z_M^{--}$ such that
  $f(\o V)$ is fixed by ${}^h (\P^-)$. Then $h^{-1} f(\o V)$ is fixed by $\P^- (M\cap K(1)) (\P^+)^h = \P(1) \cap \P(1)^h$.
  Note that the natural map $N(\O)/N(\O)^h \to \P(1)/(\P(1) \cap \P(1)^h)$ is a bijection. By the definition of the Hecke $Z_M^+$-action
  on $\pi^{N(\O)}$, we have for all $\o v \in \o V$:
  \begin{equation*}
    \chi_M(h^{-1}) f(\o v) = \sum_{N(\O)/N(\O)^h} n h^{-1} f(\o v) = \sum_{\P(1)/(\P(1) \cap \P(1)^h)} \p h^{-1} f(\o v).
  \end{equation*}
  This implies that $f(\o v) \in \pi^{\P(1)}$. In particular $f$ is $\P$-linear.
  
  By Prop.~\ref{prop:param-hecke-evals} and Cor.~\ref{cor:param-hecke-evals} (or directly) the pair $(M,\chi_M)$ gives rise
  to $\chi : \HH_M(\o V) \to \k$ such that $\chi(T_h^M) = \chi_M(h)^{-1}$ for $h \in Z_M$.
  We verify that the $\P$-linear map $f : \o V \to \pi^{\P(1)} \INTO \pi$ is an $\HH(\o V)$-eigenvector with eigenvalues
  $\chi$. For $h \in Z_M^-$,
  \begin{equation*}
    (f \ast E_h)(\o v) = \sum_{\P/(\P \cap \P^h)} \p h^{-1} f(E_h(h\p^{-1}) \o v) = \sum_{N(\O)/N(\O)^h} n h^{-1} f(\o v) = \chi_M(h^{-1}) f(\o v),
  \end{equation*}
  since $N(\O)/N(\O)^h \to \P/(\P \cap \P^h)$ is a bijection and by definition of the $Z_M^+$-Hecke action on $\pi^{N(\O)}$.
  (Note that $E_h$ was defined in the proof of Lemma~\ref{lm:hecke-compat-m-p}.) Now note that $\chi(E_h) = \chi(T_h^M) = \chi_M(h^{-1})$.

  Thus we have a natural map
  \begin{equation*}
    \Hom_{M(\O)Z_M^+} (\o V, \pi^{N(\O)}) \to \Hom_G (\pind \o V \otimes_{\HH(\o V),\chi} \k,\pi).
  \end{equation*}
  Reversing the above argument we obtain an inverse map. This completes the proof.
\end{proof}

\begin{lm}\label{lm:irr-rep-as-quot-of-ss}
  Suppose that $\pi$ is an irreducible admissible $G$-representation.  Then there exists a standard parabolic $P = MN$ and
  an irreducible admissible $M$-representation~$\sigma$ such that $\opInd \sigma \onto \pi$ and such that for all triples
  $(K'_M,T'_M,B'_M)$ as in Def.~\ref{df:supersingular}, all $T'_M$-regular $K'_M$-weights of $\sigma$ are supersingular.
\end{lm}

\begin{proof}
  Suppose that $\pi$ contains a $T'$-regular $K'$-weight that is non-supersingular, for some triple $(K',T',B')$.  Relabel
  $(K',T',B')$ as $(K,T,B)$ for a moment. Pick a $T$-regular weight~$V$ that occurs in~$\pi$ with Hecke eigenvalues
  parameterised by the pair $(M,\chi_M)$, $M \ne G$.  In the notation of the proof of Thm.~\ref{thm:classification} we see
  that $\Delta_V = \varnothing$ and $\Delta_M \ne \Delta$.  By Case~1 of the proof of Thm.~\ref{thm:classification}, there is
  a \emph{proper} standard parabolic $Q = LN'$, an irreducible admissible $L$-representation $\sigma$, and a $G$-linear
  surjection $\oqInd \sigma \onto \pi$. When we change back to the original notation, $\o Q$ gets replaced by some parabolic
  subgroup, but we can easily rewrite any parabolic induction as an induction from a standard parabolic subgroup. Now induct.
\end{proof}

\begin{thm}\label{thm:irr-rep-as-ind-of-ss}
  Suppose that $\pi$ is an irreducible admissible $G$-representation such that the following condition is satisfied:
  \begin{enumerate}
  \item[($*'$)] For all triples $(K',T',B')$ as in Def.~\ref{df:supersingular} and for all simple roots $\alpha \in X^*(T')$, the restriction to
    $k\s$ via the coroot $\alpha\dual$ of the $T'(k)$-representation $(\soc_{K'} \pi)^{U'(k)}$ does not contain the trivial
    representation.
  \end{enumerate}
  Then there exists a standard parabolic $P = MN$
  and an irreducible admissible supersingular $M$-representation~$\sigma$ such that $\pi \cong \opInd \sigma$.
\end{thm}

In the situation of the theorem all Hecke eigenvalues of $\pi$ can be naturally identified.
We can compare Hecke eigenvalues in the following way. Suppose we work with a triple $(K',T',B')$.  The Hecke eigenvalues of
a $K'$-weight are then parameterised by a pair $(L',\chi_{L'})$, where $L'$ is a standard Levi with respect to $(T',B')$. The
pair $(T',B')$ can be conjugated to $(T,B)$ by an element of~$G$, which is unique up to $T$. After conjugating $(L',\chi_{L'})$
by such an element we get a well-defined pair $(L,\chi_L)$ with $L$ standard and $\chi_L : Z_L \to \k\s$ a smooth character. (To
be precise we mean $L_{/F}$: the integral structure depends on the triple we started with.)  It is now clear from Def.~\ref{df:supersingular}
and Lemma~\ref{lm:hecke-evals-parab-ind}
that in this way all Hecke eigenvalues of $\pi$ give rise to the pair $(M,\omega_\sigma)$, where $\omega_\sigma$ is the
central character of $\sigma$ (just as in Cor.~\ref{cor:uniqueness}).

\begin{proof}
  By Thm.~\ref{lm:irr-rep-as-quot-of-ss} there exists an irreducible admissible $M$-representation~$\sigma$ such that
  $\opInd \sigma \onto \pi$ and such that all $T'_M$-regular $K'_M$-weights of $\sigma$ are supersingular. It will suffice to
  show that condition ($*'$) implies that $\opInd \sigma$ is irreducible and that all $K'_M$-weights of $\sigma$ are
  $T'_M$-regular (for any triple $(K'_M,T'_M,B'_M)$).
  
  We first show that for any weight $\o V$ of $\sigma$ the unique $M$-regular weight $V$ of $G$ such that $\vnk \cong \o V$
  (see Lemma~\ref{lm:lift-to-m-reg-wt}) is a weight of $\pi$.  Let $\chi_M : Z_M \to \k\s$ be the central character of
  $\sigma$, so we have an $M(\O)$-linear map $\o V \INTO \sigma \INTO (\Ord_P \pi)^{Z_M = \chi_M}$.  By
  Prop.~\ref{prop:weights-of-ord-p} we obtain an algebra homomorphism $\chi : \HH_M(\o V) \to \k$ and a $G$-linear surjection
  $\pind \o V \otimes_{\HH(\o V),\chi} \k \onto \pi$. As $V$ is $M$-regular, we may apply
  Cor.~\ref{cor:parahoric-and-spherical-inductions} and get a $G$-linear surjection $\ind_K^G V \otimes_{\HH(\o V),\chi} \k
  \onto \pi$. Thus $V$ is a weight of~$\pi$.

  To show that $\opInd \sigma$ is irreducible it suffices to show that $\sigma$ satisfies condition ($*$) of
  Thm.~\ref{thm:irred-parab-ind}.  Suppose that $\o V$ is a weight of $\sigma$. By what we just showed, $\o V \cong \vnk$ for
  some weight $V$ of $\pi$.  So $\o V_{(\o U\cap M)(k)} \cong \o V^{(U \cap M)(k)} \cong V^{U(k)}$. Thus the condition
  follows from ($*'$).

  To complete the proof, it will be useful to have an alternative description of a triple $(K,T,B)$.
  We claim that $K$ is a hyperspecial subgroup corresponding to a point in the apartment of $T_{/F}$. Conversely we claim
  that if we start with a hyperspecial subgroup $K'$ corresponding to a point in the apartment of a torus $T'_{/F}$ and a Borel
  subgroup $B'_{/F}$ containing $T'_{/F}$, then $K' = G'(\O)$ for a unique reductive integral structure $G'$
  of $G_{/F}$. Moreover, $T'_{/F}$ extends (uniquely) to a split maximal torus $T'$ of $G'$ and $B'_{/F}$ extends (uniquely) to a Borel
  subgroup of $G'$ containing $T'$.
  
  To see the first claim, we already recalled in Section~\ref{sec:hecke-eigenv} that the reductive integral structure $G$ of $G_{/F}$ corresponds
  to a unique hyperspecial point $x$ in the building of~$G_{/F}$ (the unique fixed point of~$K$). Choose a torus ${T_x}_{/F}$ in~$G_{/F}$ 
  whose apartment contains~$x$. By Bruhat--Tits theory, the torus ${T_x}_{/F}$ extends (uniquely) to a split maximal
  torus $T_x$ of~$G$. But then $T$ and~$T_x$ are conjugate by an element of $G(\O) = K$ \cite[Exp.~XXVI,
  Prop.~6.16]{bib:SGA3}, so $x$ also lies in the apartment of $T_{/F}$. For the second claim, the corresponding hyperspecial point gives
  rise to the required reductive integral structure $G'$ of $G_{/F}$. The remaining assertions are established by Bruhat--Tits
  (\cite[\S II.4.6]{bib:BT2}; see \cite[\S 3]{bib:satake} for more details).
  
  Finally we deduce that all $K'_M$-weights of $\sigma$ are $T'_M$-regular.  By conjugating $(K'_M,T'_M,B'_M)$ by an element
  of~$M$, we may assume that $(T'_M, B'_M)$ agrees with $(T,B\cap M)$ on the generic fibre. We claim that there is a
  reductive integral structure $G'$ of $G_{/F}$ with a torus $T'$ and a Borel subgroup $B'$ containing it such that $(T',B')$
  agrees with $(T,B)$ on the generic fibre and such that $G'(\O) \cap M = K'_M$. By the above claims it suffices to find a
  hyperspecial subgroup $K'$ corresponding to a point in the apartment of $T_{/F}$ such that $K' \cap M = K'_M$. By \cite[\S
  7.6]{bib:BT1} and \cite[II.4.2.18]{bib:BT2} the extended building of $M_{/F}$ can be embedded $M$-equivariantly in the
  extended building of $G_{/F}$ such that the apartments $A_M(T_{/F})$ and $A_G(T_{/F})$ of~$T_{/F}$ coincide. Moreover any
  wall in $A_M(T_{/F})$ is a wall in $A_G(T_{/F})$. Let $y$ be a (hyperspecial) point in $A_M(T_{/F})$ whose stabiliser in
  $M$ is $K'_M$. (Note that $y$ need not be unique, since we work with the extended building.) As $y$ is hyperspecial, for each
  simple root $\alpha \in \Delta_M$ there is an $\alpha$-wall $H_\alpha$ that passes through $y$. The stabiliser in $M$ of
  any point in $\bigcap_{\alpha \in \Delta_M} H_\alpha$ is $K'_M$. (It projects to a single point in the reduced building of $M_{/F}$.)
  For each simple root $\beta \in \Delta-\Delta_M$ pick an arbitrary $\beta$-wall $H_\beta$ in $A_G(T_{/F})$. As simple roots
  are linearly independent, we see that $\bigcap_{\beta \in \Delta} H_\beta \ne \varnothing$. Let $x$ be any point in this
  intersection. It is hyperspecial and its stabiliser $K'$ in $G$ is as required.

  After relabelling $(K',T',B')$ as $(K,T,B)$ we can now assume that
  $(K'_M,T'_M,B'_M) = (M(\O),T,B\cap M)$. Then notice that we established above that condition ($*$) holds for $\sigma$,
  but even for \emph{all} simple roots $\alpha$. We show that the condition for simple roots $\alpha \in \Delta_M$ implies
  the $T$-regularity of all weights of $\sigma$. If the derived subgroup of~$M$ is simply connected, we can write $\o V \cong F^M(\nu)$ for some
  $\nu \in X^*(T)$. If $s_\alpha(\nu) = \nu$ for some $\alpha \in \Delta_M$ it follows that $\nu \circ \alpha\dual$ is trivial on
  $k\s$, contradicting the condition. In the general case one uses a $z$-extension of~$M_{/k}$.
\end{proof}

\section{Submodule structure}
\label{sec:submodule-structure}

In this section we determine the submodule structure of parabolically induced representations in two situations: for the trivial principal series
($G$ general) and for the parabolic induction of an irreducible admissible representation (when $G = \GL_n$).

We recall some elementary facts about the submodule structure of finite length modules whose constituents occur with multiplicity one.
Let $M$ be such a module over some, not necessarily commutative, ring.
Let $J$ denote the (finite) set of irreducible constituents of~$M$. For each $j \in J$ there is a unique submodule $M_j$ whose cosocle
is $j$. There is a partial order $\le$ on $J$ such that $i \le j$ if and only if $M_j$ contains $i$ as constituent.  Then the
lattice of submodules of $M$ (with respect to inclusion) is naturally isomorphic to the lattice of lower sets in $(J,\le)$
(with respect to inclusion).  (By a \emph{lower set} in a partially ordered set $(J,\le)$ we mean a subset $J' \subset J$
such that $i \in J'$ whenever $i \le j$ and $j \in J'$.)  A submodule $N$ is sent to the lower set $\JH(N)$ of
constituents of $N$, whereas a lower set $J' \subset J$ is sent to the submodule $\sum_{j \in J'} M_j$.

Recall from Cor.~\ref{cor:steinb-irred-constit} that all constituents of $\obInd 1$ occur with multiplicity one.

\begin{prop}\label{prop:submod-obind}
  Let $\Pi$ denote the set of standard parabolic subgroups.
  The lattice of submodules of $\obInd 1$ is naturally isomorphic to the lattice of lower sets in $(\Pi,\supset)$.
\end{prop}  

Note that we can identify the set of constituents of $\obInd 1$ with the set of standard parabolic subgroups via $\Sp P \leftrightarrow P$.
With this identification a submodule $N$ corresponds to the lower set $\JH(N)$.

We have an analogous result for $\opInd \Sp Q$ as in Prop.~\ref{prop:opind-spq}: the submodule lattice is isomorphic to the lattice
of lower sets in $(\Pi_Q,\supset)$, where $\Pi_Q$ consists of all standard parabolic subgroups $P'$ such that $P' \cap M = Q$. (Just
note that $\opInd \Sp Q$ is a subquotient of $\obInd 1$.)

\begin{proof}
  Let $P$ be a standard parabolic subgroup.  We will show that $\opInd 1$ is the unique submodule with cosocle $\Sp P$. This
  implies the result: we have $\JH(\opInd 1) = \{\Sp Q : Q \supset P\}$ by using~\eqref{eq:26} and induction.
  Recall from Section~\ref{sec:gen-steinberg} that each $\Sp Q$ contains a unique weight $V_Q$ and that these weights are pairwise
  distinct. Moreover we saw in the proof of Thm.~\ref{thm:steinb-irred} that $V_P$ is a weight of $\opInd 1$ and that it generates
  $\opInd 1$ as $G$-representation. Thus the image of $V_P$ generates any quotient of $\opInd 1$, so
  $\Sp P$ is the unique irreducible quotient of $\opInd 1$. 
\end{proof}

Now let $G = \GL_n$. Let $P = MN$ be a standard parabolic subgroup and $\sigma$ an irreducible admissible
$M$-representation. We will determine the submodule structure of $\opInd \sigma$.  Recall from Cor.~\ref{cor:supercuspidal}
that $\opInd \sigma$ has finite length and that the constituents occur with multiplicity one.

  We can use Thm.~\ref{thm:classification} to express $\pi := \opInd \sigma$ in the same form as in
  Thm.~\ref{thm:decomp-parab-ind-gln}. As in the proof of Thm.~\ref{thm:decomp-parab-ind-gln} we can further rewrite $\pi$ as
  $\Ind_{\o Q}^G (\tau_1 \otimes \cdots \otimes \tau_s)$, where $Q$ is the standard parabolic with Levi $\prod_i
  \GL_{m_i}$, such that for all $i$ either
  \begin{itemize}
  \item $\tau_i$ is supersingular and $m_i > 1$, or
  \item $\tau_i \cong \Ind_{\o R_i}^{\GL_{m_i}} \Sp {S_i} \otimes \eta'_i$ for some $\eta'_i : F\s \to \k\s$
  \end{itemize}
  and such that $\eta'_i \ne \eta'_{i+1}$ whenever both $\tau_i$ and $\tau_{i+1}$ fall into the second case. Here $R_i = L_i
  N_i'$ is a standard parabolic of $\GL_{m_i}$ and $S_i$ a standard parabolic of $L_i$. Let $I \subset \{1,\dots,s\}$ be the
  subset of those $i$ such that $\tau_i$ falls into the second case. Let $X$ be the set consisting of tuples $(S_i')_{i\in I}$ such
  that each $S_i'$ is a standard parabolic subgroup of $\GL_{m_i}$ satisfying $S_i' \cap L_i = S_i$. We define a partial
  order $\le_X$ on $X$ by declaring that $(S_i')_{i\in I} \le_X (S_i'')_{i\in I}$ if and only if $S_i' \supset S_i''$ for all $i\in I$.

\begin{prop}\label{prop:submod-parabind}
  With the above notation, the lattice of submodules of $\opInd \sigma$ is naturally isomorphic to the lattice of
  all lower sets in $(X, \le_X)$.
\end{prop}

\begin{proof}
  Whenever we are given $\GL_{m_i}$-representations $\tau_i'$ for
  $i \in I$, let $\varsigma(\tau'_i : i \in I)$ denote the $\prod_i \GL_{m_i}$-representation $\bigotimes_{i=1}^s \tau_i'$, where
  $\tau_i' = \tau_i$ for $i \not\in I$. Also let $\pi(\tau'_i : i \in I)$ denote $\oqInd \big(\varsigma(\tau'_i : i \in I)\big)$.
  Then the constituents of $\pi$ are the $\pi(\Sp {S_i'}: i \in I)$ for $(S_i')_{i \in I} \in X$.
  
  Note that $\pi$ is a subquotient of $\pi' := \pi(\Ind_{\o {B_i}}^{\GL_{m_i}} \eta_i' : i \in I)$, where $B_i$ is the Borel
  in $\GL_{m_i}$. The corresponding partially ordered sets $(X,\le_X)$, $(X', \le_{X'})$ are compatible in the sense that
  $\le_X$ is the restriction of $\le_{X'}$ to~$X$. Without loss of generality we can therefore assume that $\pi = \pi'$.
  (I.e., $R_i = B_i$ and $S_i = L_i = T$.) Then $\pi$ has constituents $\pi(\Sp {R_i}\otimes \eta'_i: i \in
  I)$, where the $R_i$ run over all standard parabolic subgroups of $\GL_{m_i}$ (for $i \in I$). Fix now such $R_i$. We claim
  that the submodule $\pi(\Ind_{\o{R_i}}^{\GL_{m_i}} \eta'_i : i \in I)$ has irreducible cosocle $\pi(\Sp {R_i}\otimes
  \eta'_i: i \in I)$. To see this, let $\o {V_i}$ be any weight of $\tau_i$ for $i \not\in I$ and the unique weight of $\Sp
  {R_i}\otimes \eta'_i$ for $i \in I$. Then $\o V := \bigotimes \o {V_i}$ is a weight of
  $\varsigma(\Ind_{\o{R_i}}^{\GL_{m_i}} \eta'_i : i \in I)$ and generates it as $\prod_i \GL_{m_i}$-representation (by the
  proof of Thm.~\ref{thm:steinb-irred}). Let $V$ be the $\prod_i \GL_{m_i}$-regular weight for $G$ that corresponds to $\o V$
  under the bijection of Lemma~\ref{lm:lift-to-m-reg-wt}. By Thm.~\ref{thm:cptind-parabind} it follows that $V$ generates
  $\pi(\Ind_{\o{R_i}}^{\GL_{m_i}} \eta'_i : i \in I)$ as $G$-representation. (Just pick any set of Hecke eigenvalues $\chi$
  for the weight $\o V$.) But~\eqref{eq:6} and the fact that generalised Steinberg representations can be distinguished by
  their weights show that $\pi(\Ind_{\o{R_i}}^{\GL_{m_i}} \eta'_i : i \in I)$ indeed has irreducible cosocle $\pi(\Sp
  {R_i}\otimes \eta'_i: i \in I)$.

  Finally notice that the constituents of $\pi(\Ind_{\o{R_i}}^{\GL_{m_i}} \eta'_i : i \in I)$ are all
  $\pi(\Sp {R'_i}\otimes \eta'_i: i \in I)$ with $R_i' \supset R_i$ for all $i \in I$. (Decompose each $\Ind_{\o{R_i}}^{\GL_{m_i}} \eta_i'$
  and apply Thm.~\ref{thm:irred-parab-ind-gln-2}.)
\end{proof}

\bibliography{parab-ind}

\newcommand{\etalchar}[1]{$^{#1}$}
\def\cprime{$'$}
\begin{thebibliography}{Eme10b}

\bibitem[BL94]{bib:BL-general}
L.~Barthel and R.~Livn{\'e}.
\newblock Irreducible modular representations of {${\rm GL}\sb 2$} of a local
  field.
\newblock {\em Duke Math. J.}, 75(2):261--292, 1994.

\bibitem[BL95]{bib:BL-unram}
L.~Barthel and R.~Livn{\'e}.
\newblock Modular representations of {${\rm GL}\sb 2$} of a local field: the
  ordinary, unramified case.
\newblock {\em J. Number Theory}, 55(1):1--27, 1995.

\bibitem[BP]{bib:BP}
Christophe Breuil and Vytautas Pa{\v{s}}k{\=u}nas.
\newblock Towards a modulo {$p$} {L}anglands correspondence for {${\rm GL}_2$}.
\newblock To appear in Mem. Amer. Math. Soc.

\bibitem[Bre03a]{bib:Breuil}
Christophe Breuil.
\newblock Sur quelques repr\'esentations modulaires et {$p$}-adiques de {${\rm
  GL}\sb 2(\bold Q\sb p)$}. {I}.
\newblock {\em Compositio Math.}, 138(2):165--188, 2003.

\bibitem[Bre03b]{bib:Breuil2}
Christophe Breuil.
\newblock Sur quelques repr\'esentations modulaires et {$p$}-adiques de {${\rm
  GL}\sb 2(\bold Q\sb p)$}. {II}.
\newblock {\em J. Inst. Math. Jussieu}, 2(1):23--58, 2003.

\bibitem[BT72]{bib:BT1}
F.~Bruhat and J.~Tits.
\newblock Groupes r\'eductifs sur un corps local.
\newblock {\em Inst. Hautes \'Etudes Sci. Publ. Math.}, (41):5--251, 1972.

\bibitem[BT84]{bib:BT2}
F.~Bruhat and J.~Tits.
\newblock Groupes r\'eductifs sur un corps local. {II}. {S}ch\'emas en groupes.
  {E}xistence d'une donn\'ee radicielle valu\'ee.
\newblock {\em Inst. Hautes \'Etudes Sci. Publ. Math.}, (60):197--376, 1984.

\bibitem[BZ77]{bib:BZ}
I.~N. Bernstein and A.~V. Zelevinsky.
\newblock Induced representations of reductive {${\germ p}$}-adic groups. {I}.
\newblock {\em Ann. Sci. \'Ecole Norm. Sup. (4)}, 10(4):441--472, 1977.

\bibitem[Car79]{bib:Cartier}
P.~Cartier.
\newblock Representations of {$p$}-adic groups: a survey.
\newblock In {\em Automorphic forms, representations and {$L$}-functions
  ({P}roc. {S}ympos. {P}ure {M}ath., {O}regon {S}tate {U}niv., {C}orvallis,
  {O}re., 1977), {P}art 1}, Proc. Sympos. Pure Math., XXXIII, pages 111--155.
  Amer. Math. Soc., Providence, R.I., 1979.

\bibitem[Col10]{bib:Colmez}
Pierre Colmez.
\newblock Repr\'esentations de {${\rm GL}_2(\bold Q_p)$} et
  {$(\phi,\Gamma)$}-modules.
\newblock {\em Ast\'erisque}, (330):281--509, 2010.

\bibitem[Eme10a]{bib:emerton-ordinary1}
Matthew Emerton.
\newblock Ordinary parts of admissible representations of {$p$}-adic reductive
  groups {I}. {D}efinition and first properties.
\newblock {\em Ast\'erisque}, (331):355--402, 2010.

\bibitem[Eme10b]{bib:emerton-ordinary2}
Matthew Emerton.
\newblock Ordinary parts of admissible representations of {$p$}-adic reductive
  groups {II}. {D}erived functors.
\newblock {\em Ast\'erisque}, (331):403--459, 2010.

\bibitem[G{\etalchar{+}}70]{bib:SGA3}
Alexandre Grothendieck et~al.
\newblock {\em {SGA} 3: {S}ch\'emas en Groupes I, II, III}.
\newblock Lecture Notes in Math. 151, 152, 153. Springer-Verlag, Heidelberg,
  1970.

\bibitem[GK]{bib:grosse-kloenne}
Elmar Gro{\ss}e-Kl{\"o}nne.
\newblock On special representations of {$p$}-adic reductive groups.
\newblock Preprint, version of 9/14/2009.

\bibitem[Gro98]{bib:Gross_Satake}
Benedict~H. Gross.
\newblock On the {S}atake isomorphism.
\newblock In {\em Galois representations in arithmetic algebraic geometry
  (Durham, 1996)}, volume 254 of {\em London Math. Soc. Lecture Note Ser.},
  pages 223--237. Cambridge Univ. Press, Cambridge, 1998.

\bibitem[Her09]{bib:thesis}
Florian Herzig.
\newblock The weight in a {S}erre-type conjecture for tame {$n$}-dimensional
  {G}alois representations.
\newblock {\em Duke Math. J.}, 149(1):37--116, 2009.

\bibitem[Her11]{bib:satake}
Florian Herzig.
\newblock A {S}atake isomorphism in characteristic {$p$}.
\newblock {\em Compos. Math.}, 147(1):263--283, 2011.

\bibitem[HKP10]{bib:HKP}
Thomas~J. Haines, Robert~E. Kottwitz, and Amritanshu Prasad.
\newblock Iwahori-{H}ecke algebras.
\newblock {\em J. Ramanujan Math. Soc.}, 25(2):113--145, 2010.

\bibitem[Hu10]{bib:Hu}
Yongquan Hu.
\newblock Sur quelques repr\'esentations supersinguli\`eres de {${\rm
  GL}_2(\Bbb Q_{p^f})$}.
\newblock {\em J. Algebra}, 324(7):1577--1615, 2010.

\bibitem[Iwa66]{bib:Iwahori}
Nagayoshi Iwahori.
\newblock Generalized {T}its system ({B}ruhat decompostition) on {$p$}-adic
  semisimple groups.
\newblock In {\em Algebraic {G}roups and {D}iscontinuous {S}ubgroups ({P}roc.
  {S}ympos. {P}ure {M}ath., {B}oulder, {C}olo., 1965)}, pages 71--83. Amer.
  Math. Soc., Providence, R.I., 1966.

\bibitem[Jan03]{bib:Jan-reps}
Jens~Carsten Jantzen.
\newblock {\em Representations of algebraic groups}, volume 107 of {\em
  Mathematical Surveys and Monographs}.
\newblock American Mathematical Society, Providence, RI, second edition, 2003.

\bibitem[Kis09]{bib:Kisin}
Mark Kisin.
\newblock The {F}ontaine-{M}azur conjecture for {${\rm GL}_2$}.
\newblock {\em J. Amer. Math. Soc.}, 22(3):641--690, 2009.

\bibitem[Oll06]{bib:ollivier}
Rachel Ollivier.
\newblock Crit\`ere d'irr\'eductibilit\'e pour les s\'eries principales de
  {${\rm GL}\sb n(F)$} en caract\'eristique {$p$}.
\newblock {\em J. Algebra}, 304(1):39--72, 2006.

\bibitem[Pa{\v{s}}10]{bib:pask}
Vytautas Pa{\v{s}}k{\=u}nas.
\newblock Extensions for supersingular representations of {${\rm GL}_2(\Bbb
  Q_p)$}.
\newblock {\em Ast\'erisque}, (331):317--353, 2010.

\bibitem[SS91]{bib:schneider-stuhler}
P.~Schneider and U.~Stuhler.
\newblock The cohomology of {$p$}-adic symmetric spaces.
\newblock {\em Invent. Math.}, 105(1):47--122, 1991.

\bibitem[Tit79]{bib:Tits}
J.~Tits.
\newblock Reductive groups over local fields.
\newblock In {\em Automorphic forms, representations and {$L$}-functions
  ({P}roc. {S}ympos. {P}ure {M}ath., {O}regon {S}tate {U}niv., {C}orvallis,
  {O}re., 1977), {P}art 1}, Proc. Sympos. Pure Math., XXXIII, pages 29--69.
  Amer. Math. Soc., Providence, R.I., 1979.

\bibitem[Vig04]{bib:vigneras2}
Marie-France Vign{\'e}ras.
\newblock Representations modulo {$p$} of the {$p$}-adic group {${\rm
  GL}(2,F)$}.
\newblock {\em Compos. Math.}, 140(2):333--358, 2004.

\bibitem[Vig05]{bib:Vig-ss}
Marie-France Vign{\'e}ras.
\newblock Pro-{$p$}-{I}wahori {H}ecke ring and supersingular {$\overline{\bold
  F}_p$}-representations.
\newblock {\em Math. Ann.}, 331(3):523--556, 2005.

\bibitem[Vig08]{bib:vigneras}
Marie-France Vign{\'e}ras.
\newblock S\'erie principale modulo {$p$} de groupes r\'eductifs {$p$}-adiques.
\newblock {\em Geom. Funct. Anal.}, 17(6):2090--2112, 2008.

\bibitem[Zel80]{bib:BZ2}
A.~V. Zelevinsky.
\newblock Induced representations of reductive {${\germ p}$}-adic groups. {II}.
  {O}n irreducible representations of {${\rm GL}(n)$}.
\newblock {\em Ann. Sci. \'Ecole Norm. Sup. (4)}, 13(2):165--210, 1980.

\bibitem[Zel81]{bib:zele-kl}
A.~V. Zelevinski{\u\i}.
\newblock The {$p$}-adic analogue of the {K}azhdan-{L}usztig conjecture.
\newblock {\em Funktsional. Anal. i Prilozhen.}, 15(2):9--21, 96, 1981.

\end{thebibliography}
\bibliographystyle{halpha}

\end{document}